\documentclass[12pt]{article}
\usepackage{amsmath}
	\allowdisplaybreaks
\usepackage{amsthm,amssymb,mathrsfs}
\usepackage{stmaryrd}
\usepackage{dsfont}
\usepackage{mathtools}
\usepackage[shortlabels]{enumitem}
	\setlist[enumerate]{resume}
	\setlist[enumerate,itemize]{topsep=0em,itemsep=.5em,partopsep=0em,parsep=0em}	
\usepackage{graphics}
\usepackage[margin=3cm]{geometry}
\usepackage{changepage}
\usepackage{multicol}

\usepackage[colorlinks=true,
    bookmarksnumbered=true,
	 linkcolor = blue, citecolor = magenta, urlcolor = magenta]{hyperref}
\usepackage{bookmark}

\usepackage{tikz}

\usepackage{tikz-cd}
	\tikzcdset{every label/.append style = {font = \small}}
	
\usepackage{stmaryrd}

\renewcommand{\thesubsection}{} 

\usepackage{titlesec}
	\titleformat{\section}{\bf\filcenter}{\thesection.}{.5em}{}
	\titleformat{\subsection}{\bf}{\thesubsection}{0em}{}
	\titlespacing*{\section}{0pt}{1em}{.5em}
	\titlespacing*{\subsection}{0pt}{.5em}{.25em}

\usepackage{tocloft}

	\cftsetindents{subsection}{1.5em}{1em}

\theoremstyle{definition}
\newtheorem{definition}{Definition}[section]
\newtheorem{context}[definition]{\indent}

\numberwithin{equation}{definition}

\title{\large\bf On endomorphism algebras of Gelfand-Graev representations}
\author{\normalsize Tzu-Jan Li\footnote{Institut de Mathématiques de Jussieu-Paris Rive Gauche (IMJ-PRG), 4 place Jussieu, 75252 Paris cedex 05, FRANCE. Email: {\tt tzu-jan.li@imj-prg.fr}}}
\date{\normalsize June 17, 2021}

\begin{document}
\maketitle
$\left.\right.$\\[-18mm]
\begin{abstract}
For a connected reductive group $G$ defined over $\mathbb{F}_q$ and equipped with the induced Frobenius endomorphism $F$, we study the relation among the following three $\mathbb{Z}$-algebras: (i) the $\mathbb{Z}$-model $\mathsf{E}_G$ of endomorphism algebras of Gelfand-Graev representations of $G^F$; (ii) the Grothendieck group $\mathsf{K}_{G^\ast}$ of the category of representations of $G^{\ast F^\ast}$ over $\overline{\mathbb{F}_q}$ (Deligne-Lusztig dual side); (iii) the ring $\mathsf{B}_{G^\vee}$ of the scheme $(T^\vee\sslash W)^{F^\vee}$ over $\mathbb{Z}$ (Langlands dual side). The comparison between (i) and (iii) is motivated by recent advances in the local Langlands program.
\end{abstract}

\noindent{\bf Table of contents}\\[-12mm]
\tableofcontents

{$\left.\right.$}\\[-2mm]
\begin{minipage}{0.10\linewidth}   
\includegraphics[scale=0.06]{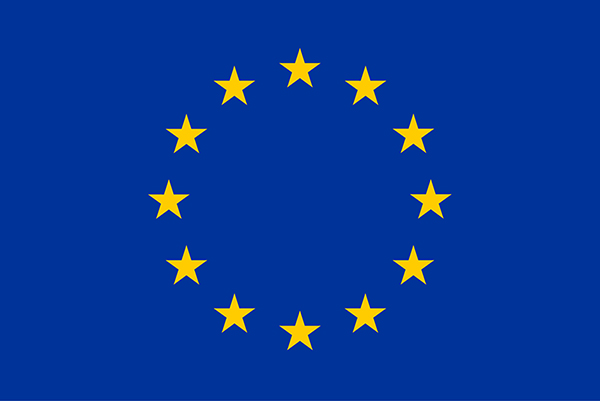}
\end{minipage} 
\begin{minipage}{0.90\linewidth}
{\footnotesize\it This project has received funding from the
European Union’s Horizon 2020 research and \\[-1mm] innovation programme under the Marie Skłodowska-Curie grant agreement No 754362.}   
\end{minipage}

{\scriptsize$\left.\right.$}

\newpage

\section{Introduction and main results}

\begin{context}\label{the-problem}
{\it The problem.}\;--- In this article, we would like to study the relation among three algebras $\mathsf{E}_G$, $\mathsf{K}_{G^\ast}$ and $\mathsf{B}_{G^\vee}$ coming from the group representation theory and from the invariant theory. Let us now describe these algebras.

So let $p$ be a prime number, let $q=p^r$ for some $r\in\mathbb{N}^\ast$, and fix an algebraic closure $\overline{\mathbb{F}_q}$ of the finite field $\mathbb{F}_q$. Let $G$ be a connected reductive group defined over $\mathbb{F}_q$, and let $F:G\longrightarrow G$ be the Frobenius endomorphism associated to the $\mathbb{F}_q$-structure of $G$. We write $G(\overline{\mathbb{F}_q})$ simply as $G$, so $G^F=\{x\in G:F(x)=x\}$ is a finite group. Fix an $F$-stable maximal torus $T$, and let $W=N_G(T)/T$ be the Weyl group of $(G,T)$; then $W$ is $F$-stable. Let $(G^\ast,T^\ast,F^\ast)$ be a Deligne-Lusztig dual (defined over $\overline{\mathbb{F}_q}$; see \S\,\ref{setup-global}) of $(G,T,F)$, and let $(G^\vee,T^\vee)$ be a Langlands dual (defined and split over $\mathbb{Z}$; see \S\,\ref{Langlands-dual}) of $(G,T)$. The endomorphism $F$ induces an endomorphism $F^\vee$ on the character group $X(T^\vee)=\mathrm{Hom}_{\mathrm{alg}}(T^\vee,\mathbb{G}_m)$ (\S\,\ref{tau-F}). In addition, let $\overline{\mathbb{Q}}$ be the field of algebraic numbers, and let $\overline{\mathbb{Z}}$ be the ring of algebraic integers. With this setup:

\begin{itemize}
\item $\mathsf{E}_G$ is the $\mathbb{Z}$-model of endomorphism algebras of Gelfand-Graev representations of $G^F$ (see \S\S\,\ref{E-definition}\,-\,\ref{E-Zalgebra}).
\item $\mathsf{K}_{G^\ast}$ is the Grothendieck group of the category of finite-dimensional representations of $G^{\ast F^\ast}$ over $\overline{\mathbb{F}_q}$ (\S\,\ref{K-definition}); it is a $\mathbb{Z}$-algebra whose multiplication comes from the tensor product.
\item $\mathsf{B}_{G^\vee}$ is the ring of functions of the affine $\mathbb{Z}$-scheme $(T^\vee\sslash W)^{F^\vee}$ (\S\,\ref{B-definition}); more precisely, $\mathsf{B}_{G^\vee}=\mathbb{Z}[X(T^\vee)]^W/I$ where $I$ is the ideal of $\mathbb{Z}[X(T^\vee)]^W$ generated by the set $\{F^\vee f-f:f\in\mathbb{Z}[X(T^\vee)]^W\}$.
\end{itemize}
Our goal is to compare these three algebras over rational and integral coefficients.

\end{context}

\begin{context}\label{the-analysis}
{\it First analysis.}\;--- In the case of $\overline{\mathbb{Q}}$-coefficients, Curtis-Deligne-Lusztig's theory (\S\S\,\ref{Gelfand-Graev}\,-\,\ref{E-decomposition}) and Brauer theory (\S\,\ref{Brauer-character}) give us two $\overline{\mathbb{Q}}$-algebra isomorphisms:
\begin{equation}\label{analysis-eq0}
\overline{\mathbb{Q}}\mathsf{E}_G\simeq \overline{\mathbb{Q}}^{(G^{\ast F^\ast})_{\mathrm{ss}}/\sim} \quad \mbox{and} \quad \overline{\mathbb{Q}}\mathsf{K}_{G^\ast}\simeq \overline{\mathbb{Q}}^{(G^{\ast F^\ast})_{p'}/\sim}.
\end{equation}
In (\ref{analysis-eq0}), the first isomorphism only depends on a choice of identifications $(\mathbb{Q}/\mathbb{Z})_{p'}\simeq\overline{\mathbb{F}_q}^\times\hookrightarrow\overline{\mathbb{Q}}^\times$, and the second only depends on a choice of embedding $\overline{\mathbb{F}_q}^\times\hookrightarrow\overline{\mathbb{Q}}^\times$. Via (\ref{analysis-eq0}), the natural maps of finite sets
\begin{equation}\label{analysis-eq1}
\begin{tikzcd}
((G^{\ast F^\ast})_{p'}/\!\sim)\arrow[equal,r] & ((G^{\ast F^\ast})_{\mathrm{ss}}/\!\sim) \arrow[twoheadrightarrow,rr,"\S\,\ref{counting-lemma}","\S\,\ref{counting-derived}"'] &&(T^\vee\sslash W)^{F^\vee}(\overline{\mathbb{Q}})
\end{tikzcd}
\end{equation} 
then induce the following $\overline{\mathbb{Q}}$-algebra homomorphisms:
\begin{equation}\label{analysis-eq2}
\begin{tikzcd}
\overline{\mathbb{Q}}\mathsf{K}_{G^\ast} \arrow[leftrightarrow,r,"\sim"]&\overline{\mathbb{Q}}\mathsf{E}_G\arrow[hookleftarrow,r]&(\overline{\mathbb{Q}}\mathsf{B}_{G^\vee})_{\mathrm{red}}\arrow[equal,r,"\S\,\ref{B-decomposition}"]&\overline{\mathbb{Q}}\mathsf{B}_{G^\vee,\,\mathrm{red}}\arrow[twoheadleftarrow,r]  &  \overline{\mathbb{Q}}\mathsf{B}_{G^\vee},
\end{tikzcd}
\end{equation}
where for a ring $A$ we denote by $A_{\mathrm{red}}$ the associated reduced ring. 

It is then natural to ask: when are the maps in (\ref{analysis-eq2}) all isomorphisms? could we change the $\overline{\mathbb{Q}}$-coefficients in (\ref{analysis-eq2}) by $\Lambda$-coefficients with $\Lambda$ a subring of $\overline{\mathbb{Q}}$, such that we still have analogous $\Lambda$-algebra homomorphisms? 

In \cite{Bonnafe-Kessar} (see also \S\S\,\ref{symmetrizing-lemma}\,-\,\ref{Curtis-saturation}), Bonnaf{\'e} and Kessar have proved a ``saturatedness property" for the Curtis homomorphism, where they used a lemma on symmetrizing forms to descend an equality of algebras over rational coefficients into an equality of algebras over integral coefficients. This idea will be the start point of our present work.

\end{context}

\begin{context}\label{the-results}
{\sc Main results.}\;---
{\it Keep the notation and assumptions introduced so far.
\begin{enumerate}[\normalfont (a)]
\item \mbox{\normalfont [see \S\S\,\ref{K-main-theorem}\,-\,\ref{K-main-corollary}]} The $\overline{\mathbb{Q}}$-algebra isomorphism $\overline{\mathbb{Q}}\mathsf{E}_G\simeq\overline{\mathbb{Q}}\mathsf{K}_{G^\ast}$ in (\ref{analysis-eq2}) descends to a ${\mathbb{Z}}[\frac{1}{p|W|}]$-algebra isomorphism ${\mathbb{Z}}[\frac{1}{p|W|}]\mathsf{E}_G\simeq{\mathbb{Z}}[\frac{1}{p|W|}]\mathsf{K}_{G^\ast}$. 

\item \mbox{\normalfont [see \S\,\ref{B-general} and \S\,\ref{q-restricted}]} $\mathsf{B}_{G^\vee,\,\mathrm{red}}$ is a free $\mathbb{Z}$-module of rank 
\begin{equation}\label{main-eq1}
\mathrm{rank}_{\mathbb{Z}}\,\mathsf{B}_{G^\vee,\,\mathrm{red}}=|(T^\vee\sslash W)^{F^\vee}(\overline{\mathbb{Q}})|=|(G^\ast_{\mathrm{ss}}/\!\sim)^{F^\ast}|;
\end{equation}
moreover, if $G_{\mathrm{der}}^\ast$ (or equivalently $G_{\mathrm{der}}^\vee$) is simply-connected, then $\mathsf{B}_{G^\vee}$ is a reduced ring (so $\mathsf{B}_{G^\vee}=\mathsf{B}_{G^\vee,\,\mathrm{red}}$), and the above rank (\ref{main-eq1}) is also equal to $|(G^{\ast F^\ast})_{\mathrm{ss}}/\!\sim|$.
\item \mbox{\normalfont [see \S\S\,\ref{KB-proposition}\,-\,\ref{KB-isomorphism}]} The formal character isomorphism 
\[
\mathrm{ch}:\mathsf{K}(\mathrm{Rep}_{\mathrm{alg}}(G^\ast))\xrightarrow{\;\;\sim\;\;}\mathbb{Z}[X(T^\vee)]^W
\] induces an inclusion of rings $\mathsf{B}_{G^\vee,\,\mathrm{red}}\hookrightarrow \mathsf{K}_{G^\ast}$ which is compatible with the identifications in (\ref{analysis-eq2}); moreover, if $G_{\mathrm{der}}^\ast$ is simply-connected, then this inclusion of rings becomes a ring isomorphism: $\mathsf{B}_{G^\vee}=\mathsf{B}_{G^\vee,\,\mathrm{red}}\simeq\mathsf{K}_{G^\ast}$.
\end{enumerate}
We expect that (a) is true over $\mathbb{Z}[\frac{1}{p}]$ and that the ring $\mathsf{B}_{G^\vee}$ in (b) is always reduced; at the moment, they are achieved for the following special cases:
\begin{enumerate}[\normalfont(a),resume]
\item \mbox{\normalfont [see \S\,\ref{GL2-exceptional} and \S\,\ref{PGL2-exceptional}]} If $G=\mathrm{GL}_2(\overline{\mathbb{F}_q})$ or $\mathrm{PGL}_2(\overline{\mathbb{F}_q})$, the isomorphism ${\mathbb{Z}}[\frac{1}{p|W|}]\mathsf{E}_G\simeq{\mathbb{Z}}[\frac{1}{p|W|}]\mathsf{K}_{G^\ast}$ in (a) descends to a ${\mathbb{Z}}[\frac{1}{p}]$-algebra isomorphism ${\mathbb{Z}}[\frac{1}{p}]\mathsf{E}_G\simeq{\mathbb{Z}}[\frac{1}{p}]\mathsf{K}_{G^\ast}$.
\item \mbox{\normalfont [see \S\,\ref{SO2n-basis}]} When $G=\mathrm{SO}_{2n}(\overline{\mathbb{F}_q})$ (so $G_{\mathrm{der}}^\ast=\mathrm{SO}_{2n}(\overline{\mathbb{F}_q})$ and $G_{\mathrm{der}}^\vee=\mathrm{SO}_{2n}$ are not simply-connected), $\mathsf{B}_{G^\vee}$ is still a reduced ring (so $\mathsf{B}_{G^\vee}=\mathsf{B}_{G^\vee,\,\mathrm{red}}$) and is a free $\mathbb{Z}$-module of rank as described in (\ref{main-eq1}); explicitly, $\mathrm{rank}_{\mathbb{Z}}\,\mathsf{B}_{G^\vee}=2q^n-2q^{n-1}+q^{n-2}$.
\end{enumerate}
}

\end{context}

\begin{context}\label{the-motivation}
{\it Motivation.}\;--- The problem of comparison between $\mathsf{E}_G$ and $\mathsf{B}_{G^\vee}$ is a finite group analogue of the ``local Langlands correspondence in families" (LLIF) conjecture in \cite{DHKM}, which asserts the existence of an isomorphism from the ring of functions of the moduli stack of Langlands parameters of a $p$-adic reductive group to the endomorphism ring of its Whittaker space. Note that the LLIF conjecture is implied by the Fargues-Scholze conjecture in \cite[Conj.\;I.10.2]{Fargues-Scholze}.

In the case of a reductive group $\mathbf{G}$ defined and split over a $p$-adic field $F$ whose residue field is $\mathbb{F}_q$, let $O_F$ be the ring of integers of $F$ and write $G=\mathbf{G}(\overline{\mathbb{F}_q})$; then a special case of the LLIF conjecture predicts a ring isomorphism 
\[
\mathcal{O}_{\mathcal{Z}^1_{\mathrm{tame}}\sslash \mathbf{G}^\vee}\simeq\mathrm{End}(\mbox{c-Ind}\,_{\mathbf{G}(O_F)}^{\mathbf{G}(F)}\Gamma_{G})
\]
which is compatible with the classical local Langlands conjectures; in the above, 
\[
\mathcal{Z}^1_{\mathrm{tame}}=\{(s,\mathcal{F})\in\mathbf{G}^\vee\times\mathbf{G}^\vee:\mathcal{F}s\mathcal{F}^{-1}=s^q\}
\]
is the moduli space of tame Langlands parameters for $\mathbf{G}^\vee$, $\mathcal{O}_{\mathcal{Z}^1_{\mathrm{tame}}\sslash \mathbf{G}^\vee}$ is the ring of functions of $\mathcal{Z}^1_{\mathrm{tame}}\sslash \mathbf{G}^\vee$, $\Gamma_G$ is a Gelfand-Graev representation of $G$, and $\mbox{c-Ind}\,_{\mathbf{G}(O_F)}^{\mathbf{G}(F)}\Gamma_{G}$ is (isomorphic to) the depth-zero part of the corresponding Whittaker space. The first projection map $\mathbf{G}^\vee\times\mathbf{G}^\vee\longrightarrow \mathbf{G}^\vee$ induces a morphism $\mathcal{Z}^1_{\mathrm{tame}}\sslash \mathbf{G}^\vee\longrightarrow (\mathbf{G}^\vee\sslash \mathbf{G}^\vee)^{(\cdot)^q}\simeq(T^\vee\sslash W)^{(\cdot)^q}$ and hence a ring homomorphism $\mathsf{B}_{G^\vee}\longrightarrow \mathcal{O}_{\mathcal{Z}^1_{\mathrm{tame}}\sslash \mathbf{G}^\vee}$. We may thus draw the following diagram:
\[
\begin{tikzcd}[column sep = large]
\mathcal{O}_{\mathcal{Z}^1_{\mathrm{tame}}\sslash \mathbf{G}^\vee}\arrow[leftrightarrow,r,"\mathrm{LLIF}","\sim"']&\mathrm{End}(\mbox{c-Ind}\,_{\mathbf{G}(O_F)}^{\mathbf{G}(F)}\Gamma_{G})\\
\mathsf{B}_{G^\vee}\arrow[u] \arrow[dashed,no head,r,"?" description]&\mathsf{E}_G\arrow[hookrightarrow,u]
\end{tikzcd}
\]

If the center of $G$ is connected, for an integral extension $\Lambda\supset\overline{\mathbb{Z}}[\frac{1}{p}]$, the ring $\Lambda\mathsf{B}_{G^\vee}$ (resp.\;$\Lambda\mathsf{E}_G$) should be thought as the integral closure of the scalars $\Lambda$ in $\mathcal{O}_{\mathcal{Z}^1_{\mathrm{tame}}\sslash \mathbf{G}^\vee}$ (resp.\;in $\mathrm{End}(\mbox{c-Ind}\,_{\mathbf{G}(O_F)}^{\mathbf{G}(F)}\Gamma_{G})$), so the LLIF conjecture should give us a $\Lambda$-algebra isomorphism 
\begin{equation}\label{conjecture-eq1}
\Lambda\mathsf{B}_{G^\vee}\simeq\Lambda\mathsf{E}_G.
\end{equation}

In \cite{Helm} and \cite{Helm-Moss}, the LLIF has been proved for $\mathbf{G}=\mathrm{GL}_n(F)$, and, as a corollary, a $\Lambda$-algebra isomorphism (\ref{conjecture-eq1}) has also been deduced for $G=\mathrm{GL}_n(\overline{\mathbb{F}_q})$ and for $\Lambda$ being the ring of Witt vectors of $\overline{\mathbb{F}_\ell}$. However, as (\ref{conjecture-eq1}) is basically a result of finite groups, an argument without $p$-adic techniques (such as those used for the LLIF) is expected, and the present article provides such an argument under some additional hypotheses:

{\sc Theorem.}\;(\S\,\ref{the-results}(a)(c), or \S\,\ref{EB-isomorphism}) --- {\it If $G_{\mathrm{der}}^\ast$ is simply-connected and if $\Lambda$ is an integral domain in which both $p$ and $|W|$ are invertible, then $\Lambda\mathsf{B}_{G^\vee}\simeq\Lambda\mathsf{E}_G$ as $\Lambda$-algebras.}

\end{context}

\begin{context}\label{the-plan}
{\it Plan of the article.}\;--- In Section \ref{E-section}, we review properties of endomorphism algebras of a modular Gelfand-Graev representation, construct the $\mathbb{Z}$-model $\mathsf{E}_G$, and study Curtis homomorphisms; in Section \ref{K-section}, we introduce the Grothendieck group $\mathsf{K}_{G^\ast}$ on the Deligne-Lusztig dual side, review the Brauer theory, and then compare $\mathsf{K}_{G^\ast}$ with the algebra $\mathsf{E}_G$. In Section \ref{B-section}, using methods from the algebraic geometry, the combinatorics of root datum and the algebraic representation theory, we turn to study the algebra $\mathsf{B}_{G^\vee}$ on the Langlands-dual side and then compare $\mathsf{B}_{G^\vee}$ with $\mathsf{K}_{G^\ast}$ and $\mathsf{E}_G$. In the final Section \ref{Example-section}, we use examples to discuss possible improvements of the results obtained in the previous three sections.
\end{context}

\begin{context}\label{graphical-summary}
{\it A graphical summary.}\;--- Let $\Lambda$ be an integral domain such that $\mathbb{Z}[\frac{1}{p}]\subset\Lambda\subset\overline{\mathbb{Q}}$. The following diagram, which will eventually be shown to be commutative (see \S\,\ref{Curtis-geometry} and \S\,\ref{KB-proposition}), summarizes the relations among the principal objects studied in this article. 

Notations used in the diagram: $(G^{\ast F^\ast})_{\mathrm{ss}}/\!\sim$ (resp.\;$(G^{\ast F^\ast})_{p'}/\!\sim$) denotes the set of semisimple (resp.\;$p$-regular) conjugacy classes in $G^{\ast F^\ast}$; ``$\mathrm{Res}$" means the obvious restriction maps; $T_w:={}^g T=gTg^{-1}$ is the $F$-stable maximal torus of $G$ whose $G^F$-conjugacy class corresponds to $w\in W$ (so $g\in G$ and the image of $g^{-1}F(g)\in N_G(T)$ in $W$ is $w$); recall that the $G^F$-conjugacy classes of $F$-stable tori in $G$ are parametrized by elements of $W$ in this way.
\[\small
\begin{tikzcd}[column sep = 0.2em, row sep = 0.2em]
\prod\limits_{w\in W} {\overline{\mathbb{Q}}}^{(T_w)^{\ast F^\ast}} \arrow[hookleftarrow,rr,"\mathrm{Res}"] &  & {\overline{\mathbb{Q}}}^{(G^{\ast F^\ast})_{\mathrm{ss}}/\sim} \arrow[hookleftarrow,rr] \arrow[leftrightarrow,rd,"\sim" sloped] &  & {\overline{\mathbb{Q}}}^{(T^\vee\sslash W)^{F^\vee}(\overline{\mathbb{Q}})}  &  \\
  &  &  & {\overline{\mathbb{Q}}}^{(G^{\ast F^\ast})_{p'}/\sim} \arrow[hookleftarrow,ru]  &  &  \\
\prod\limits_{w\in W}\overline{\mathbb{Q}}(T_w)^F \arrow[leftrightarrow,rd,"\sim" sloped, shorten <= -1em, near start] \arrow[leftrightarrow,uu,"\sim" sloped] \arrow[hookleftarrow,rr,"\mathrm{Cur}^G"] &  & \overline{\mathbb{Q}}\mathsf{E}_G  \arrow[hookleftarrow,rr] \arrow[leftrightarrow,rd,"\sim" sloped] \arrow[leftrightarrow,uu,"\sim" sloped] &  & \overline{\mathbb{Q}}\mathsf{B}_{G^\vee,\,\mathrm{red}} \arrow[leftrightarrow,uu,"\sim" sloped] &  \overline{\mathbb{Q}}\mathsf{B}_{G^\vee} \arrow[twoheadrightarrow,l] \\
  & \prod\limits_{w\in W}\overline{\mathbb{Q}}\mathsf{K}_{(T_w)^\ast} &  & \overline{\mathbb{Q}}\mathsf{K}_{G^\ast} \arrow[leftrightarrow, uu, crossing over,"\sim" sloped, near start] \arrow[hookleftarrow,ru] &  &  \\
\prod\limits_{w\in W}\Lambda(T_w)^F \arrow[hookrightarrow,uu] \arrow[hookleftarrow,rr,"\mathrm{Cur}^G"', near end] \arrow[leftrightarrow,rd,"\sim" sloped, shorten <= -1em, near start] &  & \Lambda\mathsf{E}_G \arrow[hookrightarrow,uu] \arrow[dashed,no head,rr,"?" description, near start] \arrow[dashed,no head,rd, "?" description] &  & \Lambda\mathsf{B}_{G^\vee,\,\mathrm{red}} \arrow[hookrightarrow,uu] & \Lambda\mathsf{B}_{G^\vee} \arrow[twoheadrightarrow,l] \arrow[uu] \\
  & \prod\limits_{w\in W}\Lambda\mathsf{K}_{(T_w)^\ast} \arrow[hookleftarrow,rr,"\mathrm{Res}"'] \arrow[hookrightarrow,uu,crossing over] &  & \Lambda\mathsf{K}_{G^\ast} \arrow[hookleftarrow,ru] \arrow[hookrightarrow,uu,crossing over] &  &  \\
\prod\limits_{w\in W}\mathbb{Z}(T_w)^F \arrow[hookrightarrow,uu] \arrow[leftrightarrow,rd,"\sim" sloped, shorten <= -1em, near start] &  &  &  & \mathsf{B}_{G^\vee,\,\mathrm{red}} \arrow[hookrightarrow,uu] & \mathsf{B}_{G^\vee} \arrow[twoheadrightarrow,l] \arrow[uu] \\
  & \prod\limits_{w\in W}\mathsf{K}_{(T_w)^\ast} \arrow[hookleftarrow,rr,"\mathrm{Res}"'] \arrow[hookrightarrow,uu] &  & \mathsf{K}_{G^\ast} \arrow[hookrightarrow,uu] \arrow[hookleftarrow,ru] &  &  \\
  &  &  &  & \mathbb{Z}[X(T^\vee)]^W \arrow[equal,r] \arrow[twoheadrightarrow,uu] & \mathbb{Z}[X(T^\vee)]^W \arrow[twoheadrightarrow,uu] \\
  &  &  & \mathsf{K}(\mathrm{Rep}_{\mathrm{alg}}(G^\ast)) \arrow[ru,"\sim" sloped,"\mathrm{ch}"', shorten >=-0.3em] \arrow[uu,"\mathrm{Res}"] &  &    
\arrow[hookleftarrow,from=4-2,to=4-4,"\mathrm{Res}"', near end, crossing over] 
\end{tikzcd}
\]
\end{context}

\begin{context}\label{setup-global}\label{E-setup}\label{DL-dual}
{\it Notation and convention.}\;--- The following notation and assumptions, as well as those introduced so far, will be used throughout the article. 

{\it Root data.}\;\cite{Springer} --- Fix an $F$-stable Borel subgroup $B$ of $G$ containing our $F$-stable maximal torus $T$. Let $X(T)=\mathrm{Hom}_{\mathrm{alg}}(T,\mathbb{G}_m)$ (resp.\;$Y(T)=\mathrm{Hom}_{\mathrm{alg}}(\mathbb{G}_m,T)$) be the character group (resp.\;cocharacter group) of $T$. Denote by $R\subset X(T)$ and $R^\vee\subset Y(T)$ (resp.\;$\Delta\subset X(T)$ and $\Delta^\vee\subset Y(T)$) the set of roots and the set of coroots (resp.\;the set of simple roots and the set of simple coroots) determined by $(G,T)$ (resp.\;by $(G,T,B)$). Then $(X(T),Y(T),R,R^\vee)$ the root datum of $(G,T)$. Recall that $(G,T)$ is classified (up to isomorphism) by its root datum.

{\it Deligne-Lusztig dual.}\;\cite{Digne-Michel} --- Fix a Deligne-Lusztig dual $(G^\ast,T^\ast,F^\ast)$ of $(G,T,F)$ (all choices are isomorphic), in the sense that $(G^\ast,T^\ast,F^\ast)$ is
defined over $\mathbb{F}_q$ and is obtained by assigning its character group (resp.\;cocharacter group, resp.\;set of roots, resp.\;set of coroots) as $Y(T)$ (resp.\;$X(T)$, resp.\;$R^\vee$, resp.\;$R$). In particular, we have the identifications $X(T^\ast)=Y(T)$ and $Y(T^\ast)=X(T)$, both of which are compatible with the Frobenius actions. Again, we write $G^\ast(\overline{\mathbb{F}_q})$ simply as $G^\ast$. 

Let $(G^{\ast F^\ast})_{\mathrm{ss}}$ be the set of semisimple elements of $G^{\ast F^\ast}$, and let $(G^{\ast F^\ast})_{\mathrm{ss}}/\!\sim$ be the set of $G^{\ast F^\ast}$-conjugacy classes in $(G^{\ast F^\ast})_{\mathrm{ss}}$; for each $x\in (G^{\ast F^\ast})_{\mathrm{ss}}$, we shall denote by $[x]\in (G^{\ast F^\ast})_{\mathrm{ss}}/\!\sim$ its $G^{\ast F^\ast}$-conjugacy class.

{\it Finite group representation theory.}\;\cite{Serre}\cite{Deligne-Lusztig}\cite{Digne-Michel} --- Let $\Lambda$ be a commutative ring and let $H$ be a finite group. The group ring $\Lambda[H]$ will often be written as $\Lambda H$; an element $f:H\longrightarrow\Lambda$ of $\Lambda H$ is identified with the formal sum $\sum\limits_{h\in H}f(h)h$. Let $\mathrm{Rep}_{\Lambda}(H)$ be the category of finite-dimensional representations of $H$ over $\Lambda$, and let $\mathrm{Irr}_{\Lambda}(H)$ be the set of isomorphism classes of simple objects in $\mathrm{Rep}_{\Lambda}(H)$. For $V_1,V_2\in \mathrm{Rep}_{\overline{\mathbb{Q}}}(H)$ whose characters are $\chi_1,\chi_2$ respectively, their canonical pairings $\langle V_1,V_2\rangle_{H}:=\dim_{\overline{\mathbb{Q}}}\mathrm{Hom}_{\overline{\mathbb{Q}}H}(V_1,V_2)$ and $\langle\chi_1,\chi_2\rangle_{H}:=|H|^{-1}\sum\limits_{h\in H}\chi_1(h^{-1})\chi_2(h)$ coincide. Ordinary induction functors like $\mathrm{Ind}_{U^F}^{G^F}$, ordinary restriction functors like $\mathrm{Res}_{T^{\ast F^\ast}}^{G^{\ast F^\ast}}$, and Deligne-Lusztig induction functors like $R_T^G=R_{T\subset B}^G$, are defined in the usual way. 

{\it Dualities of tori.}\;\cite{Digne-Michel} --- We fix a group isomorphism $\iota:\overline{\mathbb{F}_q}^\times\xrightarrow{\;\sim\;}(\mathbb{Q}/\mathbb{Z})_{p'}$ (choices of roots of unity) and an injective group homomorphism $\jmath:(\mathbb{Q}/\mathbb{Z})_{p'}\hookrightarrow\overline{\mathbb{Q}}^\times$. Then $\kappa:=\jmath\circ\iota:\overline{\mathbb{F}_q}^\times\hookrightarrow \overline{\mathbb{Q}}^\times$ enables us to identify $\overline{\mathbb{F}_q}^\times$ as a subgroup of $
\overline{\mathbb{Q}}^\times$. For an $F$-stable maximal torus $S$ of $G$, denote by $S^\ast$ the $F$-stable maximal torus of $G^\ast$ which is the dual of $S$. Then $\iota$ identifies $\mathrm{Irr}_{\overline{\mathbb{F}_q}}(S^F)\simeq S^{\ast F^\ast}$, and $(\iota,\kappa)$ identifies $\mathrm{Irr}_{\overline{\mathbb{Q}}}(S^F)\simeq S^{\ast F^\ast}$.

\end{context}

{\bf Acknowledgements.}\;I would like to thank Professor Jean-Fran{\c c}ois Dat for introducing me this subject and for having many stimulating discussions with me.

\section{Endomorphism algebras of Gelfand-Graev representations}\label{E-section}

\subsection{The Gelfand-Graev representation $\Gamma_{G,\psi}$}

\begin{context}\label{nondegenerate}
{\it Regular linear characters.}\;\cite[Sec.\;2]{Digne-Lehrer-Michel} --- Denote by $U_\bullet$ the subgroup of $U$ generated by the root subgroups of non-simple roots (roots in $R-\Delta$); note that $U_\bullet=[U,U]$. Let $\Delta/F$ be the set of $F$-orbits in $\Delta$. For each $i\in \Delta/F$, let $U_i$ be the product of root subgroups of (simple) roots in $i$. The quotient group $U/U_\bullet$ is canonically isomorphic to the product group $\prod\limits_{i\in \Delta/F}U_i$, and this isomorphism is $F$-stable, so that $U^F\!/U_\bullet^F\simeq\prod\limits_{i\in \Delta/F}U_i^F$ as abelian groups. Then a linear character $\psi$ on $U^F$ is called regular if $\psi$ is trivial ($=1$) on $U_\bullet^F$ and is non-trivial on every $U_i^F$ ($i\in \Delta/F$).

Let $\Psi$ be the set of $\overline{\mathbb{Z}}^\times$-valued regular linear characters of $U^F$; the group $T^F$ acts on $\Psi$ by adjoint action (for $t\in T^F$ and $\psi\in\Psi$, ${}^t\psi:=\psi({}^{t^{-1}}(\cdot))=\psi(t^{-1}(\cdot)t)\in\Psi$). Let $Z$ be the center of $G$, so the adjoint group of $G$ is $G_{\mathrm{ad}}=G/Z$. Let $T_\mathrm{ad}=T/Z$ be the image of $T$ in $ G_{\mathrm{ad}}$; the group $T_{\mathrm{ad}}^F$ also acts on $\Psi$ by adjoint action (for $x\in T_{\mathrm{ad}}^F$ and $\psi\in\Psi$, choose a $t\in T$ such that $x=tZ$, and set ${}^x\psi:={}^t\psi$), and this $T_{\mathrm{ad}}^F$-action on $\Psi$ is regular ($=$ free and transitive). Hence the first Galois cohomology group $H^1(F,Z)$, identified with $T_{\mathrm{ad}}^F/T^F$, acts regularly on $\Psi/T^F$ (the set of $T^F$-orbits in $\Psi$) by adjoint action.
\end{context}

\begin{context}\label{Gelfand-Graev-definition}
{\it Gelfand-Graev representations.}\;\cite[Ch.\;14]{Digne-Michel} --- Let $\psi:U^F\longrightarrow\overline{\mathbb{Z}}^\times$ be a regular linear character (so $\psi\in\Psi$ in the notation of \S\,\ref{nondegenerate}), and consider the primitive central idempotent
\begin{equation}\label{idempotent-eq}
e_\psi:=\displaystyle\frac{1}{|U^F|}\sum_{u\in U^F}\psi(u^{-1})u\in\overline{\mathbb{Z}}[\textstyle{\frac{1}{p}}]U^F
\end{equation}
of $\psi$ ($|U^F|$ is a power of $p$). Let $(\overline{\mathbb{Z}}[\frac{1}{p}])_{\psi}$ be the $\overline{\mathbb{Z}}[\frac{1}{p}]U^F$-module $\overline{\mathbb{Z}}[\frac{1}{p}]$ on which $U^F$ acts by $\psi$; we have $(\overline{\mathbb{Z}}[\frac{1}{p}])_{\psi}\simeq \overline{\mathbb{Z}}[\frac{1}{p}]U^Fe_\psi=\overline{\mathbb{Z}}[\frac{1}{p}]e_\psi$ as $\overline{\mathbb{Z}}[\frac{1}{p}]U^F$-modules.
Then 
\begin{equation}\label{Gelfand-Graev-eq}
\Gamma_{G,\psi}:=\mathrm{Ind}_{U^F}^{G^F}((\overline{\mathbb{Z}}\textstyle{[\frac{1}{p}]})_\psi)\simeq \overline{\mathbb{Z}}\textstyle{[\frac{1}{p}]}G^F e_\psi\in\mathrm{Rep}_{\overline{\mathbb{Z}}[\frac{1}{p}]}(G^F)
\end{equation}
is called a Gelfand-Graev representation of $G$; the character of $\Gamma_{G,\psi}$ is $\mathrm{Ind}_{U^F}^{G^F}\psi$. 

All $\Gamma_{G,\psi}$ ($\psi\in\Psi$) are conjugate by elements of $T_\mathrm{ad}^F$: indeed, for any $\psi,\psi'\in\Psi$, there is a unique $x\in T_{\mathrm{ad}}^F$ such that $\psi'={}^x\psi$ (\S\,\ref{nondegenerate}); then $e_{\psi'}={}^x(e_\psi)$, and thus, upon identifying $\Gamma_{G,\psi}=\overline{\mathbb{Z}}\textstyle{[\frac{1}{p}]}G^F e_\psi$ via (\ref{Gelfand-Graev-eq}),
\begin{equation}\label{Gelfand-Graev-eq2}
\Gamma_{G,\psi'}=\overline{\mathbb{Z}}\textstyle{[\frac{1}{p}]}G^F e_{\psi'}={}^x(\overline{\mathbb{Z}}\textstyle{[\frac{1}{p}]}G^F e_\psi)={}^x(\Gamma_{G,\psi}).
\end{equation}
If the center $Z$ of $G$ is connected (that is, if $H^1(F,Z)=0$), then $T_{\mathrm{ad}}^F=T^F/Z^F$, so all $\Gamma_{G,\psi}$ are conjugate by elements of $T^F$ and are thus all isomorphic in $\mathrm{Rep}_{\overline{\mathbb{Z}}[\frac{1}{p}]}(G^F)$.
\end{context}

\begin{context}\label{Gelfand-Graev}
{\sc Theorem.}\;\cite[Thm.\;14.49]{Digne-Michel} --- 
{\it Let $\psi\in\Psi$. Then the $\overline{\mathbb{Q}}G^F$-module $\overline{\mathbb{Q}}\Gamma_{G,\psi}$ is multiplicity-free; more precisely, we have a $\overline{\mathbb{Q}}G^F$-module decomposition
\[
\overline{\mathbb{Q}}\Gamma_{G,\psi}=\displaystyle\bigoplus_{[x]\in (G^{\ast F^\ast})_{\mathrm{ss}}/\sim}\rho_{\psi,[x]}
\]
where for each $[x]\in (G^{\ast F^\ast})_{\mathrm{ss}}/\!\sim$, $\rho_{\psi,[x]}$ is an irreducible $\overline{\mathbb{Q}}G^F$-module whose character lies in the rational Lusztig series $\mathcal{E}(G^F,[x])$. 
}

\end{context}

\begin{context}\label{Gelfand-Graev-overlap}
{\sc Lemma.}\;--- {\it Let $\psi,\psi'\in\Psi$, and suppose that $\chi\in\mathrm{Irr}_{\overline{\mathbb{Q}}}(G^F)$ is a common irreducible character of the $\overline{\mathbb{Q}}G^F$-modules $\overline{\mathbb{Q}}\Gamma_{G,\psi}$ and $\overline{\mathbb{Q}}\Gamma_{G,\psi'}$. Then for the unique $y\in T_{\mathrm{ad}}^F$ such that $\psi'={}^y\psi$ (\S\,\ref{nondegenerate}), we have ${}^y\chi=\chi$.}

{\it Proof.} By \S\,\ref{Gelfand-Graev}, $\chi$ lies in some rational Lusztig series $\mathcal{E}(G^F,[x])$; in other words, there is an $F$-stable maximal torus $S$ and a $\theta\in\mathrm{Irr}_{\overline{\mathbb{Q}}}(S^F)$ such that $\theta$ corresponds to $[x]$ under the duality $\mathrm{Irr}_{\overline{\mathbb{Q}}}(S^F)\simeq S^{\ast F^\ast}$ (\S\,\ref{E-setup}), and such that $\langle\chi,R_S^G(\theta)\rangle_{G^F}\neq 0$. 

Let us show that ${}^y(R_S^G(\theta))=R_S^G(\theta)$. Indeed, as $y\in T_{\mathrm{ad}}^F$, if we choose a $t\in T$ such that $y=tZ$, then $t^{-1}F(t)\in Z\subset S$, so the Lang-Steinberg theorem enables us to write $t^{-1}F(t)=sF(s)^{-1}$ for some $s\in S$; thus, upon setting $x:=ts\in G^F$, 
\[
{}^y(R_S^G(\theta))={}^t(R_S^G(\theta))=R_{{}^tS}^G({}^t\theta)=R_{{}^xS}^G({}^x\theta)={}^x(R_S^G(\theta))=R_S^G(\theta).
\]

By the above discussion, we have
\[
\langle{}^y\chi,R_S^G(\theta)\rangle_{G^F}=\langle{}^y\chi,{}^y(R_S^G(\theta))\rangle_{G^F}=\langle\chi,R_S^G(\theta)\rangle_{G^F}\neq 0.
\]
On the other hand, with the aide of (\ref{Gelfand-Graev-eq2}),
\[
\langle{}^y\chi,\overline{\mathbb{Q}}\Gamma_{G,\psi'}\rangle_{G^F}=\langle{}^y\chi,{}^y(\overline{\mathbb{Q}}\Gamma_{G,\psi})\rangle_{G^F}=\langle\chi,\overline{\mathbb{Q}}\Gamma_{G,\psi}\rangle_{G^F}\neq 0.
\]
Therefore, $\gamma=\chi$ and $\gamma={}^y\chi$ both verify $\langle\gamma,R_S^G(\theta)\rangle_{G^F}\neq 0$ and $\langle\gamma,\overline{\mathbb{Q}}\Gamma_{G,\psi'}\rangle_{G^F}\neq 0$; but \S\,\ref{Gelfand-Graev} tells us that there is at most one such $\gamma\in\mathrm{Irr}_{\overline{\mathbb{Q}}}(G^F)$, so ${}^y\chi$ and $\chi$ must coincide. \qed
\end{context}

\subsection{The endomorphism algebra of $\Gamma_{G,\psi}$ and its $\mathbb{Z}$-model $\mathsf{E}_G$}

\begin{context}\label{E-definition}
{\it The endomorphism algebras $\overline{\mathbb{Z}}[\frac{1}{p}]\mathsf{E}_{G,\psi}$ and $\overline{\mathbb{Q}}\mathsf{E}_{G,\psi}$.}\;--- By \S\,\ref{Gelfand-Graev-definition}, all endomorphism algebras 
\[
\overline{\mathbb{Z}}[{\textstyle\frac{1}{p}}]\mathsf{E}_{G,\psi}:=\mathrm{End}_{\overline{\mathbb{Z}}[\frac{1}{p}]G^F}(\Gamma_{G,\psi})\;\;(\psi\in\Psi)
\]
are isomorphic as $\overline{\mathbb{Z}}[\frac{1}{p}]$-algebras. Later in \S\,\ref{E-Zalgebra}, we shall introduce a $\mathbb{Z}$-algebra $\mathsf{E}_G$ which is independent of $\psi$ and is such that $\overline{\mathbb{Z}}[\frac{1}{p}]\mathsf{E}_{G}\simeq \overline{\mathbb{Z}}[\frac{1}{p}]\mathsf{E}_{G,\psi}$ for all $\psi\in\Psi$. 

Let $\overline{\mathbb{Q}}\mathsf{E}_{G,\psi}:=\overline{\mathbb{Q}}\otimes_{\overline{\mathbb{Z}}[\frac{1}{p}]}\overline{\mathbb{Z}}[\frac{1}{p}]\mathsf{E}_{G,\psi}$. By \S\,\ref{Gelfand-Graev} and Schur's lemma, we may decompose the $\overline{\mathbb{Q}}$-algebra $\overline{\mathbb{Q}}\mathsf{E}_{G,\psi}$ as:
\begin{equation}\label{E-decomp-eq1}
\overline{\mathbb{Q}}\mathsf{E}_{G,\psi}=\mathrm{End}_{\overline{\mathbb{Q}}G^F}(\overline{\mathbb{Q}}\Gamma_{G,\psi})\simeq\prod_{[x]\in (G^{\ast F^\ast})_{\mathrm{ss}}/\sim}\mathrm{End}_{\overline{\mathbb{Q}}G^F}(\rho_{\psi,[x]})\simeq \overline{\mathbb{Q}}^{(G^{\ast F^\ast})_{\mathrm{ss}}/\sim}.
\end{equation}
Thus $\overline{\mathbb{Q}}\mathsf{E}_{G,\psi}$ and $\overline{\mathbb{Z}}[\frac{1}{p}]\mathsf{E}_{G,\psi}$ are commutative reduced rings. In terms of algebraic variety, (\ref{E-decomp-eq1}) identifies $\mathrm{Specm}(\overline{\mathbb{Q}}\mathsf{E}_G)$ (consisting of maximal ideals of $\overline{\mathbb{Q}}\mathsf{E}_G$) with $(G^{\ast F^\ast})_{\mathrm{ss}}/\!\sim$.   

\end{context}

\begin{context}\label{E-decomposition}
{\it Structures of $\overline{\mathbb{Z}}[\frac{1}{p}]\mathsf{E}_{G,\psi}$ and $\overline{\mathbb{Q}}\mathsf{E}_{G,\psi}$ via idempotents.}\;--- Using the identification $\Gamma_{G,\psi}=\overline{\mathbb{Z}}[\frac{1}{p}]G^Fe_\psi$ in (\ref{Gelfand-Graev-eq}), we shall identify 
\[
\overline{\mathbb{Z}}[\textstyle{\frac{1}{p}}]\mathsf{E}_{G,\psi}=e_\psi\overline{\mathbb{Z}}[\textstyle{\frac{1}{p}}]G^Fe_\psi\subset \overline{\mathbb{Z}}[\textstyle{\frac{1}{p}}]G^F
\]
as $\overline{\mathbb{Z}}[\frac{1}{p}]$-algebras: indeed, there is a $\overline{\mathbb{Z}}[\frac{1}{p}]$-algebra anti-isomorphism 
\begin{equation}\label{endo-idemp-eq}
\overline{\mathbb{Z}}[{\textstyle\frac{1}{p}}]\mathsf{E}_{G,\psi}=\mathrm{End}_{\overline{\mathbb{Z}}[\frac{1}{p}]G^F}(\overline{\mathbb{Z}}\textstyle{[\frac{1}{p}]}G^F e_\psi)\xrightarrow{\;\;\sim\;\;}e_\psi\overline{\mathbb{Z}}[\textstyle{\frac{1}{p}}]G^Fe_\psi,\quad \theta\longmapsto\theta(e_\psi),
\end{equation}
but $\overline{\mathbb{Z}}[\frac{1}{p}]\mathsf{E}_{G,\psi}$ is a commutative ring (\S\,\ref{E-definition}), so (\ref{endo-idemp-eq}) is also an isomorphism.

(a) By \cite[Prop.\;11.30]{Curtis-Reiner}, $\overline{\mathbb{Z}}[\frac{1}{p}]\mathsf{E}_{G,\psi}=e_\psi\overline{\mathbb{Z}}[\frac{1}{p}]G^Fe_\psi$ is a free $\overline{\mathbb{Z}}[\frac{1}{p}]$-module with the following $\overline{\mathbb{Z}}[\frac{1}{p}]$-linear basis: let $\{x_1,\cdots,x_r\}\subset G^F$ be a set of representatives of $U^F\backslash G^F\!/ U^F$ (so $G^F$ is the disjoint union of $\{U^Fx_iU^F:i=1,\cdots,r\}$), and let
\[
J=\{j:1\leq j\leq r, \,{}^{x_j}\psi=\psi\;\mathrm{on}\;U^F\cap{}^{x_j}(U^F)\};
\]
for each $x\in G^F$, define $\mathrm{ind}\,x=[U^F:U^F\cap {}^x(U^F)]$ (which is a power of $p$); then $\{(\mathrm{ind}\,x_j)e_\psi x_je_\psi:j\in J\}$ is a $\overline{\mathbb{Z}}[\frac{1}{p}]$-linear basis of $\overline{\mathbb{Z}}[\frac{1}{p}]\mathsf{E}_{G,\psi}$.

(b) By \cite[Sec.\;3]{Curtis}, $\overline{\mathbb{Q}}\mathsf{E}_{G,\psi}$ may be described via idempotents: for $[x]\in (G^{\ast F^\ast})_{\mathrm{ss}}/\!\sim$, denote by $\chi_{\psi,[x]}$ the character associated to the $\overline{\mathbb{Q}}G^F$-module $\rho_{\psi,[x]}$ in \S\,\ref{Gelfand-Graev}, and set 
\[
e_{\psi,[x]}^G:=\displaystyle\frac{\chi_{\psi,[x]}(1)}{|G^F|}\sum_{g\in G}\chi_{\psi,[x]}(g^{-1})g\in\overline{\mathbb{Q}}G^F\mbox{\quad and \quad} e_{\psi,[x]}^{\mathsf{E}}:=e_{\psi,[x]}^Ge_\psi\in\overline{\mathbb{Q}}\mathsf{E}_{G,\psi},
\]
which are respectively primitive central idempotents of $\overline{\mathbb{Q}}G^F$ and of $\overline{\mathbb{Q}}\mathsf{E}_{G,\psi}$; then the inclusions $\overline{\mathbb{Q}}e_{\psi,[x]}^{\mathsf{E}}\subset\overline{\mathbb{Q}}\mathsf{E}_{G,\psi}$ induce a $\overline{\mathbb{Q}}$-algebra isomorphism 
\begin{equation}\label{E-decomp-eq2}
\overline{\mathbb{Q}}\mathsf{E}_{G,\psi}=\displaystyle\prod\limits_{[x]\in (G^{\ast F^\ast})_{\mathrm{ss}}/\sim}\overline{\mathbb{Q}}e_{\psi,[x]}^{\mathsf{E}}.
\end{equation}
Combining (\ref{E-decomp-eq1}) and (\ref{E-decomp-eq2}), we obtain $\overline{\mathbb{Q}}$-algebra isomorphisms
\[
\overline{\mathbb{Q}}\mathsf{E}_{G,\psi}\simeq \prod\limits_{[x]\in (G^{\ast F^\ast})_{\mathrm{ss}}/\sim}\overline{\mathbb{Q}}e_{\psi,[x]}^{\mathsf{E}}\simeq \overline{\mathbb{Q}}^{(G^{\ast F^\ast})_{\mathrm{ss}}/\sim},
\]
where each $e_{\psi,[x]}^{\mathsf{E}}$ corresponds to the characteristic function $\mathbf{1}_{\{[x]\}}$ on $(G^{\ast F^\ast})_{\mathrm{ss}}/\!\sim$.
\end{context}

\begin{context}\label{E-Zalgebra}
{\it Definition of the $\mathbb{Z}$-algebra $\mathsf{E}_G$.}\;--- The decomposition $U^F/U_\bullet^F=\prod\limits_{i\in\Delta/F}U_i^F$ (\S\,\ref{nondegenerate}) identifies each $U_i^F$ as a subgroup of $U^F/U_\bullet^F$. Consider the (left) $\mathbb{Z}U^F$-module
\[
\mathcal{F}_0=\{f:U^F/U_\bullet^F\longrightarrow\mathbb{Z}:\sum_{x\in U_i^F}f(x)=0\mbox{ for each }i\in\Delta/F\}
\]
on which $U^F$ acts by left translation: $u\cdot f:=f(u^{-1}(\cdot))$ ($u\in U^F$, $f\in\mathcal{F}_0$). The groups $T^F$ and $T_{\mathrm{ad}}^F$ both acts on each $U_i^F$ $(i\in\Delta/F)$ by left adjoint action; this induces left $T^F$-actions and left $T_{\mathrm{ad}}^F$-actions on $\mathcal{F}_0$, $\mathrm{Ind}_{U^F}^{G^F}(\mathcal{F}_0)$ and $\mathrm{End}_{{\mathbb{Z}}G^F}(\mathrm{Ind}_{U^F}^{G^F}(\mathcal{F}_0))$. The center of the $\mathbb{Z}$-algebra $\mathrm{End}_{{\mathbb{Z}}G^F}(\mathrm{Ind}_{U^F}^{G^F}(\mathcal{F}_0))$, denoted by $Z\left(\mathrm{End}_{{\mathbb{Z}}G^F}(\mathrm{Ind}_{U^F}^{G^F}(\mathcal{F}_0))\right)$, is fixed by the $T^F$-action on $\mathrm{End}_{{\mathbb{Z}}G^F}(\mathrm{Ind}_{U^F}^{G^F}(\mathcal{F}_0))$ and thus admits an $H^1(F,Z)$-action (recall that $H^1(F,Z)\simeq T_{\mathrm{ad}}^F/T^F$); we then set
\begin{equation}\label{E-Zmodel-def}
\mathsf{E}_G:=Z\left(\mathrm{End}_{{\mathbb{Z}}G^F}(\mathrm{Ind}_{U^F}^{G^F}(\mathcal{F}_0))\right)^{H^1(F,Z)},
\end{equation}
which consists of elements of $Z\left(\mathrm{End}_{{\mathbb{Z}}G^F}(\mathrm{Ind}_{U^F}^{G^F}(\mathcal{F}_0))\right)$ fixed by the $H^1(F,Z)$-action. This $\mathsf{E}_G$ is a $\mathbb{Z}$-algebra and is also a free $\mathbb{Z}$-module of finite rank.

After describing the $\overline{\mathbb{Z}}[\frac{1}{p}]G^F$-module $\mathrm{Ind}_{U^F}^{G^F}(\overline{\mathbb{Z}}[\frac{1}{p}]\mathcal{F}_0)$ in \S\,\ref{lemma-E-indF0}, we shall show in \S\,\ref{lemma-E-compatibility} that the $\mathbb{Z}$-algebra $\mathsf{E}_G$ in (\ref{E-Zmodel-def}) is consistent with the definition of $\overline{\mathbb{Z}}[\frac{1}{p}]\mathsf{E}_{G,\psi}$ in \S\,\ref{E-definition}.
\end{context}

\begin{context}\label{lemma-E-indF0}
{\sc Lemma.}\;---
{\it 
Let $\Lambda=\overline{\mathbb{Z}}[\frac{1}{p}]$. The natural inclusion $\mathcal{F}_0\subset\mathbb{Z}U^F$ gives an  embedding $\mathrm{Ind}_{U^F}^{G^F}(\Lambda\mathcal{F}_0)\subset\Lambda G^F$, under which the actions of $T^F$ and $T_{\mathrm{ad}}^F$ on $\mathrm{Ind}_{U^F}^{G^F}(\Lambda\mathcal{F}_0)$ (\S\,\ref{E-Zalgebra}) are exactly the restrictions of their respective adjoint actions on $\Lambda G^F$.
Moreover, we have a decomposition of $\Lambda G^F$-modules:}
\[
\mathrm{Ind}_{U^F}^{G^F}(\Lambda\mathcal{F}_0)=\bigoplus_{\psi\in\Psi}\Lambda G^F e_\psi=\bigoplus_{\psi\in\Psi}\Gamma_{G,\psi}\in\mathrm{Rep}_{\Lambda}(G^F).
\]

{\it Proof.} The description of the actions of $T^F$ and $T_{\mathrm{ad}}^F$ on $\mathrm{Ind}_{U^F}^{G^F}(\Lambda\mathcal{F}_0)$ is clear from construction, so let us prove the decomposition formula directly.

(1) (Compare \cite[(2.4.8)]{Digne-Lehrer-Michel}.) Each linear character $\psi:U^F/U_\bullet^F=\prod\limits_{i\in\Delta/F}U_i^F\longrightarrow{\mathbb{Z}}^\times$ is decomposed as $\psi=\prod\limits_{i\in\Delta/F}\psi_i$ with $\psi_i=\psi|_{U_i^F}$ ($i\in\Delta/F$); such a $\psi$ gives an element in $\Psi$ if and only if each $\psi_i$ is non-trivial. 

(2) For each $\psi\in\Psi$, let $e_\psi$ be its central primitive idempotent of $\Lambda U^F$ as in (\ref{idempotent-eq}). Using (1) and the orthogonality of characters, one can check that $\{e_\psi\,|\,\psi\in\Psi\}$ is a $\overline{\mathbb{Q}}$-linearly independent subset of $\overline{\mathbb{Q}}\mathcal{F}_0$. By definition of $\mathcal{F}_0$ and (1), we have
\[
\dim_{\overline{\mathbb{Q}}}(\overline{\mathbb{Q}}\mathcal{F}_0)=\prod_{i\in\Delta/F}(|U_i^F|-1)=|\Psi|,
\]
so $\{e_\psi\,|\,\psi\in\Psi\}$ is in fact a $\overline{\mathbb{Q}}$-linear basis of $\overline{\mathbb{Q}}\mathcal{F}_0$; identifying $\overline{\mathbb{Q}}\mathcal{F}_0\subset\overline{\mathbb{Q}}U^F$, we have the following decomposition of $\overline{\mathbb{Q}}\mathcal{F}_0$:
\begin{equation}\label{F0space-eq1}
\overline{\mathbb{Q}}\mathcal{F}_0=\bigoplus_{\psi\in\Psi}\overline{\mathbb{Q}}e_\psi;\;\;\mbox{ for $f\in\overline{\mathbb{Q}}\mathcal{F}_0$, }f=\sum_{\psi\in\Psi} fe_\psi\mbox{ with each $fe_\psi\in\overline{\mathbb{Q}}e_\psi$.}
\end{equation}
As all $e_\psi$ ($\psi\in\Psi$) lie in $\Lambda U^F$, the decomposition (\ref{F0space-eq1}) also holds in $\Lambda$-coefficients, from which the desired decomposition of $\mathrm{Ind}_{U^F}^{G^F}(\Lambda\mathcal{F}_0)$ then follows. \qed
\end{context}

\begin{context}\label{lemma-E-compatibility}
{\sc Lemma.}\;--- 
{\it Let $\Lambda=\overline{\mathbb{Z}}[\frac{1}{p}]$. For each $\psi\in\Psi$, $\Lambda\mathsf{E}_G\simeq\Lambda\mathsf{E}_{G,\psi}$ as $\Lambda G^F$-algebras. Explicitly: upon identifying $\Lambda\mathsf{E}_{G,\psi}=e_\psi\Lambda G^F e_\psi$ via (\ref{endo-idemp-eq}) and identifying analogously
\[
\mathrm{End}_{{\Lambda}G^F}(\mathrm{Ind}_{U^F}^{G^F}(\Lambda\mathcal{F}_0))=\bigoplus_{\psi,\psi'\in\Psi}e_{\psi'}\Lambda G^F e_{\psi}\subset\Lambda G^F
\] 
we have a $\Lambda$-algebra isomorphism}
\[
\Lambda\mathsf{E}_{G,\psi}\xrightarrow{\;\;\sim\;\;}\Lambda\mathsf{E}_G,\quad e_\psi f e_\psi\longmapsto \sum_{y\in T_{\mathrm{ad}}^F}e_{{}^y\psi}\cdot{}^y f\cdot e_{{}^y\psi}\;\;(f\in\Lambda G^F).
\]

{\it Proof.} (1) An element in the center of $\mathrm{End}_{{\Lambda}G^F}(\mathrm{Ind}_{U^F}^{G^F}(\Lambda\mathcal{F}_0))$ commutes with all $e_\psi$ ($\psi\in\Psi$), so the orthogonality of idempotents implies that
\begin{equation}\label{center-endoF0-eq1}
Z\left(\mathrm{End}_{{\Lambda}G^F}(\mathrm{Ind}_{U^F}^{G^F}(\Lambda\mathcal{F}_0))\right)\subset \bigoplus_{\psi\in\Psi}e_\psi\Lambda G^F e_\psi.
\end{equation}

(2) Let $y\in T_{\mathrm{ad}}^F$. For $\theta\in \mathrm{End}_{{\Lambda}G^F}(\mathrm{Ind}_{U^F}^{G^F}(\Lambda\mathcal{F}_0))$, its image $y\cdot\theta$ under the action of $y$ on $\mathrm{End}_{{\Lambda}G^F}(\mathrm{Ind}_{U^F}^{G^F}(\Lambda\mathcal{F}_0))$ is given by $(y\cdot\theta)(\varphi)={}^y(\theta({}^{y^{-1}}\varphi))$ $(\varphi\in \mathrm{Ind}_{U^F}^{G^F}(\Lambda\mathcal{F}_0))$. Then the subring $\bigoplus\limits_{\psi\in\Psi}e_\psi\Lambda G^F e_\psi$ is invariant under this action of $y$: indeed, let $\theta\in\bigoplus\limits_{\psi\in\Psi}e_\psi\Lambda G^F e_\psi$ and write $\theta=\sum\limits_{\psi\in\Psi}\theta_\psi$ where each $\theta_\psi=e_\psi\theta e_\psi\in e_\psi\Lambda G^Fe_\psi$; then it can be checked that $y\cdot\theta=\sum\limits_{\psi\in\Psi}(y\cdot\theta)_\psi$ with each $(y\cdot\theta)_\psi=y\cdot(\theta_{{}^{y^{-1}}\psi})={}^y(\theta_{{}^{y^{-1}}\psi})\in e_\psi\Lambda G^F e_\psi$. We thus deduce from (\ref{center-endoF0-eq1}) that
\begin{equation}\label{center-endoF0-eq2}
\Lambda\mathsf{E}_G\subset \left(\bigoplus_{\psi\in\Psi}e_\psi\Lambda G^F e_\psi\right)^{T_{\mathrm{ad}}^F}=:\mathcal{A}.
\end{equation}

(3) For each $\theta\in \mathcal{A}$, if we write $\theta=\sum\limits_{\psi\in\Psi}\theta_\psi$ with each $\theta_\psi\in e_\psi\Lambda G^F e_\psi$ as in (2), then $\theta_\psi={}^y(\theta_{{}^{y^{-1}}\psi})$ for all $y\in T_{\mathrm{ad}}^F$; as $T_{\mathrm{ad}}^F$ acts regularly on $\Psi$ (\S\,\ref{nondegenerate}), we see that for every $\psi\in\Psi$ we have a ring isomorphism 
\begin{equation}\label{center-endoF0-eq3}
\Lambda\mathsf{E}_{G,\psi}\xrightarrow{\;\;\sim\;\;}\mathcal{A},\quad e_\psi f e_\psi\longmapsto \sum\limits_{y\in T_{\mathrm{ad}}^F}e_{{}^y\psi}\cdot {}^y f\cdot e_{{}^y\psi}\;\;(f\in\Lambda G^F).
\end{equation}

(4) By virtue of (\ref{center-endoF0-eq2}) and (\ref{center-endoF0-eq3}), in order to complete the proof of lemma, it remains to establish the inclusion
\begin{equation}\label{center-endoF0-eq4}
\overline{\mathbb{Q}}\mathcal{A}=\left(\bigoplus_{\psi\in\Psi}e_\psi\overline{\mathbb{Q}} G^F e_\psi\right)^{T_{\mathrm{ad}}^F}\subset Z\left(\bigoplus_{\psi,\psi'\in\Psi}e_{\psi'}\overline{\mathbb{Q}} G^F e_{\psi}\right).
\end{equation}
For each $\psi,\psi'\in\Psi$, \S\,\ref{Gelfand-Graev} (with notations in \S\,\ref{E-decomposition}) tells us that the space $e_{\psi'}\overline{\mathbb{Q}}G^Fe_\psi$ is the direct sum of $e_{\psi'}\overline{\mathbb{Q}}G^Fe_{\psi,[x]}^Ge_\psi$ for $[x]$ ranging over a (possibly empty) subset of $(G^{\ast F^\ast})_{\mathrm{ss}}/\!\sim$. To show (\ref{center-endoF0-eq4}), it then suffices to show that 
\begin{equation}\label{center-endoF0-eq5}
e_{{}^y\psi}fe_{\psi,[x']}^Ge_\psi\cdot e_\psi e_{\psi,[x]}^Ge_\psi= e_{{}^y\psi} ({}^ye_{\psi,[x]}^G)e_{{}^y\psi}\cdot e_{{}^y\psi}fe_{\psi,[x']}^Ge_\psi
\end{equation}
for every $\psi\in\Psi$, $y\in T_{\mathrm{ad}}^F$, $[x],[x']\in(G^{\ast F^\ast})_{\mathrm{ss}}/\!\sim$ and $f\in\overline{\mathbb{Q}}G^F$. 

To prove (\ref{center-endoF0-eq5}), suppose that the element $e_{{}^y\psi}fe_{\psi,[x']}^Ge_\psi\in \mathrm{Hom}_{\overline{\mathbb{Q}}G^F}(\overline{\mathbb{Q}}\Gamma_{{}^y\psi},\rho_{\psi,[x']})$ is not zero (as the opposite case is trivial). Then $\chi_{\psi,[x']}$ is a common irreducible character of $\overline{\mathbb{Q}}\Gamma_{{}^y\psi}$ and $\overline{\mathbb{Q}}\Gamma_{\psi}$, so \S\,\ref{Gelfand-Graev-overlap} gives us ${}^y(\chi_{\psi,[x']})=\chi_{\psi,[x']}$ and then ${}^ye_{\psi,[x']}^G=e_{\psi,[x']}^G$. We may also suppose that $[x]=[x']$, for otherwise both side (\ref{center-endoF0-eq5}) are zero, using the fact that the idempotents $e_{\psi,[x]}^G$ are central and orthogonal in $\overline{\mathbb{Q}}G^F$. Under the above assumptions, the two sides of (\ref{center-endoF0-eq5}) are both $e_{{}^y\psi}fe_{\psi,[x']}^Ge_\psi$; this justifies (\ref{center-endoF0-eq5}).

The equality (\ref{center-endoF0-eq5}) being proved, we obtain the inclusion (\ref{center-endoF0-eq4}) as well as the equality $\Lambda\mathsf{E}_G=\mathcal{A}$; thus (\ref{center-endoF0-eq3}) gives the desired $\Lambda$-algebra isomorphism. \qed
\end{context}

\begin{context}\label{Galois-action}
{\it The Galois action on $\overline{\mathbb{Q}}\mathsf{E}_G$.}\;--- The Galois group $\mathrm{Gal}(\overline{\mathbb{Q}}/\mathbb{Q})$ acts on the $\overline{\mathbb{Q}}$-coefficient of $\overline{\mathbb{Q}}\mathsf{E}_G$; this gives us an action of $\mathrm{Gal}(\overline{\mathbb{Q}}/\mathbb{Q})$ on $\overline{\mathbb{Q}}\mathsf{E}_G$, and this action coincide with that induced by the identification
\[
\overline{\mathbb{Q}}\mathsf{E}_G=Z\left(\mathrm{End}_{\overline{\mathbb{Q}}G^F}(\mathrm{Ind}_{U^F}^{G^F}(\overline{\mathbb{Q}}\mathcal{F}_0))\right)^{H^1(F,Z)}
\] 
and by the action of $\mathrm{Gal}(\overline{\mathbb{Q}}/\mathbb{Q})$ on the $\overline{\mathbb{Q}}$-coefficient of $\overline{\mathbb{Q}}\mathcal{F}_0$.

In terms of idempotents, the above Galois action on $\overline{\mathbb{Q}}\mathsf{E}_G$ is, upon fixing a $\psi\in\Psi$ and identify $\overline{\mathbb{Q}}\mathsf{E}_{G}=\overline{\mathbb{Q}}\mathsf{E}_{G,\psi}=e_\psi\overline{\mathbb{Q}}G^F e_\psi$ (\S\,\ref{lemma-E-compatibility}): 
\[
\sigma\in\mathrm{Gal}(\overline{\mathbb{Q}}/\mathbb{Q}),\, e_\psi f e_\psi\in\overline{\mathbb{Q}}\mathsf{E}_{G,\psi}\;(f\in\overline{\mathbb{Q}}G^F)\Longrightarrow \sigma\cdot(e_\psi f e_\psi)=e_\psi\cdot\sigma({}^{y_\sigma^{-1}}f)\cdot e_\psi,
\]
where $y_\sigma$ is the unique element in $T_{\mathrm{ad}}^F$ such that $\sigma\cdot\psi={}^{y_\sigma}\psi\in\Psi$ ($\sigma$ acts on $\overline{\mathbb{Z}}^\times$ and hence on $\Psi$; then recall from \S\,\ref{nondegenerate} that $T_{\mathrm{ad}}^F$ acts regularly on $\Psi$), and where, for each $h\in\overline{\mathbb{Q}}G^F$, $\sigma(h)$ denotes its image under the action of $\sigma$ on the $\overline{\mathbb{Q}}$-coefficient of $\overline{\mathbb{Q}}G^F$.
\end{context}

\subsection{Curtis homomorphisms}

\begin{context}\label{Curtis-definition}
{\it The Curtis homomorphism $\mathrm{Cur}_S^G$.}\;\cite[Sec.\;4]{Curtis} --- Hereafter, fix a $\psi\in\Psi$. For $S$ an $F$-stable maximal torus of $G$, Curtis has constructed a $\overline{\mathbb{Q}}$-algebra homomorphism 
\[
\mathrm{Cur}_S^G:\overline{\mathbb{Q}}\mathsf{E}_{G,\psi}\longrightarrow\overline{\mathbb{Q}}S^F
\]
(which he called $f_S$ in his paper), which is the unique $\overline{\mathbb{Q}}$-algebra homomorphism from $\overline{\mathbb{Q}}\mathsf{E}_{G,\psi}$ to $\overline{\mathbb{Q}}S^F$ such that, for every irreducible character $\theta\in\mathrm{Irr}_{\overline{\mathbb{Q}}}(S^F)$ corresponding to $[x]\in (G^{\ast F^\ast})_{\mathrm{ss}}/\!\sim$ under the duality $\mathrm{Irr}_{\overline{\mathbb{Q}}}(S^F)\simeq S^{\ast F^\ast}$ (\S\,\ref{setup-global}), the following diagram of $\overline{\mathbb{Q}}$-algebras commutes: 
\begin{equation}\label{Curtis-diagram}
\begin{tikzcd}
\overline{\mathbb{Q}}\mathsf{E}_{G,\psi}\arrow[rr,"\chi_{\psi,[x]}|_{\overline{\mathbb{Q}}\mathsf{E}_{G,\psi}}"]\arrow[d,swap,"\mathrm{Cur}_S^G"] & &\overline{\mathbb{Q}}\\
\overline{\mathbb{Q}}S^F\arrow[urr,swap,"\theta"] & 
\end{tikzcd}
\end{equation}

The homomorphism $\mathrm{Cur}_S^G$ is determined by \S\,\ref{E-decomposition}(a) and by the following formula:

{\it Curtis' formula.} For all $n\in G^F$, $\mathrm{Cur}_S^G(e_\psi ne_\psi)=\displaystyle\sum_{s\in S^F}\mathrm{Cur}_S^G(e_\psi ne_\psi)(s)s\in\overline{\mathbb{Q}}S^F$ with
\[
\mathrm{Cur}_S^G(e_\psi ne_\psi)(s)=\frac{1}{\langle Q_S^G,\Gamma_{G,\psi}\rangle_{G^F}}\frac{1}{|U^F|}\frac{1}{|(C_G(s)^\circ)^F|}\sum_{\substack{g\in G^F,\,u\in U^F\\ (gung^{-1})_{\mathrm{ss}}=s}}\psi(u^{-1})Q_S^{C_G(s)^\circ}((gung^{-1})_{\mathrm{u}}).
\]
In the above formula, $Q_S^G:G^F\longrightarrow\mathbb{Z}$ is the Green function, which is the function supported on the set of unipotent elements of $G^F$ defined by $Q_S^G(u)=R_S^G(\mathrm{Id})(u)$ for unipotent elements $u$ of $G^F$ (see \cite[Def.\;4.1]{Deligne-Lusztig}); $C_G(s)^\circ$ means the connected identity component of the centralizer of $s$ in $G$; for $g\in G^F$, $g_{\mathrm{ss}}$ (resp.\;$g_{\mathrm{u}}$) denotes its semisimple part (resp.\;unipotent part).

\end{context}

\begin{context}\label{Curtis-integral}
{\sc Lemma.}\;--- 
{\it 
Let $S$ be an $F$-stable maximal torus of $G$.
\begin{enumerate}[\normalfont (a)]
\item $\mathrm{Cur}_S^G$ is defined over $\overline{\mathbb{Z}}[\frac{1}{p}]$ in the sense that $\mathrm{Cur}_S^G(\overline{\mathbb{Z}}[\frac{1}{p}]\mathsf{E}_{G,\psi})\subset\overline{\mathbb{Z}}[\frac{1}{p}]S^F$.
\item Identify $\overline{\mathbb{Z}}[\frac{1}{p}]\mathsf{E}_{G}=\overline{\mathbb{Z}}[\frac{1}{p}]\mathsf{E}_{G,\psi}$ (\S\,\ref{lemma-E-compatibility}). Then $\mathrm{Cur}_S^G:\overline{\mathbb{Q}}\mathsf{E}_G\longrightarrow\overline{\mathbb{Q}} S^F$ is equivariant under the action of $\mathrm{Gal}(\overline{\mathbb{Q}}/\mathbb{Q})$ (\S\,\ref{Galois-action}). Thus (with the aide of (a)) $\mathrm{Cur}_S^G$ is defined over $\mathbb{Z}[\frac{1}{p}]$: $\mathrm{Cur}_S^G(\mathbb{Z}[\frac{1}{p}]\mathsf{E}_{G})\subset\mathbb{Z}[\frac{1}{p}]S^F$.
\end{enumerate}
}

{\it Proof.} (a) Using the $\overline{\mathbb{Z}}[\frac{1}{p}]$-linear basis of $\overline{\mathbb{Z}}[\frac{1}{p}]\mathsf{E}_{G,\psi}$ in \S\,\ref{E-decomposition}(a), it suffices to show that $\mathrm{Cur}_S^G(e_\psi n e_\psi)\in\overline{\mathbb{Z}}[\frac{1}{p}]S^F$ for every $n\in G^F$. In Curtis' formula for $\mathrm{Cur}_S^G(e_\psi n e_\psi)$ (\S\,\ref{Curtis-definition}), observe that $\langle Q_S^G,\Gamma_{G,\psi}\rangle_{G^F}=\langle R_S^G(\mathrm{Id}),\Gamma_{G,\psi}\rangle_{G^F}=\pm1$ (\S\,\ref{Gelfand-Graev}) and that $|U^F|$ is a power of $p$, so to prove this lemma it remains to show that the number
\[
C(s):=\displaystyle\frac{1}{|(C_G(s)^\circ)^F|}\sum_{\substack{g\in G^F,\,u\in U^F\\ (gung^{-1})_{\mathrm{ss}}=s}}\psi(u^{-1})Q_S^{C_G(s)^\circ}((gung^{-1})_{\mathrm{u}})
\]
lies in $\overline{\mathbb{Z}}[\frac{1}{p}]$ for all $s\in S^F$. If $(gung^{-1})_{\mathrm{ss}}=s$ ($g,n\in G^F$, $u\in U^F$), for every $c\in C_G(s)^\circ$ we have $((cg)un(cg)^{-1})_{\mathrm{ss}}=s$ and $((cg)un(cg)^{-1})_{\mathrm{u}}=c.(gung^{-1})_{\mathrm{u}}.c^{-1}$, so
\[
C(s)=\sum_{\substack{g\in (C_G(s)^\circ)^F\backslash G^F,\,u\in U^F\\ (gung^{-1})_{\mathrm{ss}}=s}}\psi(u^{-1})Q_S^{C_G(s)^\circ}((gung^{-1})_{\mathrm{u}})\in\overline{\mathbb{Z}}[{\textstyle\frac{1}{p}}].
\]

(b) We have to show that 
\begin{equation}\label{Curtis-Galois-eq}
\sigma(\mathrm{Cur}_S^G(e_\psi ne_\psi)(s))=\mathrm{Cur}_S^G(\sigma \cdot(e_\psi ne_\psi))(s)
\end{equation}
for every $n\in G^F$, $\sigma\in\mathrm{Gal}(\overline{\mathbb{Q}}/\mathbb{Q})$ and $s\in S^F$. Again, we apply Curtis' formula: using \S\,\ref{Galois-action} and the notation therein, we have $\sigma\cdot\psi={}^{y_\sigma}\psi$ with $y_\sigma\in T_{\mathrm{ad}}^F$; choose a $t\in T$ such that $y_\sigma=tZ$, and choose a $c\in S$ such that $x:=t^{-1} c\in G^F$ (see \S\,\ref{Gelfand-Graev-overlap}); then
\begin{align*}
&\;\langle Q_S^G,\Gamma_{G,\psi}\rangle_{G^F}|U^F||(C_G(s)^{\circ})^F|\cdot\sigma(\mathrm{Cur}_S^G(e_\psi ne_\psi)(s))\\
=&\sum_{\substack{g\in G^F,\,u\in U^F\\ (gung^{-1})_{\mathrm{ss}}=s}}\sigma(\psi(u^{-1}))Q_S^{C_G(s)^\circ}((gung^{-1})_{\mathrm{u}})
\\
=&\sum_{\substack{g\in G^F,\,u\in U^F\\ (gung^{-1})_{\mathrm{ss}}=s}}({}^{y_\sigma}\psi)(u^{-1})Q_S^{C_G(s)^\circ}((gung^{-1})_{\mathrm{u}})\quad\\
=&\sum_{\substack{g\in G^F,\,u\in U^F\\ (gu.{}^{y_\sigma^{-1}}n.g^{-1})_{\mathrm{ss}}={}^{y_\sigma^{-1}}s}}\psi(u^{-1})Q_S^{C_G(s)^\circ}({}^{y_\sigma}((gu.{}^{y_\sigma^{-1}}n.g^{-1})_{\mathrm{u}}))\quad(u\mapsto{}^{y_\sigma}u;\;g\mapsto {}^{y_\sigma}g)\\
=&\sum_{\substack{g\in G^F,\,u\in U^F\\ {}^x((gu.{}^{y_\sigma^{-1}}n.g^{-1})_{\mathrm{ss}})={}^{t^{-1}}s}}\psi(u^{-1})Q_S^{C_G(s)^\circ}({}^{y_\sigma x}((gu.{}^{y_\sigma^{-1}}n.g^{-1})_{\mathrm{u}}))\quad(g\mapsto xg)\\
=&\sum_{\substack{g\in G^F,\,u\in U^F\\(gu.{}^{y_\sigma^{-1}}n.g^{-1})_{\mathrm{ss}}=s}}\psi(u^{-1})Q_S^{C_G(s)^\circ}((gu.{}^{y_\sigma^{-1}}n.g^{-1})_{\mathrm{u}})\quad({}^{x^{-1}{y_\sigma^{-1}}}(Q_S^{C_G(s)^\circ})=Q_S^{C_G(s)^\circ};\mbox{ see \S\,\ref{Gelfand-Graev-overlap}})\\
=&\;\langle Q_S^G,\Gamma_{G,\psi}\rangle_{G^F}|U^F||(C_G(s)^{\circ})^F|\cdot\mathrm{Cur}_S^G(\sigma \cdot(e_\psi ne_\psi))(s). \tag*{$\square$}
\end{align*}

\end{context}

\begin{context}\label{Curtis-idempotent}\label{Curtis-geometry}
{\sc Lemma.}\;--- 
{\it 
Let $S$ be an $F$-stable maximal torus of $G$. For $z\in S^{\ast F^\ast}$, let $\widehat{z}:S^F\longrightarrow\overline{\mathbb{Q}}$ be its corresponding irreducible character of $S^F$ obtained from the duality $S^{\ast F^\ast}\simeq\mathrm{Irr}_{\overline{\mathbb{Q}}}(S^F)$ (\S\,\ref{setup-global}), and let $e_{\widehat{z}}^S:=\displaystyle\frac{1}{|S^F|}\sum_{t\in S^F}\widehat{z}(t^{-1})t\in\overline{\mathbb{Q}}S^F$ be the primitive central idempotent associated to $\widehat{z}$. Then 
\begin{equation}\label{Curtis-eq1}
\mathrm{Cur}_S^G(e_{\psi,[x]}^{\mathsf{E}})=\displaystyle\sum_{z\in [x]\cap S^{\ast F^\ast}}e_{\widehat{z}}^S\in \overline{\mathbb{Q}}S^F \quad\mathit{\;for\;every\;} [x]\in (G^{\ast F^\ast})_{\mathrm{ss}}/\!\sim.
\end{equation}
As a consequence, $\mathrm{Cur}_S^G$ may be interpreted as a restriction map: using the identification $\overline{\mathbb{Q}}\mathsf{E}_{G,\psi}\simeq \overline{\mathbb{Q}}^{(G^{\ast F^\ast})_{\mathrm{ss}}/\sim}$ (\S\,\ref{E-decomposition}) as well as a similar identification $\overline{\mathbb{Q}}S^F=\overline{\mathbb{Q}}^{S^{\ast F^\ast}}$ (such that $e_{\widehat{z}}^S$ corresponds to the characteristic function $\mathbf{1}_{\{z\}}$ on $S^{\ast F^\ast}$), we have the following commutative diagram of $\overline{\mathbb{Q}}$-algebras: (``\:{\normalfont Res}" means the restriction map)
\[
\begin{tikzcd}
\overline{\mathbb{Q}}\mathsf{E}_{G,\psi} \arrow[leftrightarrow,r,"\sim"] \arrow[d,"\mathrm{Cur}_S^G"'] & \overline{\mathbb{Q}}^{(G^{\ast F^\ast})_{\mathrm{ss}}/\sim} \arrow[d,"\mathrm{Res}"] \\
\overline{\mathbb{Q}}S^F \arrow[leftrightarrow,r,"\sim"] & \overline{\mathbb{Q}}^{S^{\ast F^\ast}}
\end{tikzcd}
\]
}

The formula (\ref{Curtis-eq1}) gives an equivalent definition of $\mathrm{Cur}_S^G$: the unique $\overline{\mathbb{Q}}$-linear map from $\overline{\mathbb{Q}}\mathsf{E}_{G,\psi}$ to $\overline{\mathbb{Q}}S^F$ satisfying (\ref{Curtis-eq1}) makes the diagram (\ref{Curtis-diagram}) commute. 

{\it Proof of lemma.} As we have fixed a $\psi\in\Psi$, let us drop the subscript $\psi$ here. To prove this lemma, it suffices to establish the equality $\widehat{y}(\mathrm{Cur}_S^G(e_{[x]}^{\mathsf{E}}))=\widehat{y}\left(\sum\limits_{z\in [x]\cap S^{\ast F^\ast}}e_{\widehat{z}}^S\right)$ for every $[x]\in(G^{\ast F^\ast})_{\mathrm{ss}}/\!\sim$ and every $y\in S^{\ast F^\ast}$. By the defining property of $\mathrm{Cur}_S^G$, we have $\widehat{y}(\mathrm{Cur}_S^G(e_{[x]}^{\mathsf{E}}))=\chi_{[y]}(e_{[x]}^{\mathsf{E}})$ (where $[y]$ denotes the image of $y$ in $(G^{\ast F^\ast})_{\mathrm{ss}}/\!\sim$). The orthogonality of irreducible characters tell us that $ge_{[x]}^G=\chi_{[x]}(g)e_{[x]}^G$ for every $g\in G^F$, so $ge_{[x]}^{\mathsf{E}}=\chi_{[x]}(g)e_{[x]}^{\mathsf{E}}$ for every $g\in G^F$ and thus $\chi_{[x]}(e_{[y]}^{\mathsf{E}})e_{[x]}^{\mathsf{E}}=e_{[y]}^{\mathsf{E}}e_{[x]}^{\mathsf{E}}=\delta_{[x][y]}e_{[x]}^{\mathsf{E}}$ (here $\delta_{[x][y]}=1$ if $[x]=[y]$, and $=0$ if $[x]\neq [y]$). Therefore $\chi_{[y]}(e_{[x]}^{\mathsf{E}})=\delta_{[x][y]}$. In the same way, one can show that $\widehat{y}\left(\sum\limits_{z\in [x]\cap S^{\ast F^\ast}}e_{\widehat{z}}^S\right)=\delta_{[x][y]}$. \qed

\end{context}

\begin{context}\label{Curtis-injective}
{\sc Corollary} (of \S\,\ref{Curtis-idempotent}; compare \cite[Prop.\;3.2]{Bonnafe-Kessar}).\;---
{\it 
If for each $w\in W$ we choose an $F$-stable maximal torus $T_w$ of $G$ whose $G^F$-conjugacy class corresponds to $w$ (recall that the $G^F$-conjugacy classes of $F$-stable maximal tori of $G$ are parametrized by elements of $W$; see \S\,\ref{graphical-summary}), then the product map
\[
\mathrm{Cur}^G:=(\mathrm{Cur}_{T_w}^G)_{w\in W}:\overline{\mathbb{Q}}\mathsf{E}_{G,\psi}\longrightarrow\displaystyle\prod_{w\in W}\overline{\mathbb{Q}}(T_w)^F
\]
is an injective $\overline{\mathbb{Q}}$-algebra homomorphism.
}

\end{context}

\begin{context}\label{symmetric-algebra}
{\it Symmetric algebras and symmetrizing forms.}\;\cite[Sec.\;2]{Broue}--- Let $\Lambda$ be a commutative ring with unity. By definition, a symmetric $\Lambda$-algebra is a pair $(A,\tau)$ where $A$ is a $\Lambda$-algebra and $\tau:A\longrightarrow \Lambda$ is a $\Lambda$-module homomorphism such that the conditions (i) and (ii) below are satisfied:
\begin{enumerate}[(i)]
\item considered as a $\Lambda$-module, $A$ is finitely generated and projective;
\item $\tau:A\longrightarrow \Lambda$ is central in the sense that $\tau(ab)=\tau(ba)$ for all $a,b\in A$; the map 
\[
\widehat{\tau}:A\longrightarrow\mathrm{Hom}_\Lambda(A,\Lambda),\quad a\longmapsto [\widehat{\tau}(a):A\longrightarrow \Lambda,\;\;b\longmapsto \tau(ab)]
\]
is a $\Lambda$-module isomorphism. In this case, $\tau$ is called a symmetrizing form on $A$.
\end{enumerate}

Now let $(A,\tau)$ be a symmetric algebra, let $P$ be a finitely generated projective $A$-module and let $E=\mathrm{End}_A(P)$ be its endomorphism ring. Then $\tau$ induces a symmetrizing form $\tau_E:E\longrightarrow \Lambda$ via
$
\tau_E=[E\simeq\mathrm{Hom}_A(P,A)\otimes_A P\xrightarrow{\mathrm{\;natural\;pairing\;}} A\xrightarrow{\;\tau\;}\Lambda],
$
so that $(E,\tau_E)$ is a symmetric $\Lambda$-algebra. 

{\it Application to $\Lambda\mathsf{E}_{G,\psi}$.} Let $\Lambda$ be an integral domain containing $\overline{\mathbb{Z}}[\frac{1}{p}]$. The evaluation at identity $\mathrm{ev}_1:\Lambda G^F\longrightarrow\Lambda$ is a symmetrizing form on $\Lambda G^F$; the Gelfand-Graev representation $\Gamma_{G,\psi}=\Lambda G^Fe_\psi$ being a projective $\Lambda G^F$-module, the previous discussion implies that $\mathrm{ev}_1$ induces a symmetrizing form $\tau$ on $\Lambda\mathsf{E}_{G,\psi}=\mathrm{End}_{\Lambda G^F}(\Gamma_{G,\psi})$; under the identification $\Lambda\mathsf{E}_{G,\psi}=e_\psi\Lambda G^Fe_\psi\subset\Lambda G^F$, the form $\tau:\Lambda\mathsf{E}_{G,\psi}\longrightarrow\Lambda$ is just the restriction of $\mathrm{ev}_1$ on $\Lambda\mathsf{E}_{G,\psi}$.
\end{context}

\begin{context}\label{symmetrizing-lemma}
{\sc Lemma.}\;\cite[Lem.\;3.8]{Bonnafe-Kessar} --- 
{\it
Let $\Lambda$ be an integral domain and let $K$ be its field of fractions. Let $(A,\tau)$ be a symmetric $\Lambda$-algebra, and denote the $K$-linear extension of $\tau$ to $KA$ again by $\tau:KA\longrightarrow K$. Suppose that $A'$ is a subring of $KA$ such that $A\subset A'$ and $\tau(A')\subset\Lambda.$ Then $A=A'$.
}
\end{context}

\begin{context}\label{Curtis-saturation}
{\sc Theorem.}\;\cite[Thm.\;3.7]{Bonnafe-Kessar} --- 
{\it 
Recall the injective Curtis homomorphism
\[
\mathrm{Cur}^G=(\mathrm{Cur}_{T_w}^G)_{w\in W}:\overline{\mathbb{Q}}\mathsf{E}_{G,\psi}\hookrightarrow\displaystyle\prod_{w\in W}\overline{\mathbb{Q}}(T_w)^F
\]
in \S\,\ref{Curtis-injective}. The map $\mathrm{Cur}^G$ is saturated over $\Lambda=\overline{\mathbb{Z}}[\frac{1}{p|W|}]$ (with respect to its field of fractions $\overline{\mathbb{Q}}$) in the sense that 
\[
\mathrm{Cur}^G(\Lambda\mathsf{E}_{G,\psi})=\mathrm{Cur}^G(\overline{\mathbb{Q}}\mathsf{E}_{G,\psi})\cap\left(\prod_{w\in W}\Lambda(T_w)^F\right).
\]
}

{\it Remark.} There are examples of $G$ making $\mathrm{Cur}^G$ not saturated over $\overline{\mathbb{Z}}[\frac{1}{p}]$, such as $G=\mathrm{SL}_2(\overline{\mathbb{F}_q})$ or $\mathrm{GL}_2(\overline{\mathbb{F}_q})$, both of which with $q$ odd (see \cite[Rmk.\;3.9]{Bonnafe-Kessar} and \S\,\ref{GL2-non-saturation}).

For our later use, let us sketch a proof of this theorem. Consider the two $\Lambda$-algebras $A:=\mathrm{Cur}^G(\Lambda\mathsf{E}_{G,\psi})$ and $A':=\mathrm{Cur}^G(\overline{\mathbb{Q}}\mathsf{E}_{G,\psi})\cap\left(\prod_{w\in W}\Lambda(T_w)^F\right)$; to prove that $A=A'$, the idea is to construct a symmetrizing form on $A$ and then apply \S\,\ref{symmetrizing-lemma}. 

We have $A\subset A'$ since each $\mathrm{Cur}_{T_w}^G$ is defined over $\Lambda$ (\S\,\ref{Curtis-integral}). As $\overline{\mathbb{Q}}A=\mathrm{Cur}^G(\overline{\mathbb{Q}}\mathsf{E}_{G,\psi})$, we get the inclusions $A\subset A'\subset\overline{\mathbb{Q}}A$. Consider
\[
\tau^{\mathsf{E}}:=|U^F|\mathrm{ev}_{1_{G^F}}:\Lambda\mathsf{E}_{G,\psi}\longrightarrow\Lambda\quad (\Lambda\mathsf{E}_{G,\psi}=e_\psi\Lambda G^Fe_\psi\subset\Lambda G^F),
\]
so that $(\Lambda\mathsf{E}_{G,\psi},\tau^{\mathsf{E}})$ is a symmetric $\Lambda$-algebra (\S\,\ref{symmetric-algebra}). Denote the $\overline{\mathbb{Q}}$-linear extension of $\tau^{\mathsf{E}}$ again by $
\tau^{\mathsf{E}}:\overline{\mathbb{Q}}\mathsf{E}_{G,\psi}\longrightarrow\overline{\mathbb{Q}}$,
and set 
\[
\tau^W:=\dfrac{1}{|W|}\displaystyle\sum_{w\in W}\mathrm{ev}_{1_{(T_w)^F}}\circ\mathrm{pr}_{\overline{\mathbb{Q}}(T_w)^F}:\prod_{w\in W}\overline{\mathbb{Q}}(T_w)^F\longrightarrow\overline{\mathbb{Q}}
\]
which is a symmetrizing form on $\displaystyle\prod_{w\in W}\overline{\mathbb{Q}}(T_w)^F$ and verifies the relation $\tau^{\mathsf{E}}=\tau^W\circ\mathrm{Cur}^G$ on $\overline{\mathbb{Q}}\mathsf{E}_G$ (see \cite[Sec.\;3.B]{Bonnafe-Kessar}).

The last relation and the injectivity of $\mathrm{Cur}^G$ together imply that $(A,\tau^W|_A)$ is a symmetric $\Lambda$-algebra; on the other hand, since $|W|$ is invertible in $\Lambda$, the map $\tau^W$ is in fact defined over $\Lambda$ and then we have $\tau^W(A')\subset\Lambda$. So \S\,\ref{symmetrizing-lemma} implies that $A=A'$, and this completes the proof of theorem.
\end{context}

\begin{context}\label{Curtis-Galois-corollary}
{\sc Corollary} (of \S\,\ref{Curtis-saturation} and \S\,\ref{Curtis-integral}(b))\;--- 
{
Identify $\overline{\mathbb{Z}}[\frac{1}{p}]\mathsf{E}_G=\overline{\mathbb{Z}}[\frac{1}{p}]\mathsf{E}_{G,\psi}$ (\S\,\ref{lemma-E-compatibility}). Then the map $\mathrm{Cur}^G$ is saturated over $\mathbb{Z}[\frac{1}{p|W|}]$ (with respect to its field of fractions $\mathbb{Q}$):
\[
\mathrm{Cur}^G(\mathbb{Z}[{\textstyle\frac{1}{p|W|}}]\mathsf{E}_{G})=\mathrm{Cur}^G(\mathbb{Q}\mathsf{E}_{G})\cap\left(\prod_{w\in W}\mathbb{Z}[{\textstyle\frac{1}{p|W|}}](T_w)^F\right).
\]
}
\end{context}

\section{Grothendieck groups of representations over defining characteristic}\label{K-section}

\subsection{The Grothendieck group $\mathsf{K}_{G^\ast}$}

\begin{context}\label{K-setup}
{\it The $p$-regular elements.}\;--- Let $(G^{\ast F^\ast})_{p'}$ be the set of $p$-regular elements of $G^{\ast F^\ast}$; recall that an element of $G^{\ast F^\ast}$ is called $p$-regular if its order in $G^{\ast F^\ast}$ is not divided by $p$. We have $(G^{\ast F^\ast})_{p'}=(G^{\ast F^\ast})_{\mathrm{ss}}$ by the Jordan decomposition. Denote also by $(G^{\ast F^\ast})_{p'}/\!\sim$ the set of $G^{\ast F^\ast}$-conjugacy classes in $(G^{\ast F^\ast})_{p'}$; thus $((G^{\ast F^\ast})_{p'}/\!\sim
)=((G^{\ast F^\ast})_{\mathrm{ss}}/\!\sim)$.
\end{context}

\begin{context}\label{K-definition}
{\it Definition of the algebra $\mathsf{K}_{G^\ast}$.}\;--- Consider $\mathrm{Rep}_{\overline{\mathbb{F}_q}}(G^{\ast F^\ast})$, the category of finite-dimensional representations of the finite group $G^{\ast F^\ast}$ over the field $\overline{\mathbb{F}_q}$, or equivalently, the category of finitely generated $\overline{\mathbb{F}_q}G^{\ast F^\ast}$-modules. We define $\mathsf{K}_{G^\ast}$ to be the Grothendieck group of the category $\mathrm{Rep}_{\overline{\mathbb{F}_q}}(G^{\ast F^\ast})$. The tensor product on $\mathrm{Rep}_{\overline{\mathbb{F}_q}}(G^{\ast F^\ast})$ induces a multiplication $\otimes$ on $\mathsf{K}_{G^\ast}$, and then we shall consider $\mathsf{K}_{G^\ast}=(\mathsf{K}_{G^\ast},+,\otimes)$ as a $\mathbb{Z}$-algebra.

Denote by $[\cdot]:\mathrm{Rep}_{\overline{\mathbb{F}_q}}(G^{\ast F^\ast})\longrightarrow\mathsf{K}_{G^\ast}$ the natural map, which descends to a map $[\cdot]$ on the set of isomorphism classes in $\mathrm{Rep}_{\overline{\mathbb{F}_q}}(G^{\ast F^\ast})$. Let $\mathrm{Irr}_{\overline{\mathbb{F}_q}}(G^{\ast F^\ast})$ be the set of isomorphism classes of simple objects in $\mathrm{Rep}_{\overline{\mathbb{F}_q}}(G^{\ast F^\ast})$. Then $\mathsf{K}_{G^\ast}$ is a free $\mathbb{Z}$-module having the set $\{[M]\in\mathsf{K}_{G^\ast}:M\in\mathrm{Irr}_{\overline{\mathbb{F}_q}}(G^{\ast F^\ast})\}$ as a basis (see \cite[Sec.\;14.1]{Serre}).
\end{context}

\begin{context}\label{Brauer-character}
{\it The Brauer character.}\;\cite[Sec.\;18]{Serre} --- Recall that we have fixed an inclusion of multiplicative groups $\kappa:\overline{\mathbb{F}_q}^\times\hookrightarrow\overline{\mathbb{Q}}^\times$ (\S\,\ref{E-setup}). Let $M\in\mathrm{Rep}_{\overline{\mathbb{F}_q}}(G^{\ast F^\ast})$ and $g\in (G^{\ast F^\ast})_{p'}$ (a $p$-regular element, see \S\,\ref{K-setup}). Then the action $g:M\longrightarrow M$ (considered as an $\overline{\mathbb{F}_q}$-linear map) is diagonalizable with eigenvalues $\lambda_1,\cdots,\lambda_N\in\overline{\mathbb{F}_q}^\times$ (here $N:=\dim_{\overline{\mathbb{F}_q}}M$), and we set $(\mathrm{br}\,M)(g):=\kappa(\lambda_1)+\cdots+\kappa(\lambda_N)\in\overline{\mathbb{Q}}$. We thus get a map $\mathrm{br}\,M:(G^{\ast F^\ast})_{p'}\longrightarrow\overline{\mathbb{Q}}$ (the Brauer character of $M$), which descends to an element $\mathrm{br}\,M\in \overline{\mathbb{Q}}^{(G^{\ast F^\ast})_{p'}/\sim}$. 

The map $M\longmapsto \mathrm{br}\,M$ then induces a ring homomorphism
\[
\mathrm{br}:\mathsf{K}_{G^\ast}\longrightarrow \overline{\mathbb{Q}}^{(G^{\ast F^\ast})_{p'}/\sim}
\]
The unique $\overline{\mathbb{Q}}$-linear extension of the map $\mathrm{br}(\cdot)$ is a $\overline{\mathbb{Q}}$-algebra isomorphism:
\[
\mathrm{br}:\overline{\mathbb{Q}}\mathsf{K}_{G^\ast}\xrightarrow{\;\:\sim\;\:} \overline{\mathbb{Q}}^{(G^{\ast F^\ast})_{p'}/\sim}.
\]
Thus the rank of the free $\mathbb{Z}$-module $\mathsf{K}_{G^\ast}$ is
\[
\mathrm{rank}_{\mathbb{Z}}\,\mathsf{K}_{G^\ast}=|\mathrm{Irr}_{\overline{\mathbb{F}_q}}(G^{\ast F^\ast})|=|(G^{\ast F^\ast})_{p'}/\!\sim|=|(G^{\ast F^\ast})_{\mathrm{ss}}/\!\sim|.
\]
The above Brauer isomorphism and the canonical inclusion $\mathsf{K}_{G^\ast}\subset\overline{\mathbb{Q}}\mathsf{K}_{G^\ast}$ show that $\mathsf{K}_{G^\ast}$ is reduced and is injected into $\overline{\mathbb{Q}}^{(G^{\ast F^\ast})_{p'}/\sim}$ by the Brauer map $\mathrm{br}(\cdot)$.

\end{context}

\begin{context}\label{projective-ideal}
{\it Projective objects in $\mathrm{Rep}_{\overline{\mathbb{F}_q}}(G^{\ast F^\ast})$.}\;\cite[III$^{e}$ partie]{Serre} --- Let $\mathrm{Proj}_{\overline{\mathbb{F}_q}}(G^{\ast F^\ast})$ be the category of projective $\overline{\mathbb{F}_q}G^{\ast F^\ast}$-modules of finite dimension over $\overline{\mathbb{F}_q}$, and let $\mathsf{P}_{G^\ast}$ be the Grothendieck group of the category $\mathrm{Proj}_{\overline{\mathbb{F}_q}}(G^{\ast F^\ast})$. For $M\in\mathrm{Rep}_{\overline{\mathbb{F}_q}}(G^{\ast F^\ast})$, denote by $P_M\in \mathrm{Proj}_{\overline{\mathbb{F}_q}}(G^{\ast F^\ast})$ the projective cover of $M$ (which is unique up to isomorphism). Then $\mathsf{P}_{G^\ast}$ is a free $\mathbb{Z}$-module having $\{[P_M]:M\in\mathrm{Irr}_{\overline{\mathbb{F}_q}}(G^{\ast F^\ast})\}$ as a basis. (The image of $P\in \mathrm{Proj}_{\overline{\mathbb{F}_q}}(G^{\ast F^\ast})$ in $\mathsf{P}_{G^\ast}$ is denoted by $[P]$.) Two objects $P$ and $P'$ of $\mathrm{Proj}_{\overline{\mathbb{F}_q}}(G^{\ast F^\ast})$ are isomorphic if and only if $[P]=[P']$ in $\mathsf{P}_{G^\ast}$. Furthermore:
\begin{enumerate}[\normalfont (a)]
\item The pairing
\[
\langle\cdot,\cdot\rangle_{\mathsf{P},\mathsf{K}}:\mathrm{Proj}_{\overline{\mathbb{F}_q}}(G^{\ast F^\ast})\times \mathrm{Rep}_{\overline{\mathbb{F}_q}}(G^{\ast F^\ast})\longrightarrow\mathbb{Z},\quad (P,\pi)\longmapsto \dim_{\overline{\mathbb{F}_q}}\,\mathrm{Hom}_{\overline{\mathbb{F}_q}G^{\ast F^\ast}}(P,\pi)
\]
descends to a $\mathbb{Z}$-bilinear perfect pairing $\langle\cdot,\cdot\rangle_{\mathsf{P},\mathsf{K}}:\mathsf{P}_{G^\ast}\times\mathsf{K}_{G^\ast}\longrightarrow\mathbb{Z}$, with $\{[P_M]:M\in\mathrm{Irr}_{\overline{\mathbb{F}_q}}(G^{\ast F^\ast})\}$ and $\{[M]:M\in\mathrm{Irr}_{\overline{\mathbb{F}_q}}(G^{\ast F^\ast})\}$ being dual bases with each other.
\item Let $c:\mathsf{P}_{G^\ast}\longrightarrow\mathsf{K}_{G^\ast}$ be the Cartan homomorphism, which is the natural map induced by the inclusion $\mathrm{Proj}_{\overline{\mathbb{F}_q}}(G^{\ast F^\ast})\subset \mathrm{Rep}_{\overline{\mathbb{F}_q}}(G^{\ast F^\ast})$. Then $c(\mathsf{P}_{G^\ast})\supset|G^{\ast F^\ast}|_p\mathsf{K}_{G^\ast}$, the cokernel $\mathrm{coker}(c)=\mathsf{K}_{G^\ast}/c(\mathsf{P}_{G^\ast})$ is a finite $p$-group, and $c$ is an injective map. We shall use this map $c$ to identify $\mathsf{P}_{G^\ast}$ as an ideal of the ring $\mathsf{K}_{G^\ast}$.
\end{enumerate}
\end{context}

\subsection{Comparison between $\mathbb{Z}[\frac{1}{p|W|}]\mathsf{K}_{G^\ast}$ and $\mathbb{Z}[\frac{1}{p|W|}]\mathsf{E}_G$}

\begin{context}\label{KE-identification}
{\it The identification $\overline{\mathbb{Q}}\mathsf{E}_G\simeq \overline{\mathbb{Q}}\mathsf{K}_{G^\ast}$.}\;--- From now on, let us fix a $\psi\in\Psi$ and then identify $\overline{\mathbb{Z}}[\frac{1}{p}]\mathsf{E}_G=\overline{\mathbb{Z}}[\frac{1}{p}]\mathsf{E}_{G,\psi}$ (\S\,\ref{lemma-E-compatibility}). We shall then write the idempotent $e_{\psi,[x]}^{\mathsf{E}}$ (see \S\,\ref{E-decomposition}(b)) simply as $e_{[x]}^{\mathsf{E}}$.

Using the Brauer isomorphism in \S\,\ref{Brauer-character}, the equality $(G^{\ast F^\ast})_{p'}=(G^{\ast F^\ast})_{\mathrm{ss}}$ as well as the algebro-geometric interpretation of $\overline{\mathbb{Q}}\mathsf{E}_G$ in \S\,\ref{E-decomposition}, we obtain the following identifications of $\overline{\mathbb{Q}}$-algebras:
\[
\left\lbrace
\begin{array}{ccccccc}
\overline{\mathbb{Q}}\mathsf{E}_{G}&\simeq & \overline{\mathbb{Q}}^{(G^{\ast F^\ast})_{\mathrm{ss}}/\sim}&=& \overline{\mathbb{Q}}^{(G^{\ast F^\ast})_{p'}/\sim}& \simeq &\overline{\mathbb{Q}}\mathsf{K}_{G^\ast}\\
e_{[x]}^{\mathsf{E}}&\leftrightarrow & \mathbf{1}_{\{[x]\}}&=&\mathbf{1}_{\{[x]\}} & \leftrightarrow & \mathrm{br}^{-1}(\mathbf{1}_{\{[x]\}})
\end{array}
\right\rbrace.
\]

When $G=S$ is a torus, the identification $\overline{\mathbb{Q}}\mathsf{E}_S\simeq\overline{\mathbb{Q}}\mathsf{K}_{S^\ast}$ is already true over $\mathbb{Z}$: indeed, we have $\mathsf{E}_S\simeq\mathbb{Z}S^F$, and the chosen duality $S^F\simeq\mathrm{Irr}_{\overline{\mathbb{F}_q}}(S^{\ast F^\ast})$ (\S\,\ref{E-setup}) induces a canonical identification $\mathbb{Z}S^F\simeq \mathsf{K}_{S^\ast}$.

In the general case, we first consider the coefficients $\overline{\mathbb{Z}}[\frac{1}{p}]$ for the sake of the structure theory of $\overline{\mathbb{Z}}[\frac{1}{p}]\mathsf{E}_G$, and it is expected that we still have a $\overline{\mathbb{Z}}[\frac{1}{p}]$-algebra isomorphism $\overline{\mathbb{Z}}[\frac{1}{p}]\mathsf{E}_G\simeq\overline{\mathbb{Z}}[\frac{1}{p}]\mathsf{K}_{G^\ast}$. The first idea is to reconstruct the expected isomorphism from the above case of tori, with the aide of the following functoriality principle:

Let $L$ be an $F$-stable Levi subgroup of $G$ and consider the Curtis homomorphism $\mathrm{Cur}_L^G:\overline{\mathbb{Q}}\mathsf{E}_G\longrightarrow\overline{\mathbb{Q}}\mathsf{E}_L$ defined by $\mathrm{Cur}_L^G(e_{[x]}^{\mathsf{E}_G})=\displaystyle\sum_{\substack{[y]\in(L^{\ast F^\ast})_{\mathrm{ss}}/\sim\\ y\in[x]}}e_{[y]}^{\mathsf{E}_L}$ for all $[x]\in (G^{\ast F^\ast})_{\mathrm{ss}}/\sim$ (compare \cite[Sec.\;2A]{Bonnafe-Kessar} and \S\,\ref{Curtis-idempotent}); as in \S\,\ref{Curtis-geometry}, $\mathrm{Cur}_L^G$ may be interpreted as the restriction map $\overline{\mathbb{Q}}^{(G^{\ast F^\ast})_{\mathrm{ss}}/\sim}\longrightarrow \overline{\mathbb{Q}}^{(L^{\ast F^\ast})_{\mathrm{ss}}/\sim}$; then the isomorphisms $\overline{\mathbb{Q}}\mathsf{E}_G\simeq \overline{\mathbb{Q}}\mathsf{K}_{G^\ast}$ and $\overline{\mathbb{Q}}\mathsf{E}_L\simeq \overline{\mathbb{Q}}\mathsf{K}_{L^\ast}$ constructed above make the following diagram commutative: 
\[
\begin{tikzcd}
\overline{\mathbb{Q}}\mathsf{E}_G\arrow[r,leftrightarrow,"\sim"]\arrow[d,swap,"\mathrm{Cur}_{L}^G"] & \overline{\mathbb{Q}}\mathsf{K}_{G^\ast}\arrow[d,"\mathrm{Res}_{L^{\ast F^\ast}}^{G^{\ast F^\ast}}"]\\
\overline{\mathbb{Q}}\mathsf{E}_L\arrow[r,leftrightarrow,"\sim"] & \overline{\mathbb{Q}}\mathsf{K}_{L^\ast}
\end{tikzcd}
\] 

\end{context}

\begin{context}\label{tau-K}
{\sc Proposition.}\;---
{\it Let $\tau^\mathsf{K}:\overline{\mathbb{Q}}\mathsf{K}_{G^\ast}\longrightarrow\overline{\mathbb{Q}}$ be the $\overline{\mathbb{Q}}$-linear form induced by the symmetrizing form $\tau^\mathsf{E}=|U^F|\mathrm{ev}_{1_{G^F}}:\overline{\mathbb{Q}}\mathsf{E}_G\longrightarrow\overline{\mathbb{Q}}$ (\S\,\ref{Curtis-saturation}) via the canonical identification $\overline{\mathbb{Q}}\mathsf{E}_G\simeq \overline{\mathbb{Q}}\mathsf{K}_{G^\ast}$ in \S\,\ref{KE-identification}. Then, with the perfect pairing $\langle\cdot,\cdot\rangle_{\mathsf{P},\mathsf{K}}$ from \S\,\ref{projective-ideal}(a),
\[
\tau^{\mathsf{K}}(\pi)=\frac{1}{|W|}\sum_{w\in W}\langle\mathbf{1},\mathrm{Res}_{(T_w)^{\ast F^\ast}}^{G^{\ast F^\ast}}\pi\rangle_{\mathsf{P},\mathsf{K}}\quad(\pi\in\mathsf{K}_{G^\ast}).
\]
}

{\it Proof.} Under the canonical identification $\overline{\mathbb{Q}}\mathsf{K}_{G^\ast}=\overline{\mathbb{Q}}\mathsf{E}_G$, $\pi\in\mathsf{K}_{G^\ast}$ corresponds to $\displaystyle\sum_{x\in (G^{\ast F^\ast})_{\mathrm{ss}}/\sim}(\mathrm{br}\,\pi)(x)e_{[x]}^{\mathsf{E}}\in\overline{\mathbb{Q}}\mathsf{E}_G$. With the aide of the formula $\tau^{\mathsf{E}}=\tau^W\circ\mathrm{Cur}^G$ (\S\,\ref{Curtis-saturation}) as well as \S\,\ref{Curtis-idempotent}, we have
\[
\tau^{\mathsf{E}}(e_{[x]}^{\mathsf{E}})=\frac{1}{|W|}\sum_{w\in W}\mathrm{Cur}_{T_w}^G(e_{[x]}^{\mathsf{E}})(1)=\frac{1}{|W|}\sum_{w\in W}\sum_{z\in[x]\cap(T_w)^{\ast F^\ast}}\frac{1}{|(T_w)^F|},
\]
so that
\begin{align*}
\tau^{\mathsf{K}}(\pi)&=\sum_{x\in (G^{\ast F^\ast})_{\mathrm{ss}}/\sim}(\mathrm{br}\,\pi)(x)\cdot\tau^{\mathsf{E}}(e_{[x]}^{\mathsf{E}})
=\frac{1}{|W|}\sum_{x\in (G^{\ast F^\ast})_{\mathrm{ss}}/\sim}\sum_{w\in W}\sum_{z\in[x]\cap(T_w)^{\ast F^\ast}}\frac{(\mathrm{br}\,\pi)(x)}{|(T_w)^F|}\\
&=\frac{1}{|W|}\sum_{w\in W}\frac{1}{{|(T_w)^F|}}\sum_{x\in (G^{\ast F^\ast})_{\mathrm{ss}}/\sim}\sum_{z\in[x]\cap(T_w)^{\ast F^\ast}}(\mathrm{br}\,\pi)(z)\\
&=\frac{1}{|W|}\sum_{w\in W}\frac{1}{{|(T_w)^F|}}\sum_{z\in(T_w)^{\ast F^\ast}}(\mathrm{br}\,\pi)(z)=\frac{1}{|W|}\sum_{w\in W}\langle\mathbf{1},\mathrm{Res}_{(T_w)^{\ast F^\ast}}^{G^{\ast F^\ast}}\pi\rangle_{\mathsf{P},\mathsf{K}}.\tag*{\qed}
\end{align*}
\end{context}

\begin{context}\label{Steinberg-character}
{\it The Steinberg character.}\;---
Let $\mathrm{St}_{G^\ast}\in\mathrm{Rep}_{\overline{\mathbb{Q}_p}}(G^{\ast F^\ast})$ be the Steinberg character of $G^{\ast F^\ast}$, and let $\overline{\mathrm{St}}_{G^\ast}\in \mathrm{Rep}_{\overline{\mathbb{F}_q}}(G^{\ast F^\ast})$ be its image under the mod-$p$ reduction map (the decomposition map) $\mathrm{Rep}_{\overline{\mathbb{Q}_p}}(G^{\ast F^\ast})\longrightarrow \mathrm{Rep}_{\overline{\mathbb{F}_q}}(G^{\ast F^\ast})$. Recall that $\overline{\mathrm{St}}_{G^\ast}$ is projective and is isomorphic to its dual representation $\overline{\mathrm{St}}_{G^\ast}^\vee$. In addition:

\begin{enumerate}[\normalfont (a)]
\item {\normalfont\cite[Cor.\;7.15]{Deligne-Lusztig}} 
In $\mathbb{Z}[\frac{1}{|W|}]\mathsf{K}_{G^\ast}$, we have:
\begin{align*}
[\overline{\mathrm{St}}_{G^\ast}]&=\displaystyle\frac{1}{|W|}\sum_{w\in W}(-1)^{\ell(w)}\left[\mathrm{Ind}_{(T_w)^{\ast F^\ast}}^{G^{\ast F^\ast}}\mathbf{1}\right];\\
[\overline{\mathrm{St}}_{G^\ast}]
\otimes[\overline{\mathrm{St}}_{G^\ast}]&=\frac{1}{|W|}\displaystyle\sum_{w\in W}\left[\mathrm{Ind}_{(T_w)^{\ast F^\ast}}^{G^{\ast F^\ast}}\mathbf{1}\right].
\end{align*}
($\mathbf{1}$ is the trivial representation over $\overline{\mathbb{F}_q}$; $\ell(w)\in\mathbb{N}$ is the length of $w\in W$.)
\item {\normalfont\cite{Lusztig}} The map $(\cdot)\otimes[\overline{\mathrm{St}}_{G^\ast}]:\mathsf{K}_{G^\ast}\longrightarrow\mathsf{P}_{G^\ast}$ is a $\mathbb{Z}$-module isomorphism.
\end{enumerate}
\end{context}

\begin{context}\label{K-symmetrizing}
{\sc Proposition.}\;--- 
{\it 
Let $\tau^\mathsf{K}:\overline{\mathbb{Q}}\mathsf{K}_{G^\ast}\longrightarrow\overline{\mathbb{Q}}$ be the $\overline{\mathbb{Q}}$-linear form in \S\,\ref{tau-K}.  
\begin{enumerate}[\normalfont (a)]
\item For every $\pi\in\mathsf{K}_{G^\ast}$ we have
$
\tau^{\mathsf{K}}(\pi)=\langle[\overline{\mathrm{St}}_{G^\ast}]\otimes [\overline{\mathrm{St}}_{G^\ast}],\pi\rangle_{\mathsf{P},\mathsf{K}}\in\mathbb{Z}.
$
Thus $\tau^{\mathsf{K}}$ is defined over $\mathbb{Z}$, and we shall still denote the restriction of $\tau^{\mathsf{K}}$ to $\mathsf{K}_{G^\ast}$ by $\tau^{\mathsf{K}}:\mathsf{K}_{G^\ast}\longrightarrow\mathbb{Z}$.
\item The $\mathbb{Z}$-linear form $\tau^{\mathsf{K}}:\mathsf{K}_{G^\ast}\longrightarrow\mathbb{Z}$ obtained in (a) induces a $\mathbb{Z}$-bilinear form 
\[
b^{\mathsf{K}}:\mathsf{K}_{G^\ast}\times \mathsf{K}_{G^\ast}\longrightarrow\mathbb{Z},\quad b^{\mathsf{K}}(x,y):=\tau^{\mathsf{K}}(x\otimes y)\quad(x,y\in\mathsf{K}_{G^\ast}).
\] 
The discriminant $\mathrm{disc}\,b^{\mathsf{K}}$ of $b^{\mathsf{K}}$ is then an integer and satisfies
\[\left|\mathrm{disc}\,b^{\mathsf{K}}\right|=\left|\det\left[(\cdot)\otimes[\overline{\mathrm{St}}_{G^\ast}]:\mathsf{K}_{G^\ast}\longrightarrow\mathsf{K}_{G^\ast}\right]\right|=[\mathsf{K}_{G^\ast}:c(\mathsf{P}_{G^\ast})]=p^{m}
\] for some $m\in\mathbb{N}$. (Here $c:\mathsf{P}_{G^\ast}\longrightarrow\mathsf{K}_{G^\ast}$ is the Cartan homomorphism in \S\,\ref{projective-ideal}(b).)
\item $\tau^{\mathsf{K}}:\overline{\mathbb{Z}}[\frac{1}{p}]\mathsf{K}_{G^\ast}\longrightarrow\overline{\mathbb{Z}}[\frac{1}{p}]$ (well-defined by (a)) is a symmetrizing form on $\overline{\mathbb{Z}}[\frac{1}{p}]\mathsf{K}_{G^\ast}$, so $(\overline{\mathbb{Z}}[\frac{1}{p}]\mathsf{K}_{G^\ast},\tau^{\mathsf{K}})$ is a symmetric $\overline{\mathbb{Z}}[\frac{1}{p}]$-algebra.
\end{enumerate}
}

{\it Proof of proposition.} (a) The formula of $\tau^{\mathsf{K}}(\pi)$ follows from \S\,\ref{tau-K}, Frobenius reciprocity and \S\,\ref{Steinberg-character}(a): (we extend $\langle\cdot,\cdot\rangle_{\mathsf{P},\mathsf{K}}$ $\overline{\mathbb{Q}}$-bilinearly)
\begin{align*}
\tau^{\mathsf{K}}(\pi)&=\frac{1}{|W|}\sum_{w\in W}\langle\mathbf{1},\mathrm{Res}_{(T_w)^{\ast F^\ast}}^{G^{\ast F^\ast}}\pi\rangle_{\mathsf{P},\mathsf{K}}
=\frac{1}{|W|}\sum_{w\in W}\langle\mathrm{Ind}_{(T_w)^{\ast F^\ast}}^{G^{\ast F^\ast}}\mathbf{1},\pi\rangle_{\mathsf{P},\mathsf{K}} \\
&=\left\langle\frac{1}{|W|}\sum_{w\in W}\left[\mathrm{Ind}_{(T_w)^{\ast F^\ast}}^{G^{\ast F^\ast}}\mathbf{1}\right],\pi\right\rangle_{\mathsf{P},\mathsf{K}} =\langle[\overline{\mathrm{St}}_{G^\ast}]\otimes [\overline{\mathrm{St}}_{G^\ast}],\pi\rangle_{\mathsf{P},\mathsf{K}} .
\end{align*}

(b) By (a), adjunction and \S\,\ref{Steinberg-character}, for all $x,y\in\mathsf{K}_{G^\ast}$ we have
\[
b^{\mathsf{K}}(x,y)=\tau^{\mathsf{K}}(x\otimes y)=\langle[\overline{\mathrm{St}}_{G^\ast}]\otimes [\overline{\mathrm{St}}_{G^\ast}],x\otimes y\rangle_{\mathsf{P},\mathsf{K}}=
\langle x^\vee\otimes[\overline{\mathrm{St}}_{G^\ast}], y\otimes [\overline{\mathrm{St}}_{G^\ast}]\rangle_{\mathsf{P},\mathsf{K}}.
\]
Consider the basis $\beta:=\{[x]:x\in\mathrm{Irr}_{\overline{\mathbb{F}_q}}(G^{\ast F^\ast})\}$ (resp.\;$P_\beta:=\{[P_x]:x\in\mathrm{Irr}_{\overline{\mathbb{F}_q}}(G^{\ast F^\ast})\}$) for the free $\mathbb{Z}$-module $\mathsf{K}_{G^\ast}$ (resp.\;$\mathsf{P}_{G^\ast}$); see \S\,\ref{projective-ideal}). Using the matrices associated to (bi-)linear forms, we have the decomposition
\[
\left[b^\mathsf{K}\right]_{\beta\times\beta}={}^tX\cdot Y\cdot Z
\]
where:
\begin{align*}
X&:=\left[(\cdot)^\vee\otimes[\overline{\mathrm{St}}_{G^\ast}]:\mathsf{K}_{G^\ast}\longrightarrow\mathsf{P}_{G^\ast}\right]_{\beta}^{P_\beta};\\
Y&:=\left[\langle\cdot,\cdot\rangle_{\mathsf{P},\mathsf{K}}:\mathsf{P}_{G^\ast}\times \mathsf{K}_{G^\ast}\longrightarrow \mathbb{Z}\right]_{P_\beta,\beta};\\
Z&:=\left[(\cdot)\otimes[\overline{\mathrm{St}}_{G^\ast}]:\mathsf{K}_{G^\ast}\longrightarrow\mathsf{K}_{G^\ast}\right]_{\beta}^{\beta}.
\end{align*}
By \S\,\ref{Steinberg-character}(b), $(\cdot)^\vee\otimes[\overline{\mathrm{St}}_{G^\ast}]:\mathsf{K}_{G^\ast}\longrightarrow\mathsf{P}_{G^\ast}$ is invertible, so $\det X\in\mathbb{Z}^\times=\{\pm 1\}$; by \S\,\ref{projective-ideal}(a), $\det Y=1$. We thus obtain:
\[
\mathrm{disc}\, b^{\mathsf{K}}=\det \left(\left[b^\mathsf{K}\right]_{\beta\times\beta}\right)=\pm\det Z=\pm\det\left[(\cdot)\otimes[\overline{\mathrm{St}}_{G^\ast}]:\mathsf{K}_{G^\ast}\longrightarrow\mathsf{K}_{G^\ast}\right].
\]
Recall that for a free $\mathbb{Z}$-module $A$ of finite rank and for a $\mathbb{Z}$-module homomorphism $\varphi:A\longrightarrow A$, we have $|\det\varphi|=[A:\varphi(A)]$ (this can be deduced via the Smith normal form of $\varphi$). Using this fact, \S\,\ref{Steinberg-character}(b) and \S\,\ref{projective-ideal}(b), we get, for some $m\in\mathbb{N}$,
\[
\left|\det\left[(\cdot)\otimes[\overline{\mathrm{St}}_{G^\ast}]:\mathsf{K}_{G^\ast}\longrightarrow\mathsf{K}_{G^\ast}\right]\right|=[\mathsf{K}_{G^\ast}:c(\mathsf{P}_{G^\ast})]=p^m.
\]

(c) By (b), $\mathrm{disc}\, b^{\mathsf{K}}=\pm p^m$ is invertible in $\overline{\mathbb{Z}}[\frac{1}{p}]$, and this is equivalent to saying that $\tau^{\mathsf{K}}:\overline{\mathbb{Z}}[\frac{1}{p}]\mathsf{K}_{G^\ast}\longrightarrow\overline{\mathbb{Z}}[\frac{1}{p}]$ is a symmetrizing form on $\overline{\mathbb{Z}}[\frac{1}{p}]\mathsf{K}_{G^\ast}$. \qed

{\it Remark.} The equality 
\[
\det\left[(\cdot)\otimes[\overline{\mathrm{St}}_{G^\ast}]:\mathsf{K}_{G^\ast}\longrightarrow\mathsf{K}_{G^\ast}\right]=\pm p^m
\]
(a weaker version of (b) but will be sufficient for the proof of \S\,\ref{K-main-theorem}) can be directly verified as follows: the map $(\cdot)\otimes[\overline{\mathrm{St}}_{G^\ast}]:\mathsf{K}_{G^\ast}\longrightarrow\mathsf{K}_{G^\ast}$ corresponds, under the Brauer isomorphism $\mathrm{br}:\overline{\mathbb{Q}}\mathsf{K}_{G^\ast}\xrightarrow{\;\:\sim\;\:} C((G^{\ast F^\ast})_{p'}/\!\sim)_{\overline{\mathbb{Q}}}$ (\S\,\ref{Brauer-character}), to the diagonalizable endomorphism
\[
\overline{\mathbb{Q}}^{(G^{\ast F^\ast})_{p'}/\sim}\longrightarrow \overline{\mathbb{Q}}^{(G^{\ast F^\ast})_{p'}/\sim},\;\;f\longmapsto f\cdot (\mathrm{br}\,\overline{\mathrm{St}}_{G^\ast}),
\]
whose eigenvalues are all of the form $(\mathrm{br}\,\overline{\mathrm{St}}_{G^\ast})(x)=\mathrm{St}_{G^\ast}(x)=\pm|(C_G(x)^{\circ})^F|_p$ with $x\in (G^{\ast F^\ast})_{\mathrm{ss}}$; the desired equality then follows.

\end{context}

\begin{context}\label{K-main-theorem}
{\sc Theorem.}\;--- 
{\it Consider the following commutative diagram of $\overline{\mathbb{Q}}$-algebras introduced in \S\,\ref{KE-identification}:
\begin{equation}\label{K-theorem-diag1}
\begin{tikzcd}
\overline{\mathbb{Q}}\mathsf{E}_G\arrow[r,leftrightarrow,"\sim"]\arrow[d,hookrightarrow,swap,"\mathrm{Cur}^G=\left(\mathrm{Cur}_{T_w}^G\right)_{w\in W}"]&\overline{\mathbb{Q}}\mathsf{K}_{G^\ast}\arrow[d,hookrightarrow,"\mathrm{Res}:=\left(\mathrm{Res}_{(T_w)^{\ast F^\ast}}^{G^{\ast F^\ast}}\right)_{w\in W}"]\\
\displaystyle\prod_{w\in W}\overline{\mathbb{Q}}(T_w)^F\arrow[r,leftrightarrow,"\sim"]&\displaystyle\prod_{w\in W}\overline{\mathbb{Q}}\mathsf{K}_{(T_w)^\ast}
\end{tikzcd}
\end{equation}
Let $\Lambda=\overline{\mathbb{Z}}[\frac{1}{p|W|}]$. Then the above diagram induces by restriction the following commutative diagram of $\Lambda$-algebras:
\[
\begin{tikzcd}
\Lambda\mathsf{E}_G\arrow[r,leftrightarrow,"\sim"]\arrow[d,hookrightarrow,swap,"\mathrm{Cur}^G"]&\Lambda\mathsf{K}_{G^\ast}\arrow[d,hookrightarrow,"\mathrm{Res}"]\\
\displaystyle\prod_{w\in W}\Lambda(T_w)^F\arrow[r,leftrightarrow,"\sim"]&\displaystyle\prod_{w\in W}\Lambda\mathsf{K}_{(T_w)^\ast}
\end{tikzcd}
\]
(By \S\,\ref{Curtis-integral}, $\mathrm{Cur}^G$ is defined over $\Lambda$.) Moreover, the map $\mathrm{Res}:\overline{\mathbb{Q}}\mathsf{K}_{G^\ast}\longrightarrow \displaystyle\prod_{w\in W}\overline{\mathbb{Q}}\mathsf{K}_{(T_w)^\ast}$ is saturated over $\Lambda$:
$
\mathrm{Res}(\Lambda\mathsf{K}_{G^\ast})=\displaystyle\mathrm{Res}(\overline{\mathbb{Q}}\mathsf{K}_{G^\ast})\cap\left(\prod_{w\in W}\Lambda\mathsf{K}_{(T_w)^\ast}\right)
$ (compare \S\,\ref{Curtis-saturation}).
}

In order to display the symmetry between the $\mathsf{E}$-side and the $\mathsf{K}$-side, we shall give two proofs of this theorem, in both of which the main idea is to use the symmetrizing form lemma \S\,\ref{symmetrizing-lemma}.

{\it First proof.} We first show that the map $\mathrm{Res}:\overline{\mathbb{Q}}\mathsf{K}_{G^\ast}\longrightarrow \displaystyle\prod_{w\in W}\overline{\mathbb{Q}}\mathsf{K}_{(T_w)^\ast}$ is saturated over $\Lambda$. So consider the map
\[
i:=\frac{1}{|W|}\sum_{w\in W}(-1)^{\ell(w)}\mathrm{Ind}_{(T_w)^{\ast F^\ast}}^{G^{\ast F^\ast}}:\left(\prod_{w\in W}\Lambda\mathsf{K}_{(T_w)^\ast}\right)\longrightarrow \Lambda\mathsf{K}_{G^\ast},
\]
which is well-defined because $|W|$ is invertible in $\Lambda$. Then the composition
\[
\Lambda\mathsf{K}_{G^\ast}\xrightarrow{\;\mathrm{Res}\;}\left(\prod_{w\in W}\Lambda\mathsf{K}_{(T_w)^\ast}\right)\xrightarrow{\;\;i\;\;}\Lambda\mathsf{K}_{G^\ast} 
\]
coincide with $[\overline{\mathrm{St}}_{G^\ast}]\otimes(\cdot):\Lambda\mathsf{K}_{G^\ast}\longrightarrow\Lambda\mathsf{K}_{G^\ast}$, since for every $\pi\in\Lambda\mathsf{K}_{G^\ast}$ we have 
\begin{align*}
(i\circ\mathrm{Res})(\pi)&=\frac{1}{|W|}\sum_{w\in W}(-1)^{\ell(w)}\mathrm{Ind}_{(T_w)^{\ast F^\ast}}^{G^{\ast F^\ast}}\mathrm{Res}_{(T_w)^{\ast F^\ast}}^{G^{\ast F^\ast}}\pi\\
&=\frac{1}{|W|}\sum_{w\in W}(-1)^{\ell(w)}\left(\mathrm{Ind}_{(T_w)^{\ast F^\ast}}^{G^{\ast F^\ast}}\mathbf{1}\right)\otimes\pi=[\overline{\mathrm{St}}_{G^\ast}]\otimes\pi \quad(\mbox{\S\,\ref{Steinberg-character}(a)}).
\end{align*}
\S\,\ref{K-symmetrizing}(b) (or the remark of \S\,\ref{K-symmetrizing}) then implies: for some $m\in\mathbb{N}$,
\[
\det(i\circ\mathrm{Res})=\det\left[[\overline{\mathrm{St}}_{G^\ast}]\otimes(\cdot):\overline{\mathbb{Z}}[{\textstyle\frac{1}{p}}]\mathsf{K}_{G^\ast}\longrightarrow\overline{\mathbb{Z}}[{\textstyle\frac{1}{p}}]\mathsf{K}_{G^\ast}\right]=\pm p^m;
\]
in particular $\det(i\circ\mathrm{Res})\in\Lambda^\times$. Thus $\mathrm{Res}:\overline{\mathbb{Q}}\mathsf{K}_{G^\ast}\longrightarrow \displaystyle\prod_{w\in W}\overline{\mathbb{Q}}\mathsf{K}_{(T_w)^\ast}$ is saturated over $\Lambda$.

The identifications in \S\,\ref{KE-identification} induce the following commutative diagram: 
\[
\begin{tikzcd}
\Lambda\mathsf{E}_G\arrow[rr,hookrightarrow]\arrow[d,hookrightarrow,swap,"\mathrm{Cur}^G"]&&\overline{\mathbb{Q}}\mathsf{K}_{G^\ast}\arrow[d,hookrightarrow,"\mathrm{Res}"]\\
\displaystyle\prod_{w\in W}\Lambda(T_w)^F\arrow[r,leftrightarrow,"\sim"]&\displaystyle\prod_{w\in W}\Lambda\mathsf{K}_{(T_w)^\ast}\arrow[r,hookrightarrow]&\displaystyle\prod_{w\in W}\overline{\mathbb{Q}}\mathsf{K}_{(T_w)^\ast}
\end{tikzcd}
\]
The injective map $\Lambda\mathsf{E}_G\hookrightarrow\overline{\mathbb{Q}}\mathsf{K}_{G^\ast}$ in this diagram identifies $\Lambda\mathsf{E}_G\subset\overline{\mathbb{Q}}\mathsf{K}_{G^\ast}$. Using the commutativity of the same diagram as well as the fact that the injective map $\mathrm{Res}:\overline{\mathbb{Q}}\mathsf{K}_{G^\ast}\longrightarrow \displaystyle\prod_{w\in W}\overline{\mathbb{Q}}\mathsf{K}_{(T_w)^\ast}$ is saturated, we have $\Lambda\mathsf{E}_G\subset\Lambda\mathsf{K}_{G^\ast}$. Therefore: 
\begin{enumerate}[(i)]
\item we have the inclusions of rings $\Lambda\mathsf{E}_G\subset\Lambda\mathsf{K}_{G^\ast}\subset \overline{\mathbb{Q}}\mathsf{E}_G$;
\item $(\Lambda\mathsf{E}_G,\tau^{\mathsf{E}})$ in the proof of \S\,\ref{Curtis-saturation} is a symmetric $\Lambda$-algebra;
\item $\tau^{\mathsf{E}}(\Lambda\mathsf{K}_{G^\ast})\subset\Lambda$ (\S\,\ref{K-symmetrizing}(a)).
\end{enumerate}
Thus \S\,\ref{symmetrizing-lemma} implies that $\Lambda\mathsf{E}_G=\Lambda\mathsf{K}_{G^\ast}$ (under the identification $ \overline{\mathbb{Q}}\mathsf{E}_G=\overline{\mathbb{Q}}\mathsf{K}_{G^\ast}$). \qed

{\it Second proof.} This time we prove the $\Lambda$-algebra isomorphism $\Lambda\mathsf{E}_G\simeq\Lambda\mathsf{K}_{G^\ast}$ first. The canonical identifications in \S\,\ref{KE-identification} gives the following commutative diagram: 
\[
\begin{tikzcd}
\overline{\mathbb{Q}}\mathsf{E}_G\arrow[rr,hookleftarrow]\arrow[d,hookrightarrow,swap,"\mathrm{Cur}^G"]&&\Lambda\mathsf{K}_{G^\ast}\arrow[d,hookrightarrow,"\mathrm{Res}"]\\
\displaystyle\prod_{w\in W}\overline{\mathbb{Q}}(T_w)^F\arrow[r,hookleftarrow]&\displaystyle\prod_{w\in W}\Lambda(T_w)^F\arrow[r,leftrightarrow,"\sim"]&\displaystyle\prod_{w\in W}\Lambda\mathsf{K}_{(T_w)^\ast}
\end{tikzcd}
\]
We then identify $\Lambda\mathsf{K}_{G^\ast}\subset \overline{\mathbb{Q}}\mathsf{E}_G$. By the commutativity of this diagram and the fact that the map $\mathrm{Cur}^G$ is saturated over $\Lambda$ (\S\,\ref{Curtis-saturation}), we see that $\Lambda\mathsf{K}_{G^\ast}\subset \Lambda\mathsf{E}_G$. Therefore:
\begin{enumerate}[(i)]
\item we have the inclusions of rings $\Lambda\mathsf{K}_{G^\ast}\subset\Lambda\mathsf{E}_G\subset \overline{\mathbb{Q}}\mathsf{K}_{G^\ast}$;
\item $(\Lambda\mathsf{K}_{G^\ast},\tau^{\mathsf{K}})$ is a symmetric $\Lambda$-algebra (\S\,\ref{K-symmetrizing}(c));
\item $\tau^{\mathsf{K}}(\Lambda\mathsf{E}_G)=\tau^{\mathsf{E}}(\Lambda\mathsf{E}_G)\subset\Lambda$.
\end{enumerate}
Thus \S\,\ref{symmetrizing-lemma} implies that $\Lambda\mathsf{E}_G=\Lambda\mathsf{K}_{G^\ast}$ (under the identification $ \overline{\mathbb{Q}}\mathsf{E}_G=\overline{\mathbb{Q}}\mathsf{K}_{G^\ast}$), and we also obtain the commutativity of the second diagram in the statement of the theorem. Using the saturatedness of $\mathrm{Cur}^G$ (\S\,\ref{Curtis-saturation}) and the $\Lambda$-algebra isomorphism $\Lambda\mathsf{E}_G\simeq\Lambda\mathsf{K}_{G^\ast}$, we see that the map $\mathrm{Res}:\overline{\mathbb{Q}}\mathsf{K}_{G^\ast}\longrightarrow \displaystyle\prod_{w\in W}\overline{\mathbb{Q}}\mathsf{K}_{(T_w)^\ast}$ is saturated over $\Lambda$. \qed
\end{context}

\begin{context}\label{K-main-corollary}
{\sc Corollary} (of \S\,\ref{K-main-theorem} and \S\,\ref{Curtis-integral}(b)).\;--- 
{\it 
The commutative diagram (\ref{K-theorem-diag1}) is equivariant under the action of $\mathrm{Gal}(\overline{\mathbb{Q}}/\mathbb{Q})$ and induces by restriction the following commutative diagram of $\mathbb{Z}[\frac{1}{p|W|}]$-algebras:
\[
\begin{tikzcd}
\mathbb{Z}[\frac{1}{p|W|}]\mathsf{E}_G\arrow[r,leftrightarrow,"\sim"]\arrow[d,swap,"\mathrm{Cur}^G"]&\mathbb{Z}[\frac{1}{p|W|}]\mathsf{K}_{G^\ast}\arrow[d,"\mathrm{Res}"]\\
\displaystyle\prod_{w\in W}\mathbb{Z}[{\textstyle\frac{1}{p|W|}}](T_w)^F\arrow[r,leftrightarrow,"\sim"]&\displaystyle\prod_{w\in W}\mathbb{Z}[{\textstyle\frac{1}{p|W|}}]\mathsf{K}_{(T_w)^\ast}
\end{tikzcd}
\]
Moreover, the map $\mathrm{Res}:\mathbb{Q}\mathsf{K}_{G^\ast}\longrightarrow \displaystyle\prod_{w\in W}\mathbb{Q}\mathsf{K}_{(T_w)^\ast}$ is saturated over $\mathbb{Z}[\frac{1}{p|W|}]$:}
\[
\mathrm{Res}(\mathbb{Z}[{\textstyle\frac{1}{p|W|}}]\mathsf{K}_{G^\ast})=\displaystyle\mathrm{Res}(\mathbb{Q}\mathsf{K}_{G^\ast})\cap\left(\prod_{w\in W}\mathbb{Z}[{\textstyle\frac{1}{p|W|}}]\mathsf{K}_{(T_w)^\ast}\right).
\]
\end{context}

\section{Algebras from the invariant theory}\label{B-section}

\subsection{Langlands dual groups}

\begin{context}\label{Langlands-dual}
{\it The Langlands dual.}\;--- Recall the data $(G,T,F)$ and $(X(T),R,Y(T),R^\vee)$ from \S\,\ref{setup-global}. By definition, a Langlands dual of $(G,T)$ is a pair $(G^\vee,T^\vee)$ defined and split over $\mathbb{Z}$ and obtained by assigning its character group (resp.\;cocharacter group, resp.\;set of roots, resp.\;set of coroots) as $Y(T)$ (resp.\;$X(T)$, resp.\;$R^\vee$, resp.\;$R$). From now on, we fix a such Langlands dual $(G^\vee,T^\vee)$ (all choices are isomorphic).

Note that the Deligne-Lusztig dual pair $(G^\ast,T^\ast)$ (\S\,\ref{DL-dual}) may be obtained from $(G^\vee,T^\vee)$ through the base change $\mathrm{Spec}(\mathbb{F}_q)\longrightarrow\mathrm{Spec}(\mathbb{Z})$. We have the identifications $X(T^\vee)=X(T^\ast)=Y(T)$ and $Y(T^\vee)=Y(T^\ast)=X(T)$; moreover, the Weyl groups of $(G^\vee,T^\vee)$ and of $(G^\ast,T^\ast)$ are both identified with the Weyl group $W$ of $(G,T)$.
\end{context}

\begin{context}\label{tau-F}
{\it The automorphism $\tau^\vee$ and the endomorphism $F^\vee$.}\;---
Let $\tau^\ast$ be the automorphism on $X(T^\ast)=\mathrm{Hom}_{\mathrm{alg}}(T^\ast,\mathbb{G}_m)$ induced by the arithmetic Frobenius endomorphisms $\varphi$ on $T^\ast$ and $\mathbb{G}_m$, more precisely $\tau^\ast(\lambda):=\varphi^{-1}\circ\lambda\circ\varphi$ for $\lambda\in X(T^\ast)$; via the identification $X(T^\vee)=X(T^\ast)$ we obtain an automorphism $\tau^\vee$ on $X(T^\vee)$ and hence on $T^\vee$. Using the same identification $X(T^\vee)=X(T^\ast)$, the Frobenius endomorphism $F^\ast$ on $X(T^\ast)$ induces an endomorphism $F^\vee$ on $X(T^\vee)$ and then on $T^\vee$. We have the following properties (compare \cite[Ch.\;3 \& 8]{Digne-Michel}):
\begin{enumerate}[\normalfont (a)]

\item On both $X(T^\vee)$ and $T^\vee$, the endomorphisms $\tau^\vee\circ F^\vee$ and $F^\vee\circ\tau^\vee$ are equal; these endomorphisms are $q$ (the multiplication by $q$) on $X(T^\vee)$ and are $(\cdot)^q$ on $T^\vee$.
\item The following statements are equivalent: (i) $G$ is split over $\mathbb{F}_q$; (ii) $G^\ast$ is split over $\mathbb{F}_q$; (iii) $\tau^\vee=\mathrm{id}$ on $X(T^\vee)$; (iv) $F^\vee=q$ on $X(T^\vee)$; (v) $F^\vee=(\cdot)^q$ on $T^\vee$.
\end{enumerate}

\end{context}
\subsection{The algebra $\mathsf{B}_{G^\vee}$ of the $\mathbb{Z}$-scheme $(T^\vee\sslash W)^{F^\vee}$}

\begin{context}\label{B-definition}
{\it Definition of the algebra $\mathsf{B}_{G^\vee}$}.\;--- As $X(T^\vee)$ is an abelian group of finite rank, the group ring $\mathbb{Z}[X(T^\vee)]$ is a finitely generated commutative $\mathbb{Z}$-algebra, and we know that $T^\vee=\mathrm{Spec}(\mathbb{Z}[X(T^\vee)])$. The Weyl group $W\simeq N_{G^\vee}(T^\vee)/T^\vee$ acts on $T^\vee$ by conjugation and hence on $\mathbb{Z}[X(T^\vee)]$ (by adjoint action). The $\mathbb{Z}$-algebra $\mathbb{Z}[X(T^\vee)]^W$, consisting of elements of $\mathbb{Z}[X(T^\vee)]$ fixed by the $W$-action, is also finitely generated (see \cite[pf.\;of Prop.\;III.18]{Serre2}). We then consider the categorical quotient $T^\vee{\sslash}W:=\mathrm{Spec}(\mathbb{Z}[X(T^\vee)]^W)$, which is an affine $\mathbb{Z}$-scheme. 

The endomorphism $F^\vee$ on $X(T^\vee)$ induces an endomorphism $F^\vee$ on $\mathbb{Z}[X(T^\vee)]^W$ (because $F^\vee(W)=W$) and thus an $F^\vee$-action on $T^\vee{\sslash}W$, so we have the fixed-point subscheme $(T^\vee\sslash W)^{F^\vee}$, which is also an affine $\mathbb{Z}$-scheme. We define $\mathsf{B}_{G^\vee}$ as the ring of functions of the affine scheme $(T^\vee\sslash W)^{F^\vee}$, so
\[
(T^\vee\sslash W)^{F^\vee}=\mathrm{Spec}(\mathsf{B}_{G^\vee})\mbox{\quad and \quad}\mathsf{B}_{G^\vee}=\mathbb{Z}[X(T^\vee)]^W/I
\]
where $I$ is the ideal of $\mathbb{Z}[X(T^\vee)]^W$ generated by the subset $\{F^\vee f-f:f\in\mathbb{Z}[X(T^\vee)]^W\}$; then $\mathsf{B}_{G^\vee}$ is a finitely generated commutative $\mathbb{Z}$-algebra.
\end{context}

\begin{context}\label{B-reduced}
{\it A reducedness problem and the algebra $\mathsf{B}_{G^\vee,\,\mathrm{red}}$.}\;--- In the algebro-geometric viewpoint, we wish to know whether the scheme $(T^\vee\sslash W)^{F^\vee}$ is reduced or not; that is, whether the $\mathbb{Z}$-algebra $\mathsf{B}_{G^\vee}$ is reduced or not. We shall see later (\S\,\ref{q-restricted}) that $\mathsf{B}_{G^\vee}$ is reduced when the derived subgroup $G_{\mathrm{der}}^\vee$ of $G^\vee$ is simply-connected. It is expected that $\mathsf{B}_{G^\vee}$ is always reduced even if $G_{\mathrm{der}}^\vee$ is not simply-connected, but so far we haven't obtained a definitive answer on this (compare \S\,\ref{SO2n-basis}).

At the moment, let us denote by $\mathsf{B}_{G^\vee,\,\mathrm{red}}$ the reduced ring derived from $\mathsf{B}_{G^\vee}$; we have $\mathsf{B}_{G^\vee,\mathrm{red}}=\mathbb{Z}[X(T^\vee)]^W/\sqrt{I}$, where $\sqrt{I}=\{f\in \mathbb{Z}[X(T^\vee)]^W:f^m\in I\;\mathrm{for\;some\;}m\in\mathbb{N}^\ast\}$ is the radical of $I$. So $\mathsf{B}_{G^\vee}$ is reduced if and only if $\mathsf{B}_{G^\vee,\,\mathrm{red}}$ is equal to $\mathsf{B}_{G^\vee}$. Besides, $\mathsf{B}_{G^\vee,\,\mathrm{red}}$ is also a finitely generated commutative $\mathbb{Z}$-algebra (see \S\,\ref{B-definition}).
\end{context}

\begin{context}\label{B-decomposition}
{\it Decompositions of the algebras $k\mathsf{B}_{G^\vee}$ and $k\mathsf{B}_{G^\vee,\,\mathrm{red}}$.}\;--- Let $k$ be an algebraically closed field. The set of $k$-points of $(T^\vee\sslash W)^{F^\vee}=\mathrm{Spec}(\mathsf{B}_{G^\vee})$ may be identified as: 
\[
(T^\vee\sslash W)^{F^\vee}(k)=\mathrm{Specm}(k\mathsf{B}_{G^\vee})=(T^\vee(k)/W)^{F^\vee}=\left(\bigcup_{w\in W}T^\vee(k)^{wF^\vee}\right)\Big\slash W.
\]
(Notation: when a set $X$ is equipped with a (left) $W$-action, we denote by $X/W$ the set of $W$-orbits in $X$.) As each $T^\vee(k)^{w F^\vee}$ $(w\in W)$ is a finite set, we see that $\mathrm{Specm}(k\mathsf{B}_{G^\vee})$ is also a finite set, so $k\mathsf{B}_{G^\vee}$ is a finite-dimensional vector space over $k$ and in particular an artinian $k$-algebra. Therefore, if we denote by $(k\mathsf{B}_{G^\vee})_\mathfrak{m}$ the localization of $k\mathsf{B}_{G^\vee}$ at the maximal ideal $\mathfrak{m}\in (T^\vee{\sslash}W)^{F^\vee}(k)=\mathrm{Specm}(k\mathsf{B}_{G^\vee})$, then the map
\[
k\mathsf{B}_{G^\vee}\longrightarrow\prod_{\mathfrak{m}\in (T^\vee\sslash W)^{F^\vee}(k)}(k\mathsf{B}_{G^\vee})_\mathfrak{m},\quad f\longmapsto (f)_{\mathfrak{m}\in (T^\vee\sslash W)^{F^\vee}(k)}
\]
is a $k$-algebra isomorphism, where each $(k\mathsf{B}_{G^\vee})_\mathfrak{m}\simeq(k\mathsf{B}_{G^\vee})/\mathfrak{m}^N$ as $k$-algebras for some $N\in\mathbb{N}^\ast$ depending on $\mathfrak{m}$. We have the following equivalent conditions:
\begin{align*}
k\mathsf{B}_{G^\vee}\mathrm{\;is\;reduced\;}&\Longleftrightarrow\mathrm{\;each\;}(k\mathsf{B}_{G^\vee})_\mathfrak{m}\simeq (k\mathsf{B}_{G^\vee})/\mathfrak{m}=k \mbox{ as $k$-algebras}\\
&\Longleftrightarrow k\mathsf{B}_{G^\vee}\simeq k^{(T^\vee\sslash W)^{F^\vee}(k)}\mbox{ as $k$-algebras}\\
&\Longleftrightarrow \dim_k (k\mathsf{B}_{G^\vee})=|(T^\vee\sslash W)^{F^\vee}(k)|.
\end{align*}

The above discussion also applies to the reduced version $(T^\vee\sslash W)^{F^\vee}_{\mathrm{red}}=\mathrm{Spec}(\mathsf{B}_{G^\vee,\,\mathrm{red}})$. When $k=\overline{\mathbb{Q}}$ (or other fields of characteristic zero), the reducedness of $\mathsf{B}_{G^\vee,\,\mathrm{red}}$ imply that of $\overline{\mathbb{Q}}\mathsf{B}_{G^\vee,\,\mathrm{red}}$, so that:
\begin{enumerate}[\normalfont (a)]
\item $\overline{\mathbb{Q}}\mathsf{B}_{G^\vee,\,\mathrm{red}}\simeq \overline{\mathbb{Q}}^{(T^\vee\sslash W)^{F^\vee}(\overline{\mathbb{Q}})}$ as $\overline{\mathbb{Q}}$-algebras; 
\item $\dim_{\overline{\mathbb{Q}}} (\overline{\mathbb{Q}}\mathsf{B}_{G^\vee,\,\mathrm{red}})=|(T^\vee\sslash W)^{F^\vee}(\overline{\mathbb{Q}})|$.
\end{enumerate}

{\it Remark.} In general, $k\mathsf{B}_{G^\vee,\,\mathrm{red}}$ need not be reduced. For example, let $G^\vee=T^\vee=\mathbb{G}_m$ with $F^\vee=(\cdot)^q$, so $\mathsf{B}_{G^\vee}=\mathbb{Z}[X^{\pm1}]/(X^{q-1}-1)=\mathsf{B}_{G^\vee,\;\mathrm{red}}$ ($X$ an indeterminate); if $\ell$ is a prime number dividing $q-1$, then $\overline{\mathbb{F}_\ell}\mathsf{B}_{G^\vee,\,\mathrm{red}}=\overline{\mathbb{F}_\ell}[X^{\pm1}]/((X^{\frac{q-1}{\ell}}-1)^\ell)$ is not reduced.

\end{context}

\begin{context}\label{counting-lemma}
{\sc Lemma.}\;---
{\it
Denote by $G^{\ast}_{\mathrm{ss}}/\!\sim$ the set of $G^{\ast}$-conjugacy classes of semisimple elements of $G^{\ast}$, and define $G^{\vee}(\overline{\mathbb{Q}})_{\mathrm{ss}}/\!\sim$ in a similar way. Then}
\[
|(T^\vee{\sslash}W)^{F^\vee}(\overline{\mathbb{Q}})|=|(T^\ast{\sslash}W)^{F^\ast}|=|(G^{\vee}(\overline{\mathbb{Q}})_{\mathrm{ss}}/\!\sim)^{F^\vee}|=|(G^{\ast}_{\mathrm{ss}}/\!\sim)^{F^\ast}|.
\]
{\it (As in our previous convention, we identify $(T^\ast{\sslash}W)^{F^\ast}=(T^\ast{\sslash}W)^{F^\ast}(\overline{\mathbb{F}_q})$.)
}

{\it Proof.} These equalities come from the chosen embedding $\kappa:\overline{\mathbb{F}_q}^\times\hookrightarrow\overline{\mathbb{Q}}^\times$ (\S\,\ref{E-setup}) and the following observations (stated for $G^\ast$ but also valid for $G^\vee(\overline{\mathbb{Q}})$): (i) $G^\ast_{\mathrm{ss}}$ is the union of all maximal tori of $G^\ast$; (ii) for $z\in T^\ast$, $Wz\in (T^\ast{\sslash}W)^{F^\ast}$ if and only if $z\in T^{\ast wF^\ast}$ for some $w\in W$; (iii) elements of each $T^{\ast wF^\ast}$ are of finite order prime to $p$; (iv) The canonical bijection $(G^{\ast}_{\mathrm{ss}}/\!\sim)\simeq T^\ast\sslash W$ induced by the diagonalization is compatible with respect to $F^\ast$ (compare \cite[Sec.\;3.1 \& App.\;B]{Thomas}).  \qed
\end{context}

\begin{context}\label{without-torsion}
{\sc Lemma.}\;--- {\it The $\mathbb{Z}$-module $\mathsf{B}_{G^\vee,\,\mathrm{red}}$ has no nonzero $\mathbb{Z}$-torsion elements. 
}

{\it Proof.} Let $0\neq f\in\mathsf{B}_{G^\vee,\,\mathrm{red}}$. As $\mathsf{B}_{G^\vee,\,\mathrm{red}}$ is a finitely generated $\mathbb{Z}$-algebra (\S\,\ref{B-reduced}), it is a Jacobson ring, so its Jacobson radical (intersection of maximal ideals) is equal to its nilradical (intersection of prime ideals), which is zero because $\mathsf{B}_{G^\vee,\,\mathrm{red}}$ is reduced. Thus there is a maximal ideal $\mathfrak{m}$ of $\mathsf{B}_{G^\vee,\,\mathrm{red}}$ such that $f\not\in\mathfrak{m}$. Consider the field $k=\mathsf{B}_{G^\vee,\,\mathrm{red}}/\mathfrak{m}$, which is a finite field because $\mathsf{B}_{G^\vee,\,\mathrm{red}}$ is a finitely-generated $\mathbb{Z}$-algebra; then consider the canonical quotient map $t:\mathsf{B}_{G^\vee,\,\mathrm{red}}\longrightarrow k$, which represents a $k$-point of $\mathrm{Spec}(\mathsf{B}_{G^\vee,\,\mathrm{red}})$; we then have $f(t)=t(f)\neq 0\in k$. We may and we shall replace $k$ by its algebraic closure, so we shall write $k=\overline{\mathbb{F}_d}$ for some prime number $d$.

Let us show that the $\overline{\mathbb{F}_d}$-point $t:\mathsf{B}_{G^\vee,\,\mathrm{red}}\longrightarrow \overline{\mathbb{F}_d}$ can be canonically lifted to a $\overline{\mathbb{Z}_d}$-point $t':\mathsf{B}_{G^\vee,\,\mathrm{red}}\longrightarrow \overline{\mathbb{Z}_d}$. As $t\in\mathrm{Specm}(\overline{\mathbb{F}_d}\mathsf{B}_{G^\vee,\,\mathrm{red}})=(T^\vee\sslash W)^{F^\vee}(\overline{\mathbb{F}_d})$, we may write $t=Ws$ where $s\in T^\vee(\overline{\mathbb{F}_d})^{wF^\vee}$ for some $w\in W$ (see \S\,\ref{B-decomposition}). Using the canonical lifting $i:\overline{\mathbb{F}_d}^\times \hookrightarrow\overline{\mathbb{Z}_d}^\times$ of $d'$-th roots of unity (Hensel's lemma), our $s\in  T^\vee(\overline{\mathbb{F}_d})=\mathrm{Hom}(X(T^\vee),\overline{\mathbb{F}_d}^\times)$ can be canonically lifted to some $s'\in T^\vee(\overline{\mathbb{Z}_d})=\mathrm{Hom}(X(T^\vee),\overline{\mathbb{Z}_d}^\times)$, and it can be checked that $s'\in T^\vee(\overline{\mathbb{Z}_d})^{wF^\vee}$, so $t':=Ws'\in (T^\vee\sslash W)^{F^\vee}(\overline{\mathbb{Z}_d})$; this $t'$ is a $\overline{\mathbb{Z}_d}$-point $t':\mathsf{B}_{G^\vee,\,\mathrm{red}}\longrightarrow\overline{\mathbb{Z}_d}$ lifting the $\overline{\mathbb{F}_d}$-point $t$ in the sense that the following diagram is commutative (where $r_d$ is the standard reduction map):
\[
\begin{tikzcd}
& \overline{\mathbb{Z}_d}\arrow[twoheadrightarrow,d,"r_d"]\\
\mathsf{B}_{G^\vee,\,\mathrm{red}}\arrow[ru,"t'"]\arrow[r,"t"']&\overline{\mathbb{F}_d}
\end{tikzcd}
\]

Now, for our nonzero element $f\in\mathsf{B}_{G^\vee,\,\mathrm{red}}$ at the beginning of the proof, we have seen that $t(f)\neq 0\in\overline{\mathbb{F}_d}$, so $t'(f)\neq 0\in\overline{\mathbb{Z}_d}$ thanks to the above lifting diagram. Suppose that $n\cdot f=0\in \mathsf{B}_{G^\vee,\,\mathrm{red}}$ for some $n\in\mathbb{Z}$. Then $n\cdot t'(f)= t'(n\cdot f)=0\in\overline{\mathbb{Z}_d}$ while $t'(f)\neq 0$, whence $n=0$. This shows that $f$ is not a $\mathbb{Z}$-torsion element of $\mathsf{B}_{G^\vee,\,\mathrm{red}}$. \qed
\end{context}

\begin{context}\label{B-injection}
{\sc Corollary} (of \S\,\ref{without-torsion}).\;---
{\it The natural map $\mathsf{B}_{G^\vee,\,\mathrm{red}}\longrightarrow\overline{\mathbb{Q}}\mathsf{B}_{G^\vee,\,\mathrm{red}}$ is injective; combining this with the identification $\overline{\mathbb{Q}}\mathsf{B}_{G^\vee,\,\mathrm{red}}= {\overline{\mathbb{Q}}}^{(T^\vee\sslash W)^{F^\vee}(\overline{\mathbb{Q}})}$ in \S\,\ref{B-decomposition}(a), we obtain a canonical injection of $\mathsf{B}_{G^\vee,\,\mathrm{red}}$ into ${\overline{\mathbb{Q}}}^{(T^\vee\sslash W)^{F^\vee}(\overline{\mathbb{Q}})}$.}
\end{context}

\begin{context}\label{B-general}
{\sc Proposition.}\;--- {\it The $\mathbb{Z}$-module $\mathsf{B}_{G^\vee,\,\mathrm{red}}$ is free of rank }
\[
\mathrm{rank}_{\mathbb{Z}}\,\mathsf{B}_{G^\vee,\,\mathrm{red}}=|(T^\vee\sslash W)^{F^\vee}(\overline{\mathbb{Q}})|=|(G^{\ast}_{\mathrm{ss}}/\!\sim)^{F^\ast}|.
\]

{\it Remark.} Having just seen that $\mathsf{B}_{G^\vee,\,\mathrm{red}}$ is without nonzero $\mathbb{Z}$-torsion (\S\,\ref{without-torsion}), the proposition will be proved once we have the finite-generacy of $\mathsf{B}_{G^\vee,\,\mathrm{red}}$ as a $\mathbb{Z}$-module; however, it seems difficult to prove this finite-generacy directly.

{\it Proof of proposition, admitting \S\,\ref{q-restricted}.} Let us admit for the moment the results in \S\,\ref{q-restricted}, which says that this proposition is true when the derived group $G_{\mathrm{der}}^\vee$ of $G^\vee$ is simply-connected. By general theories of algebraic groups, our $G^\vee$ fits into an $F^\vee$-equivariant exact sequence of reductive groups over $\mathbb{Z}$:
\[
1\longrightarrow S^\vee\longrightarrow H^\vee\longrightarrow G^\vee\longrightarrow 1
\]
with $H_{\mathrm{der}}^\vee$ simply-connected and with $S^\vee$ a torus (in particular, the Weyl group of $H^\vee$ is also $W$). Let $T_H^\vee$ denote the maximal torus of $H_{\mathrm{der}}^\vee$ which maps onto $T^\vee$. then the induced morphism of reduced $\mathbb{Z}$-schemes $(T_H^\vee\sslash W)^{F^\vee}_{\mathrm{red}}\longrightarrow (T^\vee\sslash W)^{F^\vee}_{\mathrm{red}}$ is surjective, so that in the level of ring we obtain an inclusion of rings $\mathsf{B}_{G^\vee,\mathrm{red}}\subset\mathsf{B}_{H^\vee,\mathrm{red}}$. As $H_{\mathrm{der}}^\vee$ is simply-connected, \S\,\ref{q-restricted} tells us that $\mathsf{B}_{H^\vee,\mathrm{red}}$ is a free $\mathbb{Z}$-module of finite rank; thus its $\mathbb{Z}$-submodule $\mathsf{B}_{G^\vee,\mathrm{red}}$ is also free of finite rank. This obtained, the rank of $\mathsf{B}_{G^\vee,\mathrm{red}}$ may be calculated by \S\S\,\ref{B-decomposition}\,-\,\ref{counting-lemma}. \qed
\end{context}

\subsection{On the reducedness of $\mathsf{B}_{G^\vee}$}

\begin{context}\label{derived-group}
{\it The derived group $G_{\mathrm{der}}$.}\;(Compare \cite[App.]{Herzig}.)\;--- Denote by $G_{\mathrm{der}}=[G,G]$ the derived subgroup of $G$. We shall freely use the following properties of $G_{\mathrm{der}}$:
\begin{enumerate}[\normalfont (a)]
\item we have $G=R(G).G_{\mathrm{der}}$ with $R(G)\cap G_{\mathrm{der}}$ being a finite set (here $R(G)$ is the reductive radical of $G$; it is the connected identity component of the center of $G$);
\item $T_{\mathrm{der}}:=T\cap G_{\mathrm{der}}$ is a maximal torus of $G_{\mathrm{der}}$, and $\overline{T}:=G/G_{\mathrm{der}}$ is a torus on which the Weyl group $W$ acts trivially;

\item we have a canonical exact sequence of tori $1\longrightarrow T_{\mathrm{der}}\longrightarrow T\longrightarrow\overline{T}\longrightarrow 1$, which induces the following exact sequences of groups:
\begin{align*}
&1\longrightarrow (T_{\mathrm{der}})^{wF}\longrightarrow T^{wF}\longrightarrow\overline{T}^{wF}(=\overline{T}^{F})\longrightarrow 1\quad (w\in W);\\
&0\longrightarrow X(\overline{T})\longrightarrow X(T)\longrightarrow  X(T_{\mathrm{der}})\longrightarrow 0;
\end{align*}
the last exact sequence gives the identifications $X(\overline{T})=X(T)^W=X^0(T)$ where $X^0(T):=\{\lambda\in X(T):\langle\lambda,\alpha^\vee\rangle=0\mbox{ for all }\alpha\in\Delta\}$.
\end{enumerate}
We shall mainly apply these results on the dual sides $G_{\mathrm{der}}^\ast:=(G^\ast)_{\mathrm{der}}$ and $G_{\mathrm{der}}^\vee:=(G^\vee)_{\mathrm{der}}$. Observe that $G_{\mathrm{der}}^\ast$ is simply-connected if and only if $G_{\mathrm{der}}^\vee$ is.
\end{context}

\begin{context}\label{centralizer-theorem}
{\sc Theorem.}\;\cite[pf.\;of Lem.\;3.9]{Steinberg}\cite[Cor.\;8.5]{Steinberg-endo} ---
{\it 
If $G_{\mathrm{der}}$ is simply-connected, then the centralizer $C_G(x)$ of every semisimple element $x$ of $G$ is connected.
}

The citations here address the case where $G$ is simply-connected, while the case of simply-connected $G_{\mathrm{der}}$ can be deduced as a corollary, with the aide of the following observation: for every semisimple element $x=zy\in G$ with $z\in R(G)$ and $y\in G_{\mathrm{der}}$ (\S\,\ref{derived-group}), we have $C_G(x)=R(G).C_{G_{\mathrm{der}}}(y)$.
\end{context}

\begin{context}\label{counting-derived}
{\sc Lemma.}\;---
{\it
Recall that $G^{\ast}_{\mathrm{ss}}/\!\sim$ denotes the set of $G^{\ast}$-conjugacy classes of semisimple elements of $G^{\ast}$ (\S\,\ref{counting-lemma}). Then:
\begin{enumerate}[\normalfont (a)]
\item {\normalfont \cite[Cor.\;3.12]{Digne-Michel}} the map $((G^{\ast F^\ast})_{\mathrm{ss}}/\!\sim)\longrightarrow (G^\ast_{\mathrm{ss}}/\!\sim)^{F^\ast}$ induced by the set inclusion $G^{\ast F^\ast}_{\mathrm{ss}}\subset (G^\ast)_{\mathrm{ss}}$ is surjective.
\item Assume that $G_{\mathrm{der}}^\ast$ is simply-connected. Then the surjection in (a) is a bijection $((G^{\ast F^\ast})_{\mathrm{ss}}/\!\sim)\xrightarrow{\;\;\sim\;\;} (G^\ast_{\mathrm{ss}}/\!\sim)^{F^\ast}$, and we have $|(T^\ast\sslash W)^{F^\ast}|=q^{\mathrm{rank}\,G_{\mathrm{der}}^\ast}\cdot|\overline{T}^{\ast F^\ast}|$.
\end{enumerate}
Similar results hold for the Langlands dual side $(G^\vee(\overline{\mathbb{Q}}),T^\vee(\overline{\mathbb{Q}}))$.}

{\it Proof of (b).} To prove the injectivity of the map $((G^{\ast F^\ast})_{\mathrm{ss}}/\!\sim)\longrightarrow (G^\ast_{\mathrm{ss}}/\!\sim)^{F^\ast}$, suppose that $x,y\in (G^{\ast F^\ast})_{\mathrm{ss}}$ are $G^\ast$-conjugate, say $y=gxg^{-1}$ for some $g\in G^\ast$; then $g^{-1}{F^\ast}(g)\in C_{G^\ast}(x)$; as $G_{\mathrm{der}}^\ast$ is simply-connected, $C_{G^\ast}(x)$ is connected (\S\,\ref{centralizer-theorem}), so the Lang-Steinberg theorem gives us an $h\in C_{G^\ast}(x)$ such that $h^{-1}{F^\ast}(h)=g^{-1}{F^\ast}(g)$; thus $y=gh^{-1}x(gh^{-1})^{-1}$ with $gh^{-1}\in G^{\ast F^\ast}$. 

By \cite[\S\,14.7]{Steinberg-endo}, $|(T^\ast\sslash W)^{F^\ast}|=\displaystyle\frac{1}{|W|}\sum_{w\in W}|T^{\ast wF^\ast}|$. From the exact sequences of finite groups $1\longrightarrow (T_{\mathrm{der}}^{\ast})^{wF^\ast}\longrightarrow T^{\ast wF^\ast}\longrightarrow\overline{T}^{\ast F^\ast}\longrightarrow 1$ for every $w\in W$ (\S\,\ref{derived-group}), we get $|T^{\ast wF^\ast}|=|(T_{\mathrm{der}}^\ast)^{wF^\ast}|\cdot|\overline{T}^{\ast F^\ast}|$ for every $w\in W$, so that the above formula of $|(T^\ast\sslash W)^{F^\ast}|$ gives the equality $|(T^\ast\sslash W)^{F^\ast}|=|(T_{\mathrm{der}}^\ast\sslash W)^{F^\ast}|\cdot|\overline{T}^{\ast F^\ast}|$. As $G_{\mathrm{der}}^\ast$ is simply-connected, we have $|(T_{\mathrm{der}}^\ast\sslash W)^{F^\ast}|=|(G_{\mathrm{der}}^{\ast F^\ast})_{\mathrm{ss}}/\!\sim|=q^{\mathrm{rank}\,G_{\mathrm{der}}^\ast}$: the first equality comes from (a) and \S\,\ref{counting-lemma}, and the second equality has been proved in \cite[Lem.\;7.3]{Steinberg}. \qed
\end{context}

\begin{context}\label{B-combinatorics}
{\it Combinatorics of root data.}\;--- Elements of $X(T^\vee)$ are also called weights; inside $X(T^\vee)$, we shall need the following two subsets:
\begin{align*}
X^+(T^\vee)&:=\{\lambda\in X(T^\vee):\langle\lambda,\alpha^\vee\rangle\geq 0\mbox{ for all }\alpha\in\Delta^\vee\}\quad(\mbox{dominant\;weights});\\
X_q^+(T^\vee)&:=\{\lambda\in X(T^\vee):0\leq\langle\lambda,\alpha^\vee\rangle<q\mbox{ for all }\alpha\in\Delta^\vee\}\quad(q\mbox{-restricted weights}).
\end{align*}

Note that $X^+(T^\vee)$ is identified with the space $X(T^\vee)/W$ of $W$-orbits in $X(T^\vee)$, in the way that every $W$-orbit in $X(T^\vee)$ contains exactly one element of $X^+(T^\vee)$.

When $G^\vee$ is semisimple, $\Delta^\vee$ is a $\mathbb{Q}$-linear basis of $X(T^\vee)_{\mathbb{Q}}:=X(T)\otimes_{\mathbb{Z}}\mathbb{Q}$, from which we introduce two additional notions:
\begin{enumerate}[\normalfont (a)]
\item let $\{\omega_\alpha:\alpha\in\Delta^\vee\}\subset X(T^\vee)_{\mathbb{Q}}$ be the set of fundamental weights of $(G^\vee,T^\vee,B^\vee)$, characterised by the relations $\langle\omega_\alpha,\beta^\vee\rangle=
\left\lbrace
\begin{array}{cc}
1&\mbox{ if }\alpha=\beta\\
0&\mbox{ if }\alpha\neq\beta
\end{array}
\right\rbrace
\mbox{ for all }\alpha,\beta\in\Delta^\vee$;
\item let $\mathrm{ht}:X(T^\vee)\longrightarrow\mathbb{Q}$ be the height function with respect to $\Delta^\vee$, defined for every $\lambda\in X(T^\vee)$ by $\mathrm{ht}(\lambda)=\displaystyle\sum_{\alpha\in\Delta^\vee}m_\alpha$ where $\lambda=\displaystyle\sum_{\alpha\in\Delta^\vee}m_\alpha\alpha$ with all $m_\alpha\in\mathbb{Q}$.
\end{enumerate} 
\end{context}

\begin{context}\label{height-lemma}
{\sc Lemma.}\;--- 
{\it 
Suppose that $G^\vee$ is semisimple, so that we have the height function $\mathrm{ht}:X(T^\vee)\longrightarrow\mathbb{Q}$ with respect to $\Delta^\vee$ (\S\,\ref{B-combinatorics}). Then:
\begin{enumerate}[\normalfont (a)]
\item for every $0\neq\lambda\in X^+(T^\vee)$ and every $\mathrm{id}\neq w\in W$, we have $\mathrm{ht}(w\lambda)<\mathrm{ht}(\lambda)$;
\item for every $0\neq\lambda\in X^+(T^\vee)$, we have $\mathrm{ht}(\lambda)>0$.
\end{enumerate}
}

{\it Proof.} (a) Choose a $W$-invariant inner product $(\cdot|\cdot)$ on the $\mathbb{R}$-vector space $X(T)_{\mathbb{R}}=X(T)\otimes_{\mathbb{Z}}\mathbb{R}$. Set $\widetilde{\rho}:=\displaystyle\sum_{\alpha\in\Delta^\vee}\dfrac{2\omega_\alpha}{(\alpha|\alpha)}\in X(T^\vee)_{\mathbb{R}}$. Then one can show that $\mathrm{ht}(\mu)=(\mu|\widetilde{\rho})$ for every $\mu\in X(T^\vee)$. For every $\mu\in X(T^\vee)$, fix a $\sigma_\mu\in W$ such that $\mathrm{ht}(\sigma_\mu \mu)=\max\limits_{\sigma\in W}\mathrm{ht}(\sigma\mu)$. Then $\sigma_\mu\mu\in X^+(T^\vee)$: for every $\alpha\in\Delta^\vee$, we have $\langle\sigma_\mu\mu,\alpha^\vee\rangle\geq 0$ because
\[
(\sigma_\mu\mu|\widetilde{\rho})\geq (s_\alpha\sigma_\mu\mu|\widetilde{\rho})=(\sigma_\mu\mu|s_\alpha\widetilde{\rho})=\left(\sigma_\mu\mu\middle|\widetilde{\rho}-\frac{2\alpha}{(\alpha|\alpha)}\right)=(\sigma_\mu\mu|\widetilde{\rho})-\langle\sigma_\mu\mu,\alpha^\vee\rangle.
\] 

Now let $0\neq\lambda\in X^+(T^\vee)$ and $\mathrm{id}\neq w\in W$. The previous discussion tells us that $\mathrm{ht}(w\lambda)\leq \mathrm{ht}(\sigma_\lambda\lambda)$ and that $\sigma_\lambda\lambda\in X^+(T^\vee)$. As $W\lambda$ contains exactly one element of $X^+(T^\vee)$, we get $\sigma_\lambda\lambda=\lambda$, and then $\sigma_\lambda=\mathrm{id}$ because $\lambda\neq 0$. Hence $\mathrm{ht}(w\lambda)\leq \mathrm{ht}(\lambda)$. But $\mathrm{ht}(w\lambda)\neq \mathrm{ht}(\lambda)$, for otherwise $w\lambda\in X^+(T^\vee)\cap W\lambda=\{\lambda\}$ and then $w=\mathrm{id}$, contradicting the hypothesis on $w$. Thus $\mathrm{ht}(w\lambda)<\mathrm{ht}(\lambda)$ as desired.

(b) Recall that $\{s_\alpha:\alpha\in\Delta^\vee\}$ determines a length function $W\rightarrow\mathbb{N}$. Let $w_\circ$ be the longest element of $W$ (which maximize the length function); then it is a fact that the action $-w_\circ:X(T^\vee)\longrightarrow X(T^\vee)$ permutes the elements of $\Delta^\vee$, so in particular $\mathrm{ht}(\lambda)=\mathrm{ht}(-w_\circ\lambda)=-\mathrm{ht}(w_\circ\lambda)$ for every $\lambda\in X(T^\vee)$. 

Now let $0\neq\lambda\in X^+(T^\vee)$. Then we have just seen that $\mathrm{ht}(\lambda)=-\mathrm{ht}(w_\circ\lambda)$, while the conclusion of (a) implies that $\mathrm{ht}(w_\circ\lambda)<\mathrm{ht}(\lambda)$ $(w_\circ\neq\mathrm{id})$, from which $\mathrm{ht}(\lambda)>0$.  \qed
\end{context}

\begin{context}\label{invariant-algebra}
{\it The canonical $\mathbb{Z}$-basis of $\mathbb{Z}[X(T^\vee)]^W$.}\;--- For each $\lambda\in X(T^\vee)$, denote by $e(\lambda)$ its image in the group algebra $\mathbb{Z}[X(T^\vee)]$; thus $e(\lambda)$ is identified with the characteristic function $\mathbf{1}_{\{\lambda\}}:X(T^\vee)\longrightarrow\mathbb{Z}$, and $\{e(\lambda):\lambda\in X(T^\vee)\}$ is a $\mathbb{Z}$-linear basis for $\mathbb{Z}[X(T^\vee)]$. For $\lambda\in X(T^\vee)$, let $W\lambda\subset X(T^\vee)$ be the $W$-orbit of $\lambda$, let $W_\lambda\subset W$ be the stabilizer of $\lambda$ under the $W$-action on $X(T^\vee)$, and set $\displaystyle r(\lambda):=\frac{1}{|W_\lambda|}\sum_{w\in W}e(w\lambda)=\sum_{\mu\in W\lambda}e(\mu)$; then each $r(\lambda)$ lies in $\mathbb{Z}[X(T^\vee)]^W$. The identification $ X(T^\vee)/W=X^+(T^\vee)$ implies that $\{r(\lambda):\lambda\in X^+(T^\vee)\}$ is a $\mathbb{Z}$-linear basis for $\mathbb{Z}[X(T^\vee)]^W$.

Observe also that $F^\vee(r(\lambda))=r(F^\vee\lambda)$ for every $\lambda\in X(T^\vee)$ (because $F^\vee(W)=W$).
\end{context}

\begin{context}\label{q-restricted}
{\sc Theorem.}\;--- 
{\it Suppose that $G_{\mathrm{der}}^\vee$ is simply-connected. Let $\omega_\alpha'\in X^+(T_{\mathrm{der}}^\vee)$ ($\alpha\in\Delta^\vee$) be the fundamental weights of $G_{\mathrm{der}}^\vee$, and choose a lifting of $\omega_\alpha'$ to $\omega_\alpha\in X^+(T^\vee)$ ($\alpha\in\Delta^\vee$) via the canonical surjection $\pi:X^+(T^\vee)\twoheadrightarrow X^+(T_{\mathrm{der}}^\vee)$ (so that $\pi(\omega_\alpha)=\omega_\alpha')$. Identify $X(\overline{T}^\vee)\subset X(T^\vee)$ (\S\,\ref{derived-group}), and let $\mathfrak{A}\subset X(\overline{T}^\vee)$ be a set of representatives of the finite $\mathbb{Z}$-module $X(\overline{T}^\vee)/(F^\vee-\mathrm{id})X(\overline{T}^\vee)$. Then $\mathsf{B}_{G^\vee}=\mathbb{Z}[X(T^\vee)]^W/I$ is a free $\mathbb{Z}$-module having the set 
\[
\mathfrak{F}:=\left\lbrace r\left(\mu+\sum_{\alpha\in\Delta^\vee}b_\alpha\omega_\alpha\right)+I:\mu\in\mathfrak{A},\,b_\alpha\in\{0,1,\cdots,q-1\}\,(\alpha\in\Delta^\vee)\right\rbrace,
\]
as its basis, and the rank of $\mathsf{B}_{G^\vee}$ over $\mathbb{Z}$ is 
\[
\mathrm{rank}_{\mathbb{Z}}\,\mathsf{B}_{G^\vee}=|\mathfrak{F}|=q^{\mathrm{rank}\,G_{\mathrm{der}}^\vee}\cdot|\overline{T}^{\vee}\!(\overline{\mathbb{Q}})^{F^\vee}|=|(T^\vee\sslash W)^{F^\vee}(\overline{\mathbb{Q}})|=|(G^{\ast F^\ast})_{\mathrm{ss}}/\!\sim|.
\]
Moreover, $\mathsf{B}_{G^\vee}$ is a reduced ring, so $I=\sqrt{I}$ and $\mathsf{B}_{G^\vee}=\mathsf{B}_{G^\vee,\,\mathrm{red}}$.
}

{\it Remark.} In the special case where $G^\vee$ is simply-connected (so that $G_{\mathrm{der}}^\vee=G^\vee$), we have $\pi=\mathrm{id}$, $\overline{T}^\vee=1$ and $X(\overline{T}^\vee)=0$; we may choose $\omega_\alpha=\omega_\alpha'$, so $\{\omega_\alpha\}_{\alpha\in\Delta^\vee}$ is a $\mathbb{Z}$-basis of $X(T^\vee)$; thus $\mathsf{B}_{G^\vee}$ is a reduced ring and is also a free $\mathbb{Z}$-module having $\mathfrak{F}=\{r(\lambda)+I:\lambda\in X_q^+(T^\vee)\}$ as its $\mathbb{Z}$-linear basis; the $\mathbb{Z}$-rank of $\mathsf{B}_{G^\vee}$ is $q^{\mathrm{rank}\,G^\vee}$.

{\it Proof of theorem.} (Compare \cite[Sec.\;5.6\,-\,5.7]{Humphreys}.) By definition, $G_{\mathrm{der}}^\vee$ is semisimple, so we have the height function $\mathrm{ht}:X(T_{\mathrm{der}}^\vee)\longrightarrow\mathbb{Q}$. Let us still denote by $\mathrm{ht}:X(T^\vee)\longrightarrow\mathbb{Q}$ the composition $\mathrm{ht}\circ\pi:X(T^\vee)\longrightarrow\mathbb{Q}$. For each $f\in\mathbb{Z}[X(T^\vee)]^W$, write $\overline{f}:=f+I\in\mathsf{B}_{G^\vee}$. Then $\{\overline{r(\lambda)}:\lambda\in X^+(T^\vee)\}$ generates $\mathsf{B}_{G^\vee}$ as a $\mathbb{Z}$-module (\S\,\ref{B-definition}). Moreover, in $\mathsf{B}_{G^\vee}$ we have (using \S\,\ref{B-definition}, \S\,\ref{invariant-algebra} and \S\,\ref{tau-F})
\[
\overline{r(F^\vee\lambda)}=\overline{r(\lambda)} \mbox{\;\; and \;\;}\overline{r(q\lambda)}=\overline{r(\tau^\vee\lambda)}\mbox{\;\;for every $\lambda\in X(T^\vee)$}.
\] 

(1) Let us prove first that $\mathsf{B}_{G^\vee}$ is generated by $\{\overline{r(\lambda)}:\lambda\in X_q^+(T^\vee)\}$ as a $\mathbb{Z}$-module. Suppose that $\lambda\in X^+(T^\vee)$ but $\lambda\not\in X_q^+(T^\vee)$. Then there exists an $\alpha\in\Delta^\vee$ such that $\langle\lambda,\alpha^\vee\rangle\geq q$. Set $\lambda':=\lambda-q\omega_\alpha\in X(T^\vee)$; by hypothesis on $\lambda$, we have in fact $\lambda'\in X^+(T^\vee)$. In $\mathbb{Z}[X(T^\vee)]^W$, the product $r(\lambda')r(q\omega_\alpha)$ can be expanded as follows: let $(\cdot)^+:X(T^\vee)\longrightarrow X^+(T^\vee)$ be the composition of the projection $X(T^\vee)\longrightarrow  X(T^\vee)/W$ by the identification $ X(T^\vee)/W\simeq X^+(T^\vee)$, so $W\lambda\cap X^+(T^\vee)=\{\lambda^+\}$ for all $\lambda\in X(T^\vee)$; put $\Omega=\{(\lambda'+\sigma(q\omega_\alpha))^+:\sigma\in W\}\subset X^+(T^\vee)$; then
$\displaystyle r(\lambda')r(q\omega_\alpha)=\sum_{\mu\in\Omega}c_\mu r(\mu)$
for some constants $c_\mu\in\mathbb{Z}\;(\mu\in\Omega)$. Using \S\,\ref{height-lemma}(a), one can show that $c_\lambda=1$ and that $\mathrm{ht}(\mu)<\mathrm{ht}(\lambda)$ for all $\lambda\neq\mu\in\Omega$; thus we obtain: 
\[
\displaystyle r(\lambda')r(q\omega_\alpha)=r(\lambda)+\sum_{\lambda\neq\mu\in\Omega}c_\mu r(\mu)\in\mathbb{Z}[X(T^\vee)]^W.
\]
Passing this relation into $\mathsf{B}_{G^\vee}$ and using the relation $\overline{r(q\omega_\alpha)}=\overline{r(\tau^\vee\omega_\alpha)}\in\mathsf{B}_{G^\vee}$, we get
\[
\overline{r(\lambda')}\,\overline{r(\tau^\vee\omega_\alpha)}=\overline{r(\lambda)}+\sum_{\lambda\neq\mu\in\Omega}c_\mu\overline{r(\mu)}\in\mathsf{B}_{G^\vee}.
\]
As $\tau^\vee\omega_\alpha\in X^+(T^\vee)$ (\S\,\ref{tau-F}), we can expand the product $r(\lambda')r(\tau^\vee\omega_\alpha)$ as in the case $r(\lambda')r(q\omega_\alpha)$ and get
\[
r(\lambda')r(\tau^\vee\omega_\alpha)=r(\lambda'+\tau^\vee\omega_\alpha)+\sum_{\lambda'+\tau^\vee\omega_\alpha\neq\nu\in\Omega'}c_\nu'r(\nu)\in\mathbb{Z}[X(T^\vee)]^W
\]
where $\Omega':=\{(\lambda'+\sigma\tau^\vee\omega_\alpha)^+:\sigma\in W\}$ and $c_\nu'\in\mathbb{Z}$ ($\lambda'+\tau^\vee\omega_\alpha\neq\nu\in\Omega'$). For every $\lambda'+\tau^\vee\omega_\alpha\neq\nu\in\Omega'$, we have again $\mathrm{ht}(\nu)<\mathrm{ht}(\lambda'+\tau^\vee\omega_\alpha)$ (\S\,\ref{height-lemma}(a)), while $\mathrm{ht}(\lambda'+\tau^\vee\omega_\alpha)<\mathrm{ht}(\lambda)$ (\S\,\ref{height-lemma}(b)). Thus we obtain the relation
\[
\overline{r(\lambda)}=\overline{r(\lambda'+\tau^\vee\omega_\alpha)}+\sum_{\lambda'+\tau^\vee\omega_\alpha\neq\nu\in\Omega'}c_\nu'\overline{r(\nu)}-\sum_{\lambda\neq\mu\in\Omega}c_\mu\overline{r(\mu)}\in\mathsf{B}_{G^\vee},
\]
which expresses $\overline{r(\lambda)}$ in terms of a $\mathbb{Z}$-linear combination of some $\overline{r(\gamma)}$ where $\gamma\in X^+(T^\vee)$ with $\mathrm{ht}(\gamma)<\mathrm{ht}(\lambda)$. On the other hand, as $X(T^\vee)$ is a free $\mathbb{Z}$-module of finite rank, we see that $\mathrm{ht}(X^+(T^\vee))\subset h^{-1}\mathbb{N}$ for some $h\in\mathbb{N^\ast}$. Thus we can repeat the above reduction process of $\overline{r(\lambda)}$, so that for every $\lambda\in X^+(T^\vee)$, we can eventually express $\overline{r(\lambda)}$ as a $\mathbb{Z}$-linear combination of those $\overline{r(\mu)}$ with $\mu\in X_q^+(T^\vee)$, as desired.

(2) Let us use (1) to prove that $\mathsf{B}_{G^\vee}$ is generated by $\mathfrak{F}$ as a $\mathbb{Z}$-module. From the canonical exact sequence $0\longrightarrow X(\overline{T}^\vee)\longrightarrow X(T^\vee)\longrightarrow  X(T_{\mathrm{der}}^\vee)\longrightarrow 0$ (\S\,\ref{derived-group}), each $\lambda\in X_q^+(T^\vee)$ may be expressed as $\lambda=\mu+\displaystyle\sum_{\alpha\in\Delta^\vee}b_\alpha\omega_\alpha$ where $\mu\in X(\overline{T}^\vee)$ and each $b_\alpha\in\{0,1,\cdots,q-1\}$. Furthermore, for each $\mu\in X(\overline{T}^\vee)\subset X(T^\vee)$, we have $r(\mu)=e(\mu)\in\mathbb{Z}[X(T^\vee)]^W$ and therefore $e((F^\vee-\mathrm{id})\mu)-1=e(-\mu)(e(F^\vee\mu)-e(\mu))\in I$. These observations and (1) together imply that the $\mathbb{Z}$-module is $\mathsf{B}_{G^\vee}$ is generated by $\mathfrak{F}$.

(3) Let us now prove that $\mathsf{B}_{G^\vee}$ is a free $\mathbb{Z}$-module having $\mathfrak{F}$ as its basis and having the desired rank formulae. In (2) we have seen that $\mathfrak{F}$ generates the $\mathbb{Z}$-module $\mathsf{B}_{G^\vee}$, so that $\mathfrak{F}$ also generates the $\overline{\mathbb{Q}}$-linear vector space $\overline{\mathbb{Q}}\mathsf{B}_{G^\vee}$. By \S\,\ref{B-decomposition}, we have 
\[
\dim_{\overline{\mathbb{Q}}}\overline{\mathbb{Q}}\mathsf{B}_{G^\vee}\geq|(T^\vee\!\sslash W)^{F^\vee}(\overline{\mathbb{Q}})|;
\]
on the other hand, we have (compare the part of duality of tori in \S\,\ref{setup-global})
\[
|\mathfrak{A}|=[X(\overline{T}^\vee):(F^\vee-\mathrm{id})X(\overline{T}^\vee)]=|\mathrm{Irr}_{\overline{\mathbb{Q}}}(\overline{T}^{\vee}\!(\overline{\mathbb{Q}})^{F^\vee})|=|\overline{T}^{\vee}\!(\overline{\mathbb{Q}})^{F^\vee}|,
\]
which implies, together with \S\,\ref{counting-derived} and \S\,\ref{counting-lemma}, that
\[
\dim_{\overline{\mathbb{Q}}}\overline{\mathbb{Q}}\mathsf{B}_{G^\vee}\leq|\mathfrak{F}|=q^{\mathrm{rank}\,G_{\mathrm{der}}^\vee}\cdot|\overline{T}^\vee\!(\overline{\mathbb{Q}})^{F^\vee}|=|(T^\vee\sslash W)^{F^\vee}(\overline{\mathbb{Q}})|=|(G^{\ast F^\ast})_{\mathrm{ss}}/\!\sim|.
\]
Therefore, $\dim_{\overline{\mathbb{Q}}}\overline{\mathbb{Q}}\mathsf{B}_{G^\vee}=|(T^\vee\!\sslash W)^{F^\vee}(\overline{\mathbb{Q}})|$, and $\mathfrak{F}$ is in fact a basis for the $\overline{\mathbb{Q}}$-linear vector space $\overline{\mathbb{Q}}\mathsf{B}_{G^\vee}$; it then follows that $\mathsf{B}_{G^\vee}$ is a free $\mathbb{Z}$-module of basis $\mathfrak{F}$. 

(4) Finally, let us prove the reducedness of $\mathsf{B}_{G^\vee}$. We have just proved in (3) that $\mathsf{B}_{G^\vee}$ is a free $\mathbb{Z}$-module, and this implies that the natual map $\mathsf{B}_{G^\vee}\longrightarrow\overline{\mathbb{Q}}\mathsf{B}_{G^\vee}$ is an inclusion. In (3) we have also shown that $\dim_{\overline{\mathbb{Q}}}\overline{\mathbb{Q}}\mathsf{B}_{G^\vee}=|(T^\vee\!\sslash W)^{F^\vee}(\overline{\mathbb{Q}})|$, so the discussion in \S\,\ref{B-decomposition} implies that $\overline{\mathbb{Q}}\mathsf{B}_{G^\vee}$ is a reduced ring. Hence $\mathsf{B}_{G^\vee}$ is also reduced. \qed
\end{context}

\subsection{Comparison between $\mathsf{B}_{G^\vee}$ and $\mathsf{K}_{G^\ast}$}

\begin{context}\label{algebraic-representation}
{\it Algebraic representations and formal characters.}\;\cite[Ch.\;I.2]{Jantzen} --- Denote by $\mathrm{Rep}_{\mathrm{alg}}(G^\ast)$ the category of $G^\ast$-modules of finite dimension (over the defining field $\overline{\mathbb{F}_q}$ of $G^\ast$); an object of $\mathrm{Rep}_{\mathrm{alg}}(G^\ast)$ is called an algebraic (or rational) representation of $G^\ast$. Let $\mathrm{Irr}_{\mathrm{alg}}(G^\ast)$ be the set of isomorphism classes of simple objects in $\mathrm{Rep}_{\mathrm{alg}}(G^\ast)$. We shall also denote by $\mathsf{K}(\mathrm{Rep}_{\mathrm{alg}}(G^\ast))$ the Grothendieck group of $\mathrm{Rep}_{\mathrm{alg}}(G^\ast)$, which is a ring with multiplication given by the tensor product. 

For $M\in\mathrm{Rep}_{\mathrm{alg}}(G^\ast)$, 
its formal character $\mathrm{ch}\,M$ is defined as follows: $M$ has a weight space decomposition (relative to $T^\ast$) $M=\displaystyle\bigoplus_{\lambda\in X(T^\ast)}M_\lambda$ where the weight spaces $M_\lambda:=\{m\in M:zm=\lambda(z)m\mbox{ for all }z\in T^\ast\}$ (we call $\lambda\in X(T^\ast)$ a weight of $M$ if $M_\lambda\neq\{0\}$); then set $
\mathrm{ch}\,M:=\displaystyle\sum_{\lambda\in X(T^\ast)}\dim_{\overline{\mathbb{F}_q}}M_\lambda\cdot e(\lambda)\in\mathbb{Z}[X(T^\ast)]^W$. (We have $M_{w\lambda}=wM_\lambda$ for every $w\in W$, whence the $W$-invariance of $\mathrm{ch}\,M$.)

\end{context}

\begin{context}\label{highest-weight}
{\sc Lemma.}\;\normalfont{\cite[Ch.\;II.2]{Jantzen}} --- 
{\it 
Identify $X(T^\ast)=X(T^\vee)$, so that for $\lambda\in X(T^\ast)$ we may also define $r(\lambda)\in \mathbb{Z}[X(T^\ast)]^W$ as in \S\,\ref{invariant-algebra}. Let $\leq$ be the standard partial ordering on $X(T^\ast)=X(T^\vee)$ determined by $\Delta^\vee$. Then:
\begin{enumerate}[\normalfont (a)]
\item every $M\in\mathrm{Irr}_{\mathrm{alg}}(G^\ast)$ admits a unique highest weight $\lambda_M$ with respect to the partial ordering $\leq$; we have $\lambda_M\in X^+(T^\ast)$ and $\dim_{\overline{\mathbb{F}_q}}M_{\lambda_M}=1$; also, every weight $\lambda\in X(T^\ast)$ of $M$ satisfies $\lambda\leq\lambda_M$; as a result, for each $M\in \mathrm{Irr}_{\mathrm{alg}}(G^\ast)$, we have $\mathrm{ch}\,M\in r(\lambda_M)+\displaystyle\sum_{\lambda\in X^+(T^\ast),\;\lambda<\lambda_M}\mathbb{Z}.r(\lambda)$;
\item for every $\lambda\in X^+(T^\ast)$, there is a unique $M\in\mathrm{Irr}_{\mathrm{alg}}(G^\ast)$ which admits $\lambda$ as its highest weight in the sense of (a); we shall denote this unique $M$ by $L(\lambda)$;
\item the map $M\longmapsto\mathrm{ch}\,M$ introduced in \S\,\ref{algebraic-representation} induces a $\mathbb{Z}$-algebra isomorphism
\[
\mathrm{ch}:\mathsf{K}(\mathrm{Rep}_{\mathrm{alg}}(G^\ast))\xrightarrow{\;\:\sim\;\:}\mathbb{Z}[X(T^\ast)]^W.
\]
\end{enumerate}
}

\end{context}

\begin{context}\label{Steinberg-twist}
{\sc Lemma.}\;--- 
{\it 
Let $\lambda\in X^+(T^\ast)$. Consider the Frobenius twist $L(\lambda)^{[F^\ast]}$: it is the $G^\ast$-module whose underlying set is $L(\lambda)$ and whose $G^\ast$-action is given by the composition $G^\ast\xrightarrow{\;{F^\ast}\;}G^\ast\xrightarrow{\;\rho\;}\mathrm{GL}(L(\lambda))$ where $\rho$ denotes the $G^\ast$-action on $L(\lambda)$. Then $L(F^\ast\lambda)\simeq L(\lambda)^{[F^\ast]}$ in $\mathrm{Rep}_{\mathrm{alg}}(G^\ast)$.
}

{\it Proof.} The $G^\ast$-module $L(\lambda)^{[F^\ast]}$ is irreducible by \cite[Thm.\;5.1]{Steinberg}; since the highest weight of $L(\lambda)^{[F^\ast]}$ is $F^\ast\lambda$, we see from \S\,\ref{highest-weight} that $L(F^\ast\lambda)\simeq L(\lambda)^{[F^\ast]}$ as $G^\ast$-modules. \qed
\end{context}

\begin{context}\label{Steinberg-restriction}
{\sc Lemma.}\;\cite[Thm.\;7.4]{Steinberg}\cite[Thm.\;3.10]{Herzig} --- 
{\it 
If $G_{\mathrm{der}}^\ast$ is simply-connected, then every $M\in\mathrm{Irr}_{\overline{\mathbb{F}_q}}(G^{\ast F^\ast})$ comes from the restriction of some $L(\lambda)\in\mathrm{Irr}_{\mathrm{alg}}(G^\ast)$ with $\lambda\in X_q^+(T^\ast)$, so that the restriction map $\mathrm{Res}_{G^{\ast F^\ast}}^{G^\ast}:\mathsf{K}(\mathrm{Rep}_{\mathrm{alg}}(G^\ast))\longrightarrow\mathsf{K}_{G^\ast}$ is surjective.
}

The first reference establishes the case where $G^\ast$ is simply-connected, which is generalized in the second reference to the case where $G_{\mathrm{der}}^\ast$ is simply-connected.
\end{context}

\begin{context}\label{KB-proposition}
{\sc Proposition.}\;--- 
{\it 
Identify $X(T^\ast)=X(T^\vee)$ and let $\pi:\mathbb{Z}[X(T^\ast)]^W\twoheadrightarrow\mathsf{B}_{G^\vee,\,\mathrm{red}}$ be the induced quotient map. Observe that the maps
\[
\begin{tikzcd}
(T^\vee\sslash W)^{F^\vee}(\overline{\mathbb{Q}})\arrow[r,"\sim","\S\,\ref{counting-lemma}"']&(G^\ast_{\mathrm{ss}}/\!\sim)^{F^\ast}\arrow[twoheadrightarrow,r,"\S\,\ref{counting-derived}\;"']&((G^{\ast F^\ast})_{\mathrm{ss}}/\!\sim)\arrow[equal,r,"\S\,\ref{K-setup}"']&((G^{\ast F^\ast})_{p'}/\!\sim)
\end{tikzcd}
\]
induce an injection of $\overline{\mathbb{Q}}$-algebras
\[
j:{\overline{\mathbb{Q}}}^{(T^\vee\sslash W)^{F^\vee}(\overline{\mathbb{Q}})} \hookrightarrow {\overline{\mathbb{Q}}}^{(G^{\ast F^\ast})_{p'}/\sim}.
\] 

Then there is an injective $\mathbb{Z}$-algebra homomorphism 
$
\overline{\gamma}:\mathsf{B}_{G^\vee,\,\mathrm{red}}\hookrightarrow\mathsf{K}_{G^\ast}
$
making the following diagram commutative:}
\begin{equation}\label{KB-diagram1}
\begin{tikzcd}
\mathsf{K}(\mathrm{Rep}_{\mathrm{alg}}(G^\ast))\arrow{r}{\mathrm{ch}}[swap]{\sim} \arrow[d,swap,"\mathrm{Res}_{G^{\ast F^\ast}}^{G^\ast}"] & \mathbb{Z}[X(T^\ast)]^W \arrow[twoheadrightarrow,d,"\pi"]\\
\mathsf{K}_{G^\ast}\arrow[hookleftarrow,r,"\overline{\gamma}"]\arrow[hookrightarrow,d,"(\S\,\ref{Brauer-character})\;\mathrm{br}"']&\mathsf{B}_{G^\vee,\,\mathrm{red}}\arrow[hookrightarrow,d,"\S\,\ref{B-injection}"]\\
{\overline{\mathbb{Q}}}^{(G^{\ast F^\ast})_{p'}/\sim}\arrow[hookleftarrow,r,"j"] & {\overline{\mathbb{Q}}}^{(T^\vee\sslash W)^{F^\vee}(\overline{\mathbb{Q}})}
\end{tikzcd}
\end{equation}

{\it Proof.} (1) Let $\lambda\in X^+(T^\ast)$. Then $L(F^\ast\lambda)\simeq L(\lambda)^{[F^\ast]}$ in $\mathrm{Rep}_{\mathrm{alg}}(G^\ast)$ (\S\,\ref{Steinberg-twist}). As $F^\ast$ acts trivially on $G^{\ast F^\ast}$, we have $\mathrm{Res}_{G^{\ast F^\ast}}^{G^\ast}(L(\lambda)^{[F^\ast]})=\mathrm{Res}_{G^{\ast F^\ast}}^{G^\ast}(L(\lambda))$ in $\mathrm{Rep}_{\overline{\mathbb{F}_q}}(G^{\ast F^\ast})$, so $\mathrm{Res}_{G^{\ast F^\ast}}^{G^\ast}L(F^\ast\lambda)=\mathrm{Res}_{G^{\ast F^\ast}}^{G^\ast}L(\lambda)$ in $\mathrm{Rep}_{\overline{\mathbb{F}_q}}(G^{\ast F^\ast})$. Therefore, if we denote by $J$ the ideal of $\mathbb{Z}[X(T^\ast)]^W$ generated by $\{\mathrm{ch}\,L(F^\ast\lambda)-\mathrm{ch}\,L(\lambda):\lambda\in X^+(T^\ast)\}$, then $J$ lies in the kernel of the composition $\gamma:=\mathrm{Res}_{G^{\ast F^\ast}}^{G^\ast}\circ\mathrm{ch}^{-1}:\mathbb{Z}[X(T^\ast)]^W\longrightarrow\mathsf{K}_{G^\ast}$. As $\mathsf{K}_{G^\ast}$ is reduced (\S\,\ref{Brauer-character}), $\gamma$ descends to $\overline{\gamma}:\mathbb{Z}[X(T^\ast)]^W/\sqrt{J}\longrightarrow \mathsf{K}_{G^\ast}$, and we obtain the following commutative diagram of $\mathbb{Z}$-algebras: 
\begin{equation}\label{KB-diagram2}
\begin{tikzcd}
\mathsf{K}(\mathrm{Rep}_{\mathrm{alg}}(G^\ast))\arrow{r}{\mathrm{ch}}[swap]{\sim} \arrow[d,swap,"\mathrm{Res}_{G^{\ast F^\ast}}^{G^\ast}"] & \mathbb{Z}[X(T^\ast)]^W \arrow[twoheadrightarrow,d,"\pi"]\\
\mathsf{K}_{G^\ast}&\mathbb{Z}[X(T^\ast)]^W/\sqrt{J}\arrow[l,"\overline{\gamma}"]
\end{tikzcd}
\end{equation}

(2) We now show that under the identification $\mathbb{Z}[X(T^\ast)]^W=\mathbb{Z}[X(T^\vee)]^W$ (induced by the identification $X(T^\ast)=X(T^\vee)$), the ideal $J\subset\mathbb{Z}[X(T^\ast)]^W$ corresponds to the ideal $I\subset \mathbb{Z}[X(T^\vee)]^W$ appearing in the definition of $\mathsf{B}_{G^\vee,\,\mathrm{red}}=\mathbb{Z}[X(T^\vee)]^W/\sqrt{I}$ (\S\,\ref{B-reduced}). By \S\,\ref{B-definition} and \S\,\ref{invariant-algebra}, the ideal $I$ is generated by $\{r(F^\vee\lambda)-r(\lambda):\lambda\in X^+(T^\vee)\}$. For $\lambda\in X^+(T^\ast)$, by \S\,\ref{highest-weight} we may write $\mathrm{ch}\,L(\lambda)=r(\lambda)+\displaystyle\sum_{\mu\in X^+(T^\ast),\;\mu<\lambda}c_\mu r(\mu)$ where only a finite number of $c_\mu$ is not zero. Then $\mathrm{ch}\,L(F^\ast\lambda)=r(F^\ast\lambda)+\displaystyle\sum_{\mu\in X^+(T^\ast),\;\mu<\lambda}c_\mu r(F^\ast\mu)$: in fact, again by \S\,\ref{Steinberg-twist} we have $L(F^\ast\lambda)\simeq L(\lambda)^{[F^\ast]}$ in $\mathrm{Rep}_{\mathrm{alg}}(G^\ast)$, thus all the weights of $L(F^\ast\lambda)$ is of the form $F^\ast\mu$ where $\mu$ is a weight of $L(\lambda)$, and moreover $L(F^\ast\lambda)_{F^\ast\mu}\simeq L(\lambda)_{\mu}$ for all $\mu\in X(T^\ast)$, whence the desired expression for $\mathrm{ch}\,L(F^\ast\lambda)$. We thus have:
\[
\mathrm{ch}\,L(F^\ast\lambda)-\mathrm{ch}\,L(\lambda)=(r(F^\ast\lambda)-r(\lambda))+\displaystyle\sum_{\mu\in X^+(T^\ast),\;\mu<\lambda}c_\mu. (r(F^\ast\mu)-r(\mu))\quad(\lambda\in X^+(T^\ast)).
\]
This implies that the transition matrix from $\{\mathrm{ch}\,L(F^\ast\lambda)-\mathrm{ch}\,L(\lambda):\lambda\in X^+(T^\ast)\}$ to $\{r(F^\ast\lambda)-r(\lambda):\lambda\in X^+(T^\ast)\}$ is triangular with all diagonal elements being 1 (with respect to the partial ordering $\leq$ on $X^+(T^\ast)$), whence $J=I$. 

(3) The equality $J=I$ established in (2) implies that $\mathbb{Z}[X(T^\ast)]^W/\sqrt{J}=\mathsf{B}_{G^\vee,\,\mathrm{red}}$ under the identification $X(T^\ast)=X(T^\vee)$, so (\ref{KB-diagram2}) is exactly the upper part of (\ref{KB-diagram1}), where we still need to show the injectivity of $\overline{\gamma}$.

(4) To establish (\ref{KB-diagram1}) (in particular the injectivity of $\overline{\gamma}$), it remains to prove that the outermost rectangle diagram in (\ref{KB-diagram1}) commutes, that is to prove the commutativity of the following diagram:
\begin{equation}\label{KB-diagram3}
\quad\quad
\begin{tikzcd}
\mathsf{K}(\mathrm{Rep}_{\mathrm{alg}}(G^\ast))\arrow{r}{\mathrm{ch}}[swap]{\sim} \arrow[d,swap,"\mathrm{br}\circ\mathrm{Res}_{G^{\ast F^\ast}}^{G^\ast}"] & \mathbb{Z}[X(T^\ast)]^W \arrow[twoheadrightarrow,d,"\mathrm{canonical}"]\\
{\overline{\mathbb{Q}}}^{(G^{\ast F^\ast})_{p'}/\sim}\arrow[hookleftarrow,r,"j"] & {\overline{\mathbb{Q}}}^{(T^\ast\sslash W)^{F^\ast}}
\end{tikzcd}
\quad\quad{\color{white}(\ast):}
\end{equation}
(We have identified $(T^\vee\sslash W)^{F^\vee}(\overline{\mathbb{Q}})=(T^\ast\sslash W)^{F^\ast}$.) To do this, we need the help of the following properties: 
\begin{enumerate}[(i)]
\item $(G^{\ast F^\ast})_{p'}=(G^{\ast F^\ast})_{\mathrm{ss}}$ is the union of all $S^{\ast F^\ast}$ where $S^\ast$ runs through elements of $\mathcal{T}^\ast:=\{\mbox{$F^\ast$-stable maximal tori of $G^\ast$}\}$;
\item $(T^\ast\sslash W)^{F^\ast}$ is the set of $W$-orbits in $\displaystyle\bigcup_{w\in W}T^{\ast w F^\ast}$;
\item there is a canonical bijection 
\[
\begin{tikzcd}
\{S^{\ast F^\ast}:S^\ast\in\mathcal{T}^\ast\}\arrow[r,leftrightarrow,"\sim"]&\{T^{\ast wF^\ast}:w\in W\};
\end{tikzcd}
\]
indeed, for each $F^\ast$-stable maximal torus $S^\ast$ of $G^\ast$, there is a $g^\ast\in G^\ast$ such that $g^\ast S^\ast (g^{\ast})^{-1}=T^\ast$, so that the map $x^\ast\longmapsto g^{\ast}x^\ast(g^\ast)^{-1}$ establishes an isomorphism $S^{\ast F^\ast}\xrightarrow{\;\;\sim\;\;}T^{\ast wF^\ast}$ where $w\in W$ is the quotient image of $g^{\ast}F^\ast(g^\ast)^{-1}\in N_{G^\ast}(T^\ast)$.
\end{enumerate}
These properties enable us to integrate (\ref{KB-diagram3}) into the following cubic diagram:
\[
\begin{tikzcd}[row sep=1.8em, column sep=-0.9em]
&\displaystyle\prod_{S^\ast\in\mathcal{T}^\ast}\mathsf{K}(\mathrm{Rep}_{\mathrm{alg}}(S^\ast))\arrow[dd,near start,"\mathrm{br}\circ\mathrm{Res}"] \arrow[rr,"\mathrm{ch}"',"\sim"]&&\displaystyle\prod_{S^\ast\in\mathcal{T}^\ast}\mathbb{Z}[X(S^\ast)] \arrow[dd,"\mathrm{canonical}"]\\
\mathsf{K}(\mathrm{Rep}_{\mathrm{alg}}(G^\ast))\arrow[rr,crossing over,near start,"\mathrm{ch}"',"\sim"] \arrow[dd,swap,"\mathrm{br}\circ\mathrm{Res}_{G^{\ast F^\ast}}^{G^\ast}"] \arrow[ru,"\mathrm{Res}"]&& \mathbb{Z}[X(T^\ast)]^W  \arrow[ru,"\mathrm{(iii)}"']&\\
&\displaystyle\prod_{S^\ast\in\mathcal{T}^\ast} {\overline{\mathbb{Q}}}^{S^{\ast F^\ast}} \arrow[equal,rr]&&\displaystyle\prod_{S^\ast\in\mathcal{T}^\ast} {\overline{\mathbb{Q}}}^{S^{\ast F^\ast}}
\\
{\overline{\mathbb{Q}}}^{(G^{\ast F^\ast})_{p'}/\!\sim}\arrow[hookleftarrow,rr,"j"] \arrow[hookrightarrow,ru,"\mathrm{Res}"]&& {\overline{\mathbb{Q}}}^{(T^\ast\sslash W)^{F^\ast}} \arrow[hookrightarrow,ru,"\mathrm{(iii)}"'] \arrow[twoheadleftarrow,uu,crossing over,near start,"\mathrm{canonical}"]&
\end{tikzcd}
\]
(Here ``Res" means the natural restriction maps, and the maps ``(iii)" on the right face are the natural maps induced by the bijection in (iii) above.) In the above cubic diagram, it can be checked that all the five faces other than the front face (\ref{KB-diagram3}) are commutative diagrams; thus the front face (\ref{KB-diagram3}) is also commutative. \qed

\end{context}

\begin{context}\label{KB-isomorphism}
{\sc Theorem.}\;---
{\it 
If $G_{\mathrm{der}}^\ast$ is simply-connected, then the formal character isomorphism $\mathrm{ch}:\mathsf{K}(\mathrm{Rep}_{\mathrm{alg}}(G^\ast))\xrightarrow{\;\:\sim\;\:}\mathbb{Z}[X(T^\ast)]^W$ (\S\,\ref{highest-weight}) and the bijection of finite varieties $((G^{\ast F^\ast})_{p'}/\!\sim)\simeq (T^\vee\sslash W)^{F^\vee}(\overline{\mathbb{Q}})$ (\S\,\ref{K-setup}, \S\,\ref{counting-lemma} and \S\,\ref{counting-derived}) both induce, via the commutative diagram in \S\,\ref{KB-proposition}, the same $\mathbb{Z}$-algebra isomorphism $ \mathsf{B}_{G^\vee}\simeq \mathsf{K}_{G^\ast}$.  
}

{\it Proof.} The simple-connectedness of $G_{\mathrm{der}}^\ast$ has two consequences: (i) the restriction map $\mathrm{Res}_{G^{\ast F^\ast}}^{G^\ast}:\mathsf{K}(\mathrm{Rep}_{\mathrm{alg}}(G^\ast))\longrightarrow \mathsf{K}_{G^\ast}$ is surjective (\S\,\ref{Steinberg-restriction}), so the injective map $\overline{\gamma}:\mathsf{B}_{G^\vee,\,\mathrm{red}}\longrightarrow\mathsf{K}_{G^\ast}$ in \S\,\ref{KB-proposition} is also surjective and therefore a $\mathbb{Z}$-algebra isomorphism; (ii) $\mathsf{B}_{G^\vee}=\mathsf{B}_{G^\vee,\,\mathrm{red}}$ (\S\,\ref{q-restricted}). Thus the map $\overline{\gamma}$, coming from the formal character map, establishes a $\mathbb{Z}$-algebra isomorphism $ \mathsf{B}_{G^\vee}\simeq \mathsf{K}_{G^\ast}$. \qed
\end{context}

\begin{context}\label{EB-isomorphism}\label{EB-overQ}
{\sc Corollary} (of \S\,\ref{K-main-corollary} and \S\,\ref{KB-isomorphism}).\;---
{\it 
The maps 
\[
\begin{tikzcd}
((G^{\ast F^\ast})_{\mathrm{ss}}/\!\sim)\arrow[equal,r] & ((G^{\ast F^\ast})_{p'}/\!\sim) \arrow[twoheadrightarrow,rr,"\S\,\ref{counting-lemma}","\S\,\ref{counting-derived}"'] &&(T^\vee\sslash W)^{F^\vee}(\overline{\mathbb{Q}})
\end{tikzcd}
\]
induce $\overline{\mathbb{Q}}$-algebra homomorphisms 
\[
\begin{tikzcd}
\overline{\mathbb{Q}}\mathsf{E}_G \arrow[leftrightarrow,r,"\sim"]&\overline{\mathbb{Q}}\mathsf{K}_{G^\ast}\arrow[hookleftarrow,r]&(\overline{\mathbb{Q}}\mathsf{B}_{G^\vee})_{\mathrm{red}}\arrow[equal,r,"\S\,\ref{B-decomposition}"]&\overline{\mathbb{Q}}\mathsf{B}_{G^\vee,\,\mathrm{red}}\arrow[twoheadleftarrow,r]  &  \overline{\mathbb{Q}}\mathsf{B}_{G^\vee}
\end{tikzcd}
\]
which descends to $\mathbb{Z}[\frac{1}{p|W|}]$-algebra homomorphisms:
\begin{equation}\label{identification-eq1}
\begin{tikzcd}
\mathbb{Z}[\frac{1}{p|W|}]\mathsf{E}_G \arrow[leftrightarrow,r,"\sim"]&\mathbb{Z}[\frac{1}{p|W|}]\mathsf{K}_{G^\ast}\arrow[hookleftarrow,r]&\mathbb{Z}[\frac{1}{p|W|}]\mathsf{B}_{G^\vee,\,\mathrm{red}}\arrow[twoheadleftarrow,r]  &  \mathbb{Z}[\frac{1}{p|W|}]\mathsf{B}_{G^\vee}.
\end{tikzcd}
\end{equation}
If $G_{\mathrm{der}}^\ast$ is simply-connected, all maps in (\ref{identification-eq1}) are ${\mathbb{Z}}[\frac{1}{p|W|}]$-algebra isomorphisms: }
\begin{equation}\label{identification-eq2}
{\mathbb{Z}}[{\textstyle\frac{1}{p|W|}}]\mathsf{E}_G\simeq{\mathbb{Z}}[{\textstyle\frac{1}{p|W|}}]\mathsf{K}_{G^\ast}\simeq{\mathbb{Z}}[{\textstyle\frac{1}{p|W|}}]\mathsf{B}_{G^\vee,\,\mathrm{red}}={\mathbb{Z}}[{\textstyle\frac{1}{p|W|}}]\mathsf{B}_{G^\vee}.
\end{equation}
\end{context}

\begin{context}\label{toric-graduation}
{\it Remark on toric graduations.}\;--- Recall from the proof of \S\,\ref{B-general} that a general connected reductive group $G^\vee$ (over $\mathbb{Z}$) fits into an $F^\vee$-equivariant exact sequence of reductive groups over $\mathbb{Z}$,
\[
1\longrightarrow S^\vee\longrightarrow H^\vee\longrightarrow G^\vee\longrightarrow 1,
\]
for a torus $S^\vee$ and a reductive group $H^\vee$ with $H_{\mathrm{der}}^\vee$ simply-connected. Denote by $G^\vee\sslash G^\vee$ be the categorical quotient induced by the adjoint action of $G^\vee$ on itself; note that $G^\vee\sslash G^\vee\simeq T^\vee\sslash W$ by Chevalley's restriction theorem. The multiplication action of $S^{\vee F^\vee}$ on $(H^\vee\sslash H^\vee)^{F^\vee}=\mathrm{Spec}(\mathsf{B}_{H^\vee})$ induces an $S^{\vee F^\vee}$-action on $\mathsf{B}_{H^\vee}$, from which the ring $\mathsf{B}_{H^\vee}$ admits an $X(S^{\vee F^\vee})$-graded decomposition $\mathsf{B}_{H^\vee}=\displaystyle\bigoplus_{\lambda\in X(S^{\vee F^\vee})}(\mathsf{B}_{H^\vee})_\lambda$. Observe that $(\mathsf{B}_{H^\vee})_{0}=(\mathsf{B}_{H^\vee})^{S^{\vee F^\vee}}$, so the canonical surjection $(H^\vee\sslash H^\vee)^{F^\vee}\sslash S^{\vee F^\vee}\twoheadrightarrow (G^\vee\sslash G^\vee)^{F^\vee}$ induces an inclusion of rings $\mathsf{B}_{G^\vee,\,\mathrm{red}}\hookrightarrow(\mathsf{B}_{H^\vee})_0$.

We have analogous discussion on the $\mathsf{K}$-side: as $S^{\ast F^\ast}$ lies in the center of $H^{\ast F^\ast}$, the map associating a character of $H^{\ast F^\ast}$ to its central character induces an $X(S^{\ast F^\ast})$-graded decomposition $\mathsf{K}_{H^\ast}=\displaystyle\bigoplus_{\lambda\in X(S^{\ast F^\ast})}(\mathsf{K}_{H^\ast})_\lambda$, while this time we have canonical ring isomorphisms $(\mathsf{K}_{H^\ast})_0\simeq\mathsf{K}_{G^\ast}$. With the identification $X(S^{\vee F^\vee})\simeq X(S^{\ast F^\ast})$, the ring isomorphism $\mathsf{B}_{H^\vee}\simeq \mathsf{K}_{H^\ast}$ established in \S\,\ref{KB-isomorphism} respects the above graded structures and restricts to a ring isomorphism $(\mathsf{B}_{H^\vee})_0\simeq(\mathsf{K}_{H^\ast})_0$. 

In summary, we have the following commutative diagram of rings:
\[
\begin{tikzcd}
\mathsf{B}_{G^\vee\,\mathrm{red}}\arrow[hookrightarrow,r]\arrow[hookrightarrow,d,"\S\,\ref{KB-proposition}"']&(\mathsf{B}_{H^\vee})_0\arrow[hookrightarrow,r]\arrow[leftrightarrow,d,"\sim"sloped]&\mathsf{B}_{H^\vee}\arrow[leftrightarrow,d,"\sim"'sloped,"\mathrm{graded}\;(\S\,\ref{KB-isomorphism})"]\\
\mathsf{K}_{G^\ast}\arrow[leftrightarrow,r,"\sim"]&(\mathsf{K}_{H^\ast})_0\arrow[hookrightarrow,r]&\mathsf{K}_{H^\ast}
\end{tikzcd}
\]
\end{context}

\section{Further discussions by examples}\label{Example-section}

\subsection{On the identification ${\mathbb{Z}}[\frac{1}{p}]\mathsf{E}_G\simeq{\mathbb{Z}}[\frac{1}{p}]\mathsf{B}_{G^\vee}$}

\begin{context}\label{B-Curtis}
{\it Curtis homomorphisms on the $\mathsf{B}$-side.}\;--- Suppose that $G_{\mathrm{der}}^\ast$ is simply-connected. Then $\mathsf{B}_{G^\vee}=\mathbb{Z}[X(T^\vee)]^W/I$ is a reduced ring (\S\,\ref{q-restricted}), and \S\,\ref{EB-overQ} gives us a $\overline{\mathbb{Q}}$-algebra isomorphism $\overline{\mathbb{Q}}\mathsf{E}_G\simeq\overline{\mathbb{Q}}\mathsf{B}_{G^\vee}$. This isomorphism, the Curtis homomorphism $\mathrm{Cur}^G$ and the canonical isomorphisms $(T_w)^F\simeq T^{wF}$ ($w\in W$) together define the map 
 $\Phi=(\Phi_w)_{w\in W}:\overline{\mathbb{Q}}\mathsf{B}_{G^\vee}\longrightarrow \displaystyle\prod_{w\in W}\overline{\mathbb{Q}}T^{wF}$ in the way that the following diagram is commutative:
\begin{equation}\label{EB-diagram1}
\begin{tikzcd}
\overline{\mathbb{Q}}\mathsf{E}_G\arrow[r,leftrightarrow,"\sim"]\arrow[d,hookrightarrow,swap,"\mathrm{Cur}^G=\left(\mathrm{Cur}_{T_w}^G\right)_{w\in W}"]&\overline{\mathbb{Q}}\mathsf{B}_{G^\vee}\arrow[d,hookrightarrow,"\Phi=(\Phi_{w})_{w\in W}"]\\
\displaystyle\prod_{w\in W}\overline{\mathbb{Q}}(T_w)^F\arrow[r,leftrightarrow,"\sim"]&\displaystyle\prod_{w\in W}\overline{\mathbb{Q}}T^{wF}
\end{tikzcd}
\end{equation}
Note that $\Phi$ is an injective $\overline{\mathbb{Q}}$-algebra homomorphism  as $\mathrm{Cur}^G$ is. 

The diagram (\ref{EB-diagram1}) (``$\mathsf{E}$-$\mathsf{B}$ diagram") is similar to (\ref{K-theorem-diag1}) (``$\mathsf{E}$-$\mathsf{K}$ diagram"), the map $\Phi$ here being the counterpart of the restriction map on the $\mathsf{K}$-side. The $\mathsf{E}$-$\mathsf{K}$ diagram was well-studied in \S\,\ref{K-main-theorem}, however it wouldn't be the case for the $\mathsf{E}$-$\mathsf{B}$ diagram, since a general study of the map $\Phi$ on the $\mathsf{B}$-side would involve complicated combinatorics. 

It is however not difficult to prove that (recall that $\mathsf{B}_{G^\vee}=\mathbb{Z}[X(T^\vee)]^W/I$)
\begin{align*}
\Phi_w(r(\lambda)+I)&=\sum_{t\in T^{wF}}\frac{1}{|W_\lambda|}\left(\sum_{\sigma\in W}\left\langle(\sigma\lambda)|_{T^{\ast wF^\ast}},\widehat{t}\;\right\rangle_{T^{\ast wF^\ast}}\right)t\\
&=\sum_{t\in T^{wF}}\left(\sum_{\sigma\in W_\lambda\backslash W}\left\langle(\sigma\lambda)|_{T^{\ast wF^\ast}},\widehat{t}\;\right\rangle_{T^{\ast wF^\ast}}\right)t
\in\mathbb{Z}T^{wF}
\end{align*}
for each $w\in W$ and for each $\lambda\in X(T^\vee)$ (compare \S\,\ref{tau-K}). In the above formula, for each $t\in T^{wF}$, $\widehat{t}:T^{\ast wF^\ast}\longrightarrow\overline{\mathbb{Q}}^\times$ is its corresponding element under the chosen duality $T^{wF}\simeq\mathrm{Irr}_{\overline{\mathbb{Q}}}(T^{\ast wF^\ast})$; also, each $\lambda\in X(T^\vee)$ is identified with a map $\lambda:T^{\ast}\longrightarrow\overline{\mathbb{Q}}^\times$ via the chosen identifications $T^\vee(\overline{\mathbb{F}_q})=T^\ast$ and $\kappa:\overline{\mathbb{F}_q}^\times\hookrightarrow\overline{\mathbb{Q}}^\times$ (\S\,\ref{E-setup}). As a consequence:
\begin{enumerate}[(a)]
\item the map $\Phi$ is in fact defined over $\mathbb{Z}$:
\begin{equation}\label{Phi-equation}
\Phi=(\Phi_w)_{w\in W}:\mathsf{B}_{G^\vee}\longrightarrow\prod_{w\in W}\mathbb{Z}T^{wF};
\end{equation}
\item the commutative diagram (\ref{EB-diagram1}) is equivariant under the action of $\mathrm{Gal}(\overline{\mathbb{Q}}/\mathbb{Q})$ (using (a) and \S\,\ref{Curtis-integral}(b)) and hence induces by restriction a commutative diagram of $\mathbb{Q}$-algebras, obtained by replacing all $\overline{\mathbb{Q}}$-coefficients in (\ref{EB-diagram1}) by $\mathbb{Q}$-coefficients.

\end{enumerate}

We now study the $\mathsf{E}$-$\mathsf{B}$ diagram for $G=\mathrm{GL}_2(\overline{\mathbb{F}_q})$ and $\mathrm{PGL}_2(\overline{\mathbb{F}_q})$, in which case the combinatorics will be tolerable and we will be able to establish a $\mathbb{Z}[\frac{1}{p}]$-algebra isomorphism $\mathbb{Z}[\frac{1}{p}]\mathsf{E}_G\simeq{\mathbb{Z}}[\frac{1}{p}]\mathsf{B}_{G^\vee}$.
\end{context}

\begin{context}\label{B-GL2-lemma}
{\sc Proposition.}\;--- 
{\it 
Let $G=\mathrm{GL}_2(\overline{\mathbb{F}_q})$, let $T$ be the diagonal maximal torus of $G$, and let $W=\{1,s\}$ be the Weyl group of $G$ (where $s$ is the image of $\left[\begin{array}{cc}0&1\\1&0\end{array}\right]\in N_G(T)$ in $W$), so that the map $\Phi$ in (\ref{Phi-equation}) is 
\[
\Phi=(\Phi_1,\Phi_s):\mathsf{B}_{G^\vee}\longrightarrow\mathbb{Z}T^F\times\mathbb{Z}T^{sF}.
\]
Denote again the linear extensions of this map by $\Phi$. Let $Z$ be the center of $G$, so that $Z^F=T^F\cap T^{sF}$. Let $\Lambda$ be an integral domain and let $K$ be the field of fractions of $\Lambda$.  Then $\Phi$ is saturated over $\Lambda$ (with respect to $K$):}
\[
\Phi({\Lambda}\mathsf{B}_{G^\vee})=\Phi(K\mathsf{B}_{G^\vee})\cap\{(f_1,f_s)\in{\Lambda}T^F\times{\textstyle\Lambda}T^{sF}:f_1|_{Z^F}-f_s|_{Z^F}\in 2{\Lambda}Z^F\}.
\]

{\it Sketch of proof.} Let us prove (b) directly. Note first that $T^\vee$ is also the diagonal maximal torus of $G^\vee=\mathrm{GL}_2$, and $X(T^\vee)=\mathbb{Z}\varepsilon_1\oplus\mathbb{Z}\varepsilon_2$ where 
\[
\varepsilon_i:T^\vee\longrightarrow\mathbb{G}_m,\quad \varepsilon_i(\mathrm{diag}(t_1,t_2)):=t_i\quad(i=1,2).
\] 
By \S\,\ref{q-restricted}, the $\mathbb{Z}$-module $\mathsf{B}_{G^\vee}=\mathbb{Z}[X(T^\vee)]^W/I$ is free of rank $q(q-1)$ and admits a basis 
\[
\mathfrak{F}=\{r_{i,j}:=r((i+j)\varepsilon_1+j\varepsilon_2)+I:i,j\in\mathbb{N},\;0\leq i\leq q-1,\;0\leq j\leq q-2\}.
\]

Let us now calculate the coefficients
\[
\Phi_w(r_{i,j})(t)=\frac{1}{|W_{(i+j)\varepsilon_1+j\varepsilon_2}|}\sum_{\sigma\in W}\left\langle(\sigma((i+j)\varepsilon_1+j\varepsilon_2))|_{T^{\ast wF^\ast}},\widehat{t}\;\right\rangle_{T^{\ast wF^\ast}}\in\mathbb{Z}
\] for $0\leq i\leq q-1$ and $0\leq j\leq q-2$ (see \S\,\ref{B-Curtis}). 

Recall the chosen group isomorphism $\iota:\overline{\mathbb{F}_q}^\times\xrightarrow{\;\sim\;}(\mathbb{Q}/\mathbb{Z})_{p'}$ (\S\,\ref{E-setup}). Set $\xi:=\iota^{-1}(\frac{1}{q^2-1})$ which is a primitive $(q^2-1)$-st root of unity generating the multiplicative group $\mathbb{F}_{q^2}^\times$, and set $\zeta:=\iota^{-1}(\frac{1}{q-1})=\xi^{q+1}$ which is a primitive $(q-1)$-st root of unity generating the multiplicative group $\mathbb{F}_q^\times$. Then we have the following group isomorphisms:
\[
\begin{array}{rr}
\theta_1:\mathbb{Z}/(q-1)\mathbb{Z}\times\mathbb{Z}/(q-1)\mathbb{Z}\xrightarrow{\;\;\sim\;\;}T^F,&\quad (a,b)\longmapsto\theta_1(a,b)=\left[\begin{array}{cc}\zeta^a&0\\0&\zeta^b\end{array}\right];\\
\theta_s:\mathbb{Z}/(q^2-1)\mathbb{Z}\xrightarrow{\;\;\sim\;\;}T^{sF},&\quad c\longmapsto\theta_s(c)=\left[\begin{array}{cc}\xi^c&0\\0&\xi^{qc}\end{array}\right].
\end{array}
\]
 
With these preparations, one can deduce the following calculations of coefficients: (for $m\in\mathbb{Z}$, we write $\overline{m}:=m+(q-1)\mathbb{Z}\in\mathbb{Z}/(q-1)\mathbb{Z}$)
\begin{enumerate}[(i)]
\item $w=1$, $t=\theta_1(a,b)\in T^F$: \\[2mm]
if $a=b\in\mathbb{Z}/(q-1)\mathbb{Z}$:
$
\Phi_1(r_{i,j})(\theta_1(a,a))=
\left\lbrace
\begin{array}{cl}
1 &((i,\overline{j})=(0,a));\\
2 &((i,\overline{j})=(q-1,a));\\
0 &(\mathrm{otherwise});
\end{array}
\right.
$\\[2mm]
if $a\neq b\in\mathbb{Z}/(q-1)\mathbb{Z}$:
$
\Phi_1(r_{i,j})(\theta_1(a,b))=
\left\lbrace
\begin{array}{cl}
1 &((\overline{i},\overline{j})=(a-b,b)\mathrm{\;or\;}(b-a,a));\\
0 &(\mathrm{otherwise}).
\end{array}
\right.
$
\item $w=s$, $t=\theta_s(c)\in T^{sF}$: for every $c\in\mathbb{Z}/(q^2-1)\mathbb{Z}$, lift it to $\dot{c}\in\mathbb{Z}$ with $0\leq \dot{c}\leq q^2-2$ ($\dot{c}$ is uniquely determined by $c$), write $\dot{c}=(q+1)\dot{u}+\dot{v}$ where $\dot{u},\dot{v}\in\mathbb{Z}$ with $0\leq \dot{v}\leq q$ and $0\leq \dot{u}\leq q-2$ (division of $\dot{c}$ by $(q+1)$); set also $u=\dot{u}+(q-1)\mathbb{Z}\in\mathbb{Z}/(q-1)\mathbb{Z}$ and $v=\dot{v}+(q-1)\mathbb{Z}\in\mathbb{Z}/(q-1)\mathbb{Z}$; then:\\[2mm]
if $\dot{v}=0$: 
$\Phi_s(r_{i,j})(\theta_s(c))=
\left\lbrace
\begin{array}{cl}
1 & ((i,j)=(0,\dot{u}));\\
0 & (\mathrm{otherwise});
\end{array}
\right.$\\[2mm]
if $\dot{v}=1\mathrm{\;or\;}q$:
$\Phi_s(r_{i,j})(\theta_s(c))=
\left\lbrace
\begin{array}{cl}
1 & ((i,j)=(1,\dot{u}));\\
0 & (\mathrm{otherwise});
\end{array}
\right.$\\[2mm]
if $2\leq \dot{v}\leq q-1$:
$\Phi_s(r_{i,j})(\theta_s(c))=
\left\lbrace
\begin{array}{cl}
1 & ((i,\overline{j})=(1,u)\mathrm{\;or\;}(q+1-\dot{v},u+v-1);\\
0 & (\mathrm{otherwise}).
\end{array}
\right.$
\end{enumerate}
These calculations will lead to the desired equalities; we omit the details here. \qed

\end{context}

\begin{context}\label{E-GL2-lemma}
{\sc Proposition.}\;---
{\it
Let $G=\mathrm{GL}_2(\overline{\mathbb{F}_q})$. Let $T$ and $W=\{1,s\}$ be as in \S\,\ref{B-GL2-lemma}, so that the Curtis homomorphism of $G$ is 
\[
\mathrm{Cur}^G=(\mathrm{Cur}_{T}^G,\mathrm{Cur}_{T_s}^G):\overline{\mathbb{Q}}\mathsf{E}_G\longrightarrow\overline{\mathbb{Q}}T^F\times\overline{\mathbb{Q}}(T_w)^F.
\]
Let $Z$ be the center of $G$ (so $Z^F=T^F\cap (T_s)^F)$. Then 
\[
\mathrm{Cur}^G(\overline{\mathbb{Z}}[{\textstyle\frac{1}{p}}]\mathsf{E}_G)=\mathrm{Cur}^G(\overline{\mathbb{Q}}\mathsf{E}_G)\cap\{(f_1,f_s)\in\overline{\mathbb{Z}}[{\textstyle\frac{1}{p}}]T^F\times\overline{\mathbb{Z}}[{\textstyle\frac{1}{p}}](T_s)^F:f_1|_{Z^F}-f_s|_{Z^F}\in 2\overline{\mathbb{Z}}[{\textstyle\frac{1}{p}}]Z^F\}.\\[-7mm]
\]
}

{\it Sketch of proof.} As in the proof of \S\,\ref{B-GL2-lemma}, the idea is again to choose a $\overline{\mathbb{Z}}[\frac{1}{p}]$-linear basis of $\overline{\mathbb{Z}}[\frac{1}{p}]\mathsf{E}_G$, calculate their image under the Curtis homomorphism and then invert the resulting system. Let $B$ be the upper-triangular Borel subgroup of $G$, so that its unipotent radical $U$ is also upper-triangular. Let $\psi_0:\mathbb{F}_q\longrightarrow\overline{\mathbb{Z}}[\frac{1}{p}]^\times$ be a non-trivial linear character, which induces, by the canonical additive group isomorphism $U^F\simeq\mathbb{F}_q$, a nondegenerate linear character $\psi:U\longrightarrow\overline{\mathbb{Z}}[\frac{1}{p}]^\times$. By \S\,\ref{E-definition}, a $\overline{\mathbb{Z}}[\frac{1}{p}]$-linear basis of $\overline{\mathbb{Z}}[\frac{1}{p}]\mathsf{E}_G$ is 
\[
\{c_a:=e_\psi n_ae_\psi:a\in\mathbb{F}_q^\times\}\sqcup\{c_{a,b}':=qe_\psi n'_{a,b}e_\psi:a,b\in\mathbb{F}_q^\times\}
\]
where $n_a:=\left[\begin{array}{cc}a&0\\0&a\end{array}\right]$ and $n_{a,b}':=\left[\begin{array}{cc}0&b\\-a^{-1}&0\end{array}\right]$. With the aide of Curtis' formula (\S\,\ref{Curtis-definition}) and reasoning as in \cite[Sec.\;5]{Curtis}, one may obtain the following formulae: 
\begin{enumerate}[(i)]
\item $\mathrm{Cur}_T^G(c_a)(t)=
\left\lbrace
\begin{array}{cl}
1 & (t=n_a);\\
0 & (\mathrm{otherwise});
\end{array}
\right.$

$\mathrm{Cur}_T^G(c_{a,b}')(t)=
\left\lbrace
\begin{array}{cl}
\psi_0(a.\mathrm{tr}(t)) & (\det(t)=\det(n_{a,b}'));\\
0 & (\mathrm{otherwise});
\end{array}
\right.$
\item $\mathrm{Cur}_{T_s}^G(c_a)(t)=
\left\lbrace
\begin{array}{cl}
1 & (t=n_a);\\
0 & (\mathrm{otherwise}).
\end{array}
\right.$

$\mathrm{Cur}_{T_s}^G(c_{a,b}')(t)=
\left\lbrace
\begin{array}{cl}
-\psi_0(a.\mathrm{tr}(t)) & (\det(t)=\det(n_{a,b}'));\\
0 & (\mathrm{otherwise}).
\end{array}
\right.$
\end{enumerate}
Again, these formulae will lead to the desired equality, and we omit the details here. \qed

\end{context}

\begin{context}\label{GL2-exceptional}
{\sc Corollary} (of \S\S\,\ref{B-Curtis}\,-\,\ref{E-GL2-lemma}).\;--- 
{\it 
For $G=\mathrm{GL}_2(\overline{\mathbb{F}_q})$, the $\mathbb{Z}[\frac{1}{p|W|}]$-algebra isomorphisms $\mathbb{Z}[\frac{1}{p|W|}]\mathsf{E}_G\simeq\mathbb{Z}[\frac{1}{p|W|}]\mathsf{B}_{G^\vee}\simeq \mathbb{Z}[\frac{1}{p|W|}]\mathsf{K}_{G^\ast}$ in (\ref{identification-eq2}) induce (by restriction) ${\mathbb{Z}}[\frac{1}{p}]$-algebra isomorphisms ${\mathbb{Z}}[\frac{1}{p}]\mathsf{E}_G\simeq{\mathbb{Z}}[\frac{1}{p}]\mathsf{B}_{G^\vee}\simeq {\mathbb{Z}}[\frac{1}{p}]\mathsf{K}_{G^\ast}$.
}
\end{context}

\begin{context}\label{GL2-non-saturation}
{\it On the non-saturation of Curtis homomorphisms.}\;--- Still consider the case of $G=\mathrm{GL}_2(\overline{\mathbb{F}_q})$ (so that $|W|=2$) and suppose that $q$ is odd. In this case, $\mathrm{Cur}^G$ and $\Phi$ are both not saturated over $\overline{\mathbb{Z}}[{\textstyle\frac{1}{p}}]$ (with respect to $\overline{\mathbb{Q}}$): indeed, by the $\overline{\mathbb{Z}}[\frac{1}{p}]$-algebra isomorphism $\overline{\mathbb{Z}}[\frac{1}{p}]\mathsf{E}_G\simeq\overline{\mathbb{Z}}[\frac{1}{p}]\mathsf{B}_{G^\vee}$ in \S\,\ref{GL2-exceptional} and the $\mathsf{E}$-$\mathsf{B}$ commutative diagram in \S\,\ref{B-Curtis}, it suffices it suffices to show that $\Phi$ is not saturated over $\overline{\mathbb{Z}}[{\textstyle\frac{1}{p}}]$; using the notation in the proof of \S\,\ref{B-GL2-lemma}, if we consider the element $f=\sum_{i,j}f_{i,j}r_{i,j}\in\overline{\mathbb{Q}}\mathsf{B}_{G^\vee}$ (the sum is over $0\leq i\leq q-1$ and $0\leq j\leq q-2$) where $f_{i,j}=1/2$ if $i=2,4,6,\cdots,q-1$ and $f_{i,j}=0$ otherwise, then the calculations made in the proof of \S\,\ref{B-GL2-lemma} will show that
\[
f\not\in\overline{\mathbb{Z}}[{\textstyle\frac{1}{p}}]\mathsf{B}_{G^\vee}\mbox{\;\;but\;\;}\Phi(f)\in\prod_{w\in W}\overline{\mathbb{Z}}[{\textstyle\frac{1}{p}}]T^{wF};
\]
thus $\Phi$ is not saturated over $\overline{\mathbb{Z}}[{\textstyle\frac{1}{p}}]$. The same argument also shows that neither $\Phi$ nor $\mathrm{Cur}^G$ is saturated over ${\mathbb{Z}}[{\textstyle\frac{1}{p}}]$ (with respect to $\mathbb{Q}$).
\end{context}

\begin{context}\label{PGL2-exceptional}
{\it The example $G=\mathrm{PGL}_2(\overline{\mathbb{F}_q})$.}\;--- For $G=\mathrm{PGL}_2(\overline{\mathbb{F}_q})$, we know that its dual is $G^\ast=\mathrm{SL}_2(\overline{\mathbb{F}_q})$, which is already simply-connected. Using a similar discussion as in the previous $\mathrm{GL}_2$-case, it can be shown that the conclusions of \S\S\,\ref{B-GL2-lemma}\,-\,\ref{GL2-exceptional} all hold for $G=\mathrm{PGL}_2(\overline{\mathbb{F}_q})$; we omit the details here.
\end{context}

\subsection{On the structure of $\mathsf{B}_{G^\vee}$ for general $G^\vee$}

\begin{context}
{\it The problem.}\;--- When $G_{\mathrm{der}}^\vee$ is simply-connected, in \S\,\ref{q-restricted} we have shown that $\mathsf{B}_{G^\vee}$ is a reduced ring and is a free $\mathbb{Z}$-module with an explicit $\mathbb{Z}$-basis $\mathfrak{F}$. We would like to see what can we say about the $\mathbb{Z}$-module structure of $\mathsf{B}_{G^\vee}$ in the case where $G_{\mathrm{der}}^\vee$ is not simply-connected; we haven't obtained a general theory so far, and here we limit ourself to the case of $G=\mathrm{SO}_{2n}(\overline{\mathbb{F}_q})$ (with $n\geq 4$): in this case, $G^\vee$ and $G_{\mathrm{der}}^\vee$ are both isomorphic to $\mathrm{SO}_{2n}$, which is not simply-connected; nevertheless, we can still show that the $\mathsf{B}_{G^\vee}$ here is a reduced ring and is a free $\mathbb{Z}$-module with an explicit $\mathbb{Z}$-basis.
\end{context}

\begin{context}\label{SO2n-setup}
{\it The group $G=\mathrm{SO}_{2n}(\overline{\mathbb{F}_q})$.}\;--- Let $n\geq 4$ and $J=\left[\begin{array}{cccc}&&1\\&$\reflectbox{$\ddots$}$&\\1&&\end{array}\right]\in\mathrm{GL}_{2n}(\overline{\mathbb{F}_q})$. Consider the special orthogonal group $G=\mathrm{SO}_{2n}(\overline{\mathbb{F}_q})=\{A\in\mathrm{GL}_{2n}(\overline{\mathbb{F}_q})\,|\,{}^tA\cdot J\cdot A=J\}$, equipped with the Frobenius action $F:G\longrightarrow G$ which is the restriction of the standard Frobenius action $(g_{ij})_{i,j}\longmapsto (g_{ij}^q)_{i,j}$ on $\mathrm{GL}_{2n}(\overline{\mathbb{F}_q})$; then $(G,F)$ is split over $\mathbb{F}_q$. The diagonal matrices (resp.\;the upper-triangular matrices) in $G$ form an $F$-stable maximal torus $T$ (resp.\;an $F$-stable Borel subgroup $B$). The Langlands dual of $(G,T,B)$ is $(G^\vee,T^\vee,B^\vee)$ with $G^\vee=\mathrm{SO}_{2n}$, $T^\vee=\{\mathrm{diagonal\;elements\;in\;}G^\vee\}$ and $B^\vee=\{\mbox{upper-triangular\;elements\;in\;}G^\vee\}$. An explicit description of elements of $T^\vee$ is:
\[
T^\vee=\{\theta(z_1,\cdots,z_n):=\mathrm{diag}\,(z_1,\cdots,z_n,z_n^{-1},\cdots,z_1^{-1})\,|\,z_1,\cdots,z_n\in\mathbb{G}_m\}\simeq\mathbb{G}_m^n.
\]

For $1\leq i\leq n$, consider the following two morphisms:
\[
\left\lbrace
\begin{array}{ll}
\varepsilon_i:T^\vee\longrightarrow\mathbb{G}_m,&\theta(z_1,\cdots,z_n)\longmapsto z_i;\\
\eta_i:\mathbb{G}_m\longrightarrow T^\vee,&z\longmapsto \theta(1,\cdots,z,\cdots,1)\;\;(z\mbox{ at the $i$-th position}).
\end{array}
\right.
\]
Then $X(T^\vee)$ (resp.\;$Y(T^\vee)$) is a free $\mathbb{Z}$-module with basis $\{\varepsilon_1,\cdots,\varepsilon_n\}$ (resp.\;$\{\eta_1,\cdots,\eta_n\}$), and the canonical perfect pairing $\langle\cdot,\cdot\rangle:X(T^\vee)\times Y(T^\vee)\longrightarrow\mathbb{Z}$ verifies $\langle\varepsilon_i,\eta_j\rangle=\delta_{ij}$ (Kronecker delta) for all $i,j$.

The set of simple roots determined by $(G^\vee,T^\vee,B^\vee)$ is $\Delta^\vee=\{\alpha_1,\cdots,\alpha_n\}\subset X(T^\vee)$ where $\alpha_i=\varepsilon_i-\varepsilon_{i+1}$ ($1\leq i\leq n-1$) and $\alpha_n=\varepsilon_{n-1}+\varepsilon_n$. The root system of $G^\vee$ is of type $D_n$. The $\mathbb{Z}$-lattice generated by $\Delta^\vee$ has index 2 in $X(T^\vee)$, so $G^\vee$ is not semisimple and hence not simply-connected. We also have $G_{\mathrm{der}}^\vee=G^\vee$.

The Weyl group $W=N_{G^\vee}(T^\vee)/T^\vee$ of $(G^\vee,T^\vee)$ may be described as follows: let $\mathfrak{S}_{\{\pm\varepsilon_1,\cdots,\pm\varepsilon_n\}}$ be the symmetric group on the set $\{\pm\varepsilon_1,\cdots,\pm\varepsilon_n\}$, and consider the following two group inclusions:
\[
\begin{array}{ll}
i_1:\mathfrak{S}_n\hookrightarrow\mathfrak{S}_{\{\pm\varepsilon_1,\cdots,\pm\varepsilon_n\}},&\sigma\mapsto (\pm\varepsilon_i\mapsto \pm\varepsilon_{\sigma(i)}\;\;(1\leq i\leq n));\\
i_2:\{\pm1\}^n\hookrightarrow\mathfrak{S}_{\{\pm\varepsilon_1,\cdots,\pm\varepsilon_n\}},&(j_1,\cdots,j_n)\mapsto (\pm\varepsilon_i\mapsto\pm j_i\varepsilon_i\;\;(1\leq i\leq n));
\end{array}
\]
set also $D=\{(j_1,\cdots,j_n)\in\{\pm1\}^n\,|\,\mbox{$(-1)$ appears an even times in $j_1,\cdots,j_n$}\}$; then we may identify $W=i_1(\mathfrak{S}_n).i_2(D)\subset\mathfrak{S}_{\{\pm\varepsilon_1,\cdots,\pm\varepsilon_n\}}$ (it is an inner semidirect product), and we have $|W|=|\mathfrak{S}_n|\cdot|D|=n!\cdot 2^{n-1}$.
\end{context}

\begin{context}\label{SO2n-algebra}
{\sc Proposition.}\;---  
{\it 
Let $G=\mathrm{SO}_{2n}(\overline{\mathbb{F}_q})$ and keep the notations in \S\,\ref{SO2n-setup}. Let $\Sigma_i=r(\varepsilon_1+\cdots+\varepsilon_i)\in\mathbb{Z}[X(T^\vee)]^W$ for $i\in\{1,\cdots,n\}$ and $\Sigma_{n}'=r(\varepsilon_1+\cdots+\varepsilon_{n-1}-\varepsilon_n)\in\mathbb{Z}[X(T^\vee)]^W$. Then $\mathbb{Z}[X(T^\vee)]^W=\mathbb{Z}[\Sigma_1,\cdots,\Sigma_n,\Sigma_n']$.
}

{\it Sketch of proof.} The proof is similar to the reduction process in the proof of \S\,\ref{q-restricted}, but here we don't have the height function for induction and we need to introduce new quantities to play the role of the height function.

For all $\lambda\in X^+(T^\vee)$, either (i) $\lambda=a_1\varepsilon_1+\cdots+a_n\varepsilon_n$ where $(a_1,\cdots,a_n)\in\mathbb{Z}^n$ with $a_1\geq\cdots\geq a_n\geq 0$, or (ii) $\lambda=a_1\varepsilon_1+\cdots+a_{n-1}\varepsilon_{n-1}-a_n\varepsilon_n$ where $(a_1,\cdots,a_n)\in\mathbb{Z}^n$ with $a_1\geq\cdots\geq a_n>0$. In either case, we define $h(\lambda)=a_1\in\mathbb{N}$, $\ell(\lambda)=\#\{i\in\{1,\cdots,n\}:a_i=h(\lambda)\}\in\mathbb{N}^\ast$ and $d(\lambda)=\displaystyle\sum_{i\in\{\ell(\lambda)+1,\cdots,n\},\,a_i\geq 1}(a_i-1)\in\mathbb{N}$.

{\small
\begin{multicols}{2}
\setlength{\columnseprule}{0.5pt}
\begin{center}
\fbox{$\lambda=a_1\varepsilon_1+\cdots+a_m\varepsilon_m$} \\[2mm]
($a_1=\cdots=a_{\ell}>a_{\ell+1}\geq\cdots\geq a_m\geq 0$)
\end{center}
\[
\begin{tikzpicture}
\draw (5.5,0) -- (-0.5,0) -- (-0.5,3) -- (1.5,3) -- (1.5,2.5) -- (2,2.5) -- (2,2) -- (2.5,2) -- (2.5,1) -- (4,1) -- (4,0);
\draw[dashed] (0,0) -- (0,3); \draw[dashed] (1,0) -- (1,3); \draw[dashed] (1.5,0) -- (1.5,2.5); \draw[dashed] (2,0) -- (2,2); \draw[dashed] (3.5,0) -- (3.5,1);
\draw (2.8,0.5) node {$\cdots$}; \draw (0.55,1.5) node {$\cdots$};
\draw (2.8,-0.3) node {$\cdots$}; \draw (0.55,-0.3) node {$\cdots$}; \draw (4.5,-0.3) node {$\cdots$};
\draw (-0.25,0) node[below] {\small $\varepsilon_1$};
\draw (1.25,0) node[below] {\small $\varepsilon_\ell$};
\draw (1.75,0) node[below] {\small $\varepsilon_{\ell+1}$};
\draw (3.75,0) node[below] {\small $\varepsilon_m$};
\draw (5.25,0) node[below] {\small $\varepsilon_n$};
\draw (-0.25,1.5) node {\small $a_1$};
\draw (1.25,1.5) node {\small $a_\ell$};
\draw (1.85,1) node {\footnotesize $a_{\ell+1}$};
\draw (3.75,0.5) node {\small $a_m$};
\draw[<->] (-0.7,0) -- (-0.7,3); \draw[<->] (-0.5,3.2) -- (1.5,3.2);
\draw (-1.2,1.5) node {\small $h(\lambda)$}; 
\draw (0.5,3.5) node {\small $\ell(\lambda)$};
\end{tikzpicture}
\]
\vspace*{\fill}
\columnbreak
\begin{center}
\fbox{$\lambda=a_1\varepsilon_1+\cdots+a_{n-1}\varepsilon_{n-1}-a_n\varepsilon_n$} \\[2mm]
($a_1=\cdots=a_{\ell}>a_{\ell+1}\geq\cdots\geq a_n>0$)
\end{center}
\[
\begin{tikzpicture}
\draw (5.5,0) -- (-0.5,0) -- (-0.5,3) -- (1.5,3) -- (1.5,2.5) -- (2.5,2.5) -- (2.5,2) -- (3,2) -- (3,1) -- (5,1) -- (5,-1) -- (5.5,-1) -- (5.5,0);
\draw[dashed] (0,0) -- (0,3); \draw[dashed] (1,0) -- (1,3); \draw[dashed] (1.5,0) -- (1.5,2.5); \draw[dashed] (2,0) -- (2,2.5); \draw[dashed] (4.5,0) -- (4.5,1);
\draw (3.3,0.5) node {$\cdots$}; \draw (0.55,1.5) node {$\cdots$};
\draw (3.3,-0.3) node {$\cdots$}; \draw (0.55,-0.3) node {$\cdots$}; 
\draw (-0.25,0) node[below] {\small $\varepsilon_1$};
\draw (1.25,0) node[below] {\small $\varepsilon_\ell$};
\draw (1.75,0) node[below] {\small $\varepsilon_{\ell+1}$};
\draw (4.65,0) node[below] {\footnotesize $\varepsilon_{n-1}$};
\draw (5.25,0) node[above] {\small $\varepsilon_n$};
\draw (-0.25,1.5) node {\small $a_1$};
\draw (1.25,1.5) node {\small $a_\ell$};
\draw (1.85,1) node {\footnotesize $a_{\ell+1}$};
\draw (4.65,0.5) node {\footnotesize $a_{n-1}$};
\draw (5.25,-0.5) node {\small $a_n$};
\draw[<->] (-0.7,0) -- (-0.7,3); \draw[<->] (-0.5,3.2) -- (1.5,3.2);
\draw (-1.2,1.5) node {\small $h(\lambda)$}; 
\draw (0.5,3.5) node {\small $\ell(\lambda)$};
\end{tikzpicture}
\]
\end{multicols}
}

The goal is to prove that $r(\lambda)\in\mathbb{Z}[\Sigma_1,\cdots,\Sigma_n,\Sigma_n']$ for all $\lambda\in X^+(T^\vee)$ (recall from \S\,\ref{invariant-algebra} that these $r(\lambda)$ form a $\mathbb{Z}$-basis for $\mathbb{Z}[X(T^\vee)]^W$). The reduction proceeds as follows:
\begin{itemize}
\item The case $\ell(\lambda)=1$: if $\lambda$ is of type (i) as described above, expand the product $r(\lambda-\varepsilon_1)\cdot\Sigma_1$ to reduce $h(\lambda)$; if $\lambda$ is of type (ii), expand instead the product $r(\lambda-(\varepsilon_1+\cdots+\varepsilon_{n-1}-\varepsilon_n))\cdot\Sigma_n'$ to reduce $h(\lambda)$ (in case that $h(\lambda)$ is not reduced, $h(\lambda)$ remains the same while $d(\lambda)$ is reduced). 
\item The case $\ell(\lambda)>1$: if $\lambda$ is of type (i), expand $r(\lambda-(\varepsilon_1+\cdots+\varepsilon_{\ell(\lambda)}))\cdot\Sigma_{\ell(\lambda)}$ to reduce $h(\lambda)$ or $\ell(\lambda)$; if $\lambda$ is of type (ii), expand $r(\lambda-(\varepsilon_1+\cdots+\varepsilon_{n-1}-\varepsilon_n))\cdot\Sigma_n'$ to reduce $h(\lambda)$ or $\ell(\lambda)$ (in case that $h(\lambda)$ and $\ell(\lambda)$ are not reduced, they remain the same while $d(\lambda)$ is reduced). \qed
\end{itemize}

\end{context}

\begin{context}\label{SO2n-basis}
{\sc Proposition.}\;--- 
{\it 
Let $G=\mathrm{SO}_{2n}(\overline{\mathbb{F}_q})$ and keep the notations in \S\,\ref{SO2n-setup}, and let $I$ be the ideal of $\mathbb{Z}[X(T^\vee)]^W$ appearing in the definition $\mathsf{B}_{G^\vee}=\mathbb{Z}[X(T^\vee)]^W/{I}$. Then:
\begin{enumerate}[\normalfont (a)]
\item $|(T^\vee(\overline{\mathbb{Q}})\sslash W)^{F^\vee}|=|(T^\ast\sslash W)^{F^\ast}|=2q^n-2q^{n-1}+q^{n-2}$;
\item $\mathsf{B}_{G^\vee}=\mathbb{Z}[X(T^\vee)]^W/I$ is a reduced ring and is also a free $\mathbb{Z}$-module of rank $|(T^\vee(\overline{\mathbb{Q}})\sslash W)^{F^\vee}|$ having the set $S=S_1\sqcup S_2\sqcup S_2'$ as a $\mathbb{Z}$-basis, where:
\begin{itemize}
\item $S_1$ consists of those $r(a_1\varepsilon_1+\cdots+a_n\varepsilon_n)+I\in\mathsf{B}_{G^\vee}$ where $a_i\in\mathbb{N}$ $(1\leq i\leq n)$, $0\leq a_i-a_{i+1}<q$ $(1\leq i\leq n-1)$ and $0\leq a_n<q$;
\item $S_2$ consists of those $r(a_1\varepsilon_1+\cdots+a_{n-1}\varepsilon_{n-1}-a_n\varepsilon_n)+I\in\mathsf{B}_{G^\vee}$ where $a_i\in\mathbb{N}$ $(1\leq i\leq n)$, $0\leq a_i-a_{i+1}<q$ $(1\leq i\leq n-2)$ and $0< a_n\leq a_{n-1}<q$;
\item $S_2'$ consists of those $r(a_1\varepsilon_1+\cdots+a_{n-1}\varepsilon_{n-1}-a_n\varepsilon_n)+I\in\mathsf{B}_{G^\vee}$ where $a_i\in\mathbb{N}$ $(1\leq i\leq n)$, $0\leq a_i-a_{i+1}<q$ $(1\leq i\leq n-2)$ and $0<a_n< a_{n-1}-q<q$.
\end{itemize}
Observe that $|S_1|=q^n$, $|S_2|=q^{n-2}\cdot\binom{q-1}{2}$ and $|S_2'|=q^{n-2}\cdot(q-1+\binom{q-1}{2})$, so $|S|=|S_1|+|S_2|+|S_2'|=2q^n-2q^{n-1}+q^{n-2}$. 
\end{enumerate}
}

{\it Sketch of proof.} (a) The first equality has been seen in \S\,\ref{counting-lemma}, so it suffices to establish the second equality $|(T^\ast\sslash W)^{F^\ast}|=2q^n-2q^{n-1}+q^{n-2}$; to do this, consider the surjection
\[
\begin{array}{cccccc}
\pi:&(T^\ast\sslash W)^{F^\ast}&\twoheadrightarrow&(T_{\mathrm{GL}_n}^\ast\sslash\mathfrak{S}_n)^{F^\ast}&\xrightarrow{\;\sim\;}&\left\lbrace\begin{array}{cc}\mbox{\small monic one-variable}\\\mbox{\small polynomials over $\mathbb{F}_q$}\\\mbox{\small of degree $n$}\end{array}\right\rbrace\\&
\theta(z_1,\cdots,z_n)&\mapsto&\mathrm{diag}(z_1+z_1^{-1},\cdots,z_n+z_n^{-1})&\\
&&&z&\mapsto&\left(\begin{array}{cc}\mbox{\small monic characteristic}\\\mbox{\small polynomial of $z$}\end{array}\right)
\end{array}
\]
and calculate $|(T^\ast\sslash W)^{F^\ast}|$ as the sum of the size of each fiber of $\pi$.

(b) The proof is a refined version of the proof of \S\,\ref{SO2n-algebra} and uses again techniques from the proofs of \S\,\ref{q-restricted} (in particular one proves first that $S$ generates the $\mathbb{Z}$-module $\mathsf{B}_{G^\vee}$), but involves no new insights otherwise; we thus omit the details here. \qed

\end{context}

\end{document}